% File jytex.tex, for jyTeX version 2.6M (June 1992)
% Copyright (c) 1991, 1992 by Jonathan P. Yamron
% For full documentation, "get jydoc" from hep-ph@xxx.lanl.gov
%   Problems?  Contact brahm@theory3.caltech.edu.

\catcode`\@=11

%*****************************************************************************

\message{Loading jyTeX fonts...}

%************************************************************
%*
%*             Available fonts
%*
%************************************************************

%************** 5-point fonts *******************************

\font\vptrm=cmr5 \font\vptmit=cmmi5 \font\vptsy=cmsy5 \font\vptbf=cmbx5

\skewchar\vptmit='177 \skewchar\vptsy='60 \fontdimen16
\vptsy=\the\fontdimen17 \vptsy

\def\vpt{\ifmmode\err@badsizechange\else
     \@mathfontinit
     \textfont0=\vptrm  \scriptfont0=\vptrm  \scriptscriptfont0=\vptrm
     \textfont1=\vptmit \scriptfont1=\vptmit \scriptscriptfont1=\vptmit
     \textfont2=\vptsy  \scriptfont2=\vptsy  \scriptscriptfont2=\vptsy
     \textfont3=\xptex  \scriptfont3=\xptex  \scriptscriptfont3=\xptex
     \textfont\bffam=\vptbf
     \scriptfont\bffam=\vptbf
     \scriptscriptfont\bffam=\vptbf
     \@fontstyleinit
     \def\rm{\vptrm\fam=\z@}%
     \def\bf{\vptbf\fam=\bffam}%
     \def\oldstyle{\vptmit\fam=\@ne}%
     \rm\fi}

%************** 6-point fonts *******************************

\font\viptrm=cmr6 \font\viptmit=cmmi6 \font\viptsy=cmsy6
\font\viptbf=cmbx6

\skewchar\viptmit='177 \skewchar\viptsy='60 \fontdimen16
\viptsy=\the\fontdimen17 \viptsy

\def\vipt{\ifmmode\err@badsizechange\else
     \@mathfontinit
     \textfont0=\viptrm  \scriptfont0=\vptrm  \scriptscriptfont0=\vptrm
     \textfont1=\viptmit \scriptfont1=\vptmit \scriptscriptfont1=\vptmit
     \textfont2=\viptsy  \scriptfont2=\vptsy  \scriptscriptfont2=\vptsy
     \textfont3=\xptex   \scriptfont3=\xptex  \scriptscriptfont3=\xptex
     \textfont\bffam=\viptbf
     \scriptfont\bffam=\vptbf
     \scriptscriptfont\bffam=\vptbf
     \@fontstyleinit
     \def\rm{\viptrm\fam=\z@}%
     \def\bf{\viptbf\fam=\bffam}%
     \def\oldstyle{\viptmit\fam=\@ne}%
     \rm\fi}
%************** 7-point fonts *******************************

\font\viiptrm=cmr7 \font\viiptmit=cmmi7 \font\viiptsy=cmsy7
\font\viiptit=cmti7 \font\viiptbf=cmbx7

\skewchar\viiptmit='177 \skewchar\viiptsy='60 \fontdimen16
\viiptsy=\the\fontdimen17 \viiptsy

\def\viipt{\ifmmode\err@badsizechange\else
     \@mathfontinit
     \textfont0=\viiptrm  \scriptfont0=\vptrm  \scriptscriptfont0=\vptrm
     \textfont1=\viiptmit \scriptfont1=\vptmit \scriptscriptfont1=\vptmit
     \textfont2=\viiptsy  \scriptfont2=\vptsy  \scriptscriptfont2=\vptsy
     \textfont3=\xptex    \scriptfont3=\xptex  \scriptscriptfont3=\xptex
     \textfont\itfam=\viiptit
     \scriptfont\itfam=\viiptit
     \scriptscriptfont\itfam=\viiptit
     \textfont\bffam=\viiptbf
     \scriptfont\bffam=\vptbf
     \scriptscriptfont\bffam=\vptbf
     \@fontstyleinit
     \def\rm{\viiptrm\fam=\z@}%
     \def\it{\viiptit\fam=\itfam}%
     \def\bf{\viiptbf\fam=\bffam}%
     \def\oldstyle{\viiptmit\fam=\@ne}%
     \rm\fi}

%************** 8-point fonts *******************************

\font\viiiptrm=cmr8 \font\viiiptmit=cmmi8 \font\viiiptsy=cmsy8
\font\viiiptit=cmti8
%\font\viiiptsl=cmsl8
\font\viiiptbf=cmbx8
%\font\viiipttt=cmtt8
%\font\viiiptss=cmss8

\skewchar\viiiptmit='177 \skewchar\viiiptsy='60 \fontdimen16
\viiiptsy=\the\fontdimen17 \viiiptsy

\def\viiipt{\ifmmode\err@badsizechange\else
     \@mathfontinit
     \textfont0=\viiiptrm  \scriptfont0=\viptrm  \scriptscriptfont0=\vptrm
     \textfont1=\viiiptmit \scriptfont1=\viptmit \scriptscriptfont1=\vptmit
     \textfont2=\viiiptsy  \scriptfont2=\viptsy  \scriptscriptfont2=\vptsy
     \textfont3=\xptex     \scriptfont3=\xptex   \scriptscriptfont3=\xptex
     \textfont\itfam=\viiiptit
     \scriptfont\itfam=\viiptit
     \scriptscriptfont\itfam=\viiptit
     \textfont\bffam=\viiiptbf
     \scriptfont\bffam=\viptbf
     \scriptscriptfont\bffam=\vptbf
     \@fontstyleinit
     \def\rm{\viiiptrm\fam=\z@}%
     \def\it{\viiiptit\fam=\itfam}%
     \def\bf{\viiiptbf\fam=\bffam}%
     \def\oldstyle{\viiiptmit\fam=\@ne}%
     \rm\fi}

%************** Optional 9-point fonts **********************

\def\getixpt{%
     \font\ixptrm=cmr9
     \font\ixptmit=cmmi9
     \font\ixptsy=cmsy9
     \font\ixptit=cmti9
%     \font\ixptsl=cmsl9
     \font\ixptbf=cmbx9
%     \font\ixpttt=cmtt9
%     \font\ixptss=cmss9
     \skewchar\ixptmit='177 \skewchar\ixptsy='60
     \fontdimen16 \ixptsy=\the\fontdimen17 \ixptsy}

\def\ixpt{\ifmmode\err@badsizechange\else
     \@mathfontinit
     \textfont0=\ixptrm  \scriptfont0=\viiptrm  \scriptscriptfont0=\vptrm
     \textfont1=\ixptmit \scriptfont1=\viiptmit \scriptscriptfont1=\vptmit
     \textfont2=\ixptsy  \scriptfont2=\viiptsy  \scriptscriptfont2=\vptsy
     \textfont3=\xptex   \scriptfont3=\xptex    \scriptscriptfont3=\xptex
     \textfont\itfam=\ixptit
     \scriptfont\itfam=\viiptit
     \scriptscriptfont\itfam=\viiptit
     \textfont\bffam=\ixptbf
     \scriptfont\bffam=\viiptbf
     \scriptscriptfont\bffam=\vptbf
     \@fontstyleinit
     \def\rm{\ixptrm\fam=\z@}%
     \def\it{\ixptit\fam=\itfam}%
     \def\bf{\ixptbf\fam=\bffam}%
     \def\oldstyle{\ixptmit\fam=\@ne}%
     \rm\fi}

%************** 10-point fonts ******************************

\font\xptrm=cmr10 \font\xptmit=cmmi10 \font\xptsy=cmsy10
\font\xptex=cmex10 \font\xptit=cmti10 \font\xptsl=cmsl10
\font\xptbf=cmbx10 \font\xpttt=cmtt10 \font\xptss=cmss10
\font\xptsc=cmcsc10 \font\xptbfs=cmb10 \font\xptbmit=cmmib10

\skewchar\xptmit='177 \skewchar\xptbmit='177 \skewchar\xptsy='60
\fontdimen16 \xptsy=\the\fontdimen17 \xptsy

\def\xpt{\ifmmode\err@badsizechange\else
     \@mathfontinit
     \textfont0=\xptrm  \scriptfont0=\viiptrm  \scriptscriptfont0=\vptrm
     \textfont1=\xptmit \scriptfont1=\viiptmit \scriptscriptfont1=\vptmit
     \textfont2=\xptsy  \scriptfont2=\viiptsy  \scriptscriptfont2=\vptsy
     \textfont3=\xptex  \scriptfont3=\xptex    \scriptscriptfont3=\xptex
     \textfont\itfam=\xptit
     \scriptfont\itfam=\viiptit
     \scriptscriptfont\itfam=\viiptit
     \textfont\bffam=\xptbf
     \scriptfont\bffam=\viiptbf
     \scriptscriptfont\bffam=\vptbf
     \textfont\bfsfam=\xptbfs
     \scriptfont\bfsfam=\viiptbf
     \scriptscriptfont\bfsfam=\vptbf
     \textfont\bmitfam=\xptbmit
     \scriptfont\bmitfam=\viiptmit
     \scriptscriptfont\bmitfam=\vptmit
     \@fontstyleinit
     \def\rm{\xptrm\fam=\z@}%
     \def\it{\xptit\fam=\itfam}%
     \def\sl{\xptsl}%
     \def\bf{\xptbf\fam=\bffam}%
     \def\tt{\xpttt}%
     \def\ss{\xptss}%
     \def\sc{\xptsc}%
     \def\bfs{\xptbfs\fam=\bfsfam}%
     \def\bmit{\fam=\bmitfam}%
     \def\oldstyle{\xptmit\fam=\@ne}%
     \rm\fi}

%************** Optional 11-point fonts *********************

\def\getxipt{%
     \font\xiptrm=cmr10  scaled\magstephalf
     \font\xiptmit=cmmi10 scaled\magstephalf
     \font\xiptsy=cmsy10 scaled\magstephalf
     \font\xiptex=cmex10 scaled\magstephalf
     \font\xiptit=cmti10 scaled\magstephalf
     \font\xiptsl=cmsl10 scaled\magstephalf
     \font\xiptbf=cmbx10 scaled\magstephalf
     \font\xipttt=cmtt10 scaled\magstephalf
     \font\xiptss=cmss10 scaled\magstephalf
     \skewchar\xiptmit='177 \skewchar\xiptsy='60
     \fontdimen16 \xiptsy=\the\fontdimen17 \xiptsy}

\def\xipt{\ifmmode\err@badsizechange\else
     \@mathfontinit
     \textfont0=\xiptrm  \scriptfont0=\viiiptrm  \scriptscriptfont0=\viptrm
     \textfont1=\xiptmit \scriptfont1=\viiiptmit \scriptscriptfont1=\viptmit
     \textfont2=\xiptsy  \scriptfont2=\viiiptsy  \scriptscriptfont2=\viptsy
     \textfont3=\xiptex  \scriptfont3=\xptex     \scriptscriptfont3=\xptex
     \textfont\itfam=\xiptit
     \scriptfont\itfam=\viiiptit
     \scriptscriptfont\itfam=\viiptit
     \textfont\bffam=\xiptbf
     \scriptfont\bffam=\viiiptbf
     \scriptscriptfont\bffam=\viptbf
     \@fontstyleinit
     \def\rm{\xiptrm\fam=\z@}%
     \def\it{\xiptit\fam=\itfam}%
     \def\sl{\xiptsl}%
     \def\bf{\xiptbf\fam=\bffam}%
     \def\tt{\xipttt}%
     \def\ss{\xiptss}%
     \def\oldstyle{\xiptmit\fam=\@ne}%
     \rm\fi}

%************** 12-point fonts ******************************

\font\xiiptrm=cmr12 \font\xiiptmit=cmmi12 \font\xiiptsy=cmsy10
scaled\magstep1 \font\xiiptex=cmex10  scaled\magstep1 \font\xiiptit=cmti12
\font\xiiptsl=cmsl12 \font\xiiptbf=cmbx12
%\font\xiipttt=cmtt12
\font\xiiptss=cmss12 \font\xiiptsc=cmcsc10 scaled\magstep1
\font\xiiptbfs=cmb10  scaled\magstep1 \font\xiiptbmit=cmmib10
scaled\magstep1

\skewchar\xiiptmit='177 \skewchar\xiiptbmit='177 \skewchar\xiiptsy='60
\fontdimen16 \xiiptsy=\the\fontdimen17 \xiiptsy

\def\xiipt{\ifmmode\err@badsizechange\else
     \@mathfontinit
     \textfont0=\xiiptrm  \scriptfont0=\viiiptrm  \scriptscriptfont0=\viptrm
     \textfont1=\xiiptmit \scriptfont1=\viiiptmit \scriptscriptfont1=\viptmit
     \textfont2=\xiiptsy  \scriptfont2=\viiiptsy  \scriptscriptfont2=\viptsy
     \textfont3=\xiiptex  \scriptfont3=\xptex     \scriptscriptfont3=\xptex
     \textfont\itfam=\xiiptit
     \scriptfont\itfam=\viiiptit
     \scriptscriptfont\itfam=\viiptit
     \textfont\bffam=\xiiptbf
     \scriptfont\bffam=\viiiptbf
     \scriptscriptfont\bffam=\viptbf
     \textfont\bfsfam=\xiiptbfs
     \scriptfont\bfsfam=\viiiptbf
     \scriptscriptfont\bfsfam=\viptbf
     \textfont\bmitfam=\xiiptbmit
     \scriptfont\bmitfam=\viiiptmit
     \scriptscriptfont\bmitfam=\viptmit
     \@fontstyleinit
     \def\rm{\xiiptrm\fam=\z@}%
     \def\it{\xiiptit\fam=\itfam}%
     \def\sl{\xiiptsl}%
     \def\bf{\xiiptbf\fam=\bffam}%
     \def\tt{\xiipttt}%
     \def\ss{\xiiptss}%
     \def\sc{\xiiptsc}%
     \def\bfs{\xiiptbfs\fam=\bfsfam}%
     \def\bmit{\fam=\bmitfam}%
     \def\oldstyle{\xiiptmit\fam=\@ne}%
     \rm\fi}

%************** Optional 13-point fonts *********************

\def\getxiiipt{%
     \font\xiiiptrm=cmr12  scaled\magstephalf
     \font\xiiiptmit=cmmi12 scaled\magstephalf
     \font\xiiiptsy=cmsy9  scaled\magstep2
     \font\xiiiptit=cmti12 scaled\magstephalf
     \font\xiiiptsl=cmsl12 scaled\magstephalf
     \font\xiiiptbf=cmbx12 scaled\magstephalf
     \font\xiiipttt=cmtt12 scaled\magstephalf
     \font\xiiiptss=cmss12 scaled\magstephalf
     \skewchar\xiiiptmit='177 \skewchar\xiiiptsy='60
     \fontdimen16 \xiiiptsy=\the\fontdimen17 \xiiiptsy}

\def\xiiipt{\ifmmode\err@badsizechange\else
     \@mathfontinit
     \textfont0=\xiiiptrm  \scriptfont0=\xptrm  \scriptscriptfont0=\viiptrm
     \textfont1=\xiiiptmit \scriptfont1=\xptmit \scriptscriptfont1=\viiptmit
     \textfont2=\xiiiptsy  \scriptfont2=\xptsy  \scriptscriptfont2=\viiptsy
     \textfont3=\xivptex   \scriptfont3=\xptex  \scriptscriptfont3=\xptex
     \textfont\itfam=\xiiiptit
     \scriptfont\itfam=\xptit
     \scriptscriptfont\itfam=\viiptit
     \textfont\bffam=\xiiiptbf
     \scriptfont\bffam=\xptbf
     \scriptscriptfont\bffam=\viiptbf
     \@fontstyleinit
     \def\rm{\xiiiptrm\fam=\z@}%
     \def\it{\xiiiptit\fam=\itfam}%
     \def\sl{\xiiiptsl}%
     \def\bf{\xiiiptbf\fam=\bffam}%
     \def\tt{\xiiipttt}%
     \def\ss{\xiiiptss}%
     \def\oldstyle{\xiiiptmit\fam=\@ne}%
     \rm\fi}

%************** 14-point fonts ******************************

\font\xivptrm=cmr12   scaled\magstep1 \font\xivptmit=cmmi12
scaled\magstep1 \font\xivptsy=cmsy10  scaled\magstep2
\font\xivptex=cmex10  scaled\magstep2 \font\xivptit=cmti12 scaled\magstep1
\font\xivptsl=cmsl12  scaled\magstep1 \font\xivptbf=cmbx12
scaled\magstep1
%\font\xivpttt=cmtt12  scaled\magstep1
\font\xivptss=cmss12  scaled\magstep1 \font\xivptsc=cmcsc10
scaled\magstep2 \font\xivptbfs=cmb10  scaled\magstep2
\font\xivptbmit=cmmib10 scaled\magstep2

\skewchar\xivptmit='177 \skewchar\xivptbmit='177 \skewchar\xivptsy='60
\fontdimen16 \xivptsy=\the\fontdimen17 \xivptsy

\def\xivpt{\ifmmode\err@badsizechange\else
     \@mathfontinit
     \textfont0=\xivptrm  \scriptfont0=\xptrm  \scriptscriptfont0=\viiptrm
     \textfont1=\xivptmit \scriptfont1=\xptmit \scriptscriptfont1=\viiptmit
     \textfont2=\xivptsy  \scriptfont2=\xptsy  \scriptscriptfont2=\viiptsy
     \textfont3=\xivptex  \scriptfont3=\xptex  \scriptscriptfont3=\xptex
     \textfont\itfam=\xivptit
     \scriptfont\itfam=\xptit
     \scriptscriptfont\itfam=\viiptit
     \textfont\bffam=\xivptbf
     \scriptfont\bffam=\xptbf
     \scriptscriptfont\bffam=\viiptbf
     \textfont\bfsfam=\xivptbfs
     \scriptfont\bfsfam=\xptbfs
     \scriptscriptfont\bfsfam=\viiptbf
     \textfont\bmitfam=\xivptbmit
     \scriptfont\bmitfam=\xptbmit
     \scriptscriptfont\bmitfam=\viiptmit
     \@fontstyleinit
     \def\rm{\xivptrm\fam=\z@}%
     \def\it{\xivptit\fam=\itfam}%
     \def\sl{\xivptsl}%
     \def\bf{\xivptbf\fam=\bffam}%
     \def\tt{\xivpttt}%
     \def\ss{\xivptss}%
     \def\sc{\xivptsc}%
     \def\bfs{\xivptbfs\fam=\bfsfam}%
     \def\bmit{\fam=\bmitfam}%
     \def\oldstyle{\xivptmit\fam=\@ne}%
     \rm\fi}

%************** 17-point fonts ******************************

\font\xviiptrm=cmr17 \font\xviiptmit=cmmi12 scaled\magstep2
\font\xviiptsy=cmsy10 scaled\magstep3 \font\xviiptex=cmex10
scaled\magstep3 \font\xviiptit=cmti12 scaled\magstep2
\font\xviiptbf=cmbx12 scaled\magstep2 \font\xviiptbfs=cmb10
scaled\magstep3

\skewchar\xviiptmit='177 \skewchar\xviiptsy='60 \fontdimen16
\xviiptsy=\the\fontdimen17 \xviiptsy

\def\xviipt{\ifmmode\err@badsizechange\else
     \@mathfontinit
     \textfont0=\xviiptrm  \scriptfont0=\xiiptrm  \scriptscriptfont0=\viiiptrm
     \textfont1=\xviiptmit \scriptfont1=\xiiptmit \scriptscriptfont1=\viiiptmit
     \textfont2=\xviiptsy  \scriptfont2=\xiiptsy  \scriptscriptfont2=\viiiptsy
     \textfont3=\xviiptex  \scriptfont3=\xiiptex  \scriptscriptfont3=\xptex
     \textfont\itfam=\xviiptit
     \scriptfont\itfam=\xiiptit
     \scriptscriptfont\itfam=\viiiptit
     \textfont\bffam=\xviiptbf
     \scriptfont\bffam=\xiiptbf
     \scriptscriptfont\bffam=\viiiptbf
     \textfont\bfsfam=\xviiptbfs
     \scriptfont\bfsfam=\xiiptbfs
     \scriptscriptfont\bfsfam=\viiiptbf
     \@fontstyleinit
     \def\rm{\xviiptrm\fam=\z@}%
     \def\it{\xviiptit\fam=\itfam}%
     \def\bf{\xviiptbf\fam=\bffam}%
     \def\bfs{\xviiptbfs\fam=\bfsfam}%
     \def\oldstyle{\xviiptmit\fam=\@ne}%
     \rm\fi}

%************** 21-point fonts ******************************

\font\xxiptrm=cmr17  scaled\magstep1
%\font\xxiptmit=cmmi12 scaled\magstep3
%\font\xxiptsy=cmsy10 scaled\magstep4
%\font\xxiptex=cmex10 scaled\magstep4
%\font\xxiptbf=cmbx12 scaled\magstep3

%\skewchar\xxiptmit='177 \skewchar\xxiptsy='60
%\fontdimen16 \xxiptsy=\the\fontdimen17 \xxiptsy

\def\xxipt{\ifmmode\err@badsizechange\else
     \@mathfontinit
%     \textfont0=\xxiptrm  \scriptfont0=\xivptrm  \scriptscriptfont0=\xptrm
%     \textfont1=\xxiptmit \scriptfont1=\xivptmit \scriptscriptfont1=\xptmit
%     \textfont2=\xxiptsy  \scriptfont2=\xivptsy  \scriptscriptfont2=\xptsy
%     \textfont3=\xxiptex  \scriptfont3=\xivptex  \scriptscriptfont3=\xptex
%     \textfont\bffam=\xxiptbf
%     \scriptfont\bffam=\xivptbf
%     \scriptscriptfont\bffam=\xptbf
     \@fontstyleinit
     \def\rm{\xxiptrm\fam=\z@}%
     \rm\fi}

%************** 25-point fonts ******************************

\font\xxvptrm=cmr17  scaled\magstep2
%\font\xxvptmit=cmmi12 scaled\magstep4
%\font\xxvptsy=cmsy10 scaled\magstep5
%\font\xxvptex=cmex10 scaled\magstep5
%\font\xxvptbf=cmbx12 scaled\magstep4

%\skewchar\xxvptmit='177 \skewchar\xxvptsy='60
%\fontdimen16 \xxvptsy=\the\fontdimen17 \xxvptsy

\def\xxvpt{\ifmmode\err@badsizechange\else
     \@mathfontinit
%     \textfont0=\xxvptrm  \scriptfont0=\xviiptrm  \scriptscriptfont0=\xiiptrm
%     \textfont1=\xxvptmit \scriptfont1=\xviiptmit \scriptscriptfont1=\xiiptmit
%     \textfont2=\xxvptsy  \scriptfont2=\xviiptsy  \scriptscriptfont2=\xiiptsy
%     \textfont3=\xxvptex  \scriptfont3=\xviiptex  \scriptscriptfont3=\xiiptex
%     \textfont\bffam=\xxvptbf
%     \scriptfont\bffam=\xviiptbf
%     \scriptscriptfont\bffam=\xiiptbf
     \@fontstyleinit
     \def\rm{\xxvptrm\fam=\z@}%
     \rm\fi}

%************** Other fonts *********************************

%\font\dummy=dummy

%******************************************************************************

\message{Loading jyTeX macros...}

%************************************************************
%*
%*              Simple modifications to plain
%*
%************************************************************
\message{modifications to plain.tex,}

% The "\outer" qualifier is removed from the definitions of \newcount through
% \newif so that they may be used in definitions.  \newif is also changed to
% make \if commands globally defined.

\def\newcount{\alloc@0\count\countdef\insc@unt}
\def\newdimen{\alloc@1\dimen\dimendef\insc@unt}
\def\newskip{\alloc@2\skip\skipdef\insc@unt}
\def\newmuskip{\alloc@3\muskip\muskipdef\@cclvi}
\def\newbox{\alloc@4\box\chardef\insc@unt}
\def\newtoks{\alloc@5\toks\toksdef\@cclvi}
\def\newhelp#1#2{\newtoks#1\global#1\expandafter{\csname#2\endcsname}}
\def\newread{\alloc@6\read\chardef\sixt@@n}
\def\newwrite{\alloc@7\write\chardef\sixt@@n}
\def\newfam{\alloc@8\fam\chardef\sixt@@n}
\def\newinsert#1{\global\advance\insc@unt by\m@ne
     \ch@ck0\insc@unt\count
     \ch@ck1\insc@unt\dimen
     \ch@ck2\insc@unt\skip
     \ch@ck4\insc@unt\box
     \allocationnumber=\insc@unt
     \global\chardef#1=\allocationnumber
     \wlog{\string#1=\string\insert\the\allocationnumber}}
\def\newif#1{\count@\escapechar \escapechar\m@ne
     \expandafter\expandafter\expandafter
          \xdef\@if#1{true}{\let\noexpand#1=\noexpand\iftrue}%
     \expandafter\expandafter\expandafter
          \xdef\@if#1{false}{\let\noexpand#1=\noexpand\iffalse}%
     \global\@if#1{false}\escapechar=\count@}

%************** Some parameter changes **********************

\newlinechar=`\^^J
\overfullrule=0pt

%************** Font-related modifications ******************

% The plain fonts are mapped onto the corresponding jyTeX fonts

% Some control sequences are disabled.

\let\itfam=\undefined

\let\bffam=\undefined

\count18=3

% German sharp s is given a new name (\ss is already taken)

\chardef\sharps="19

% The mathcode assignments of characters in the math italic font are changed to
% allow for switching to boldface.

\mathchardef\alpha="710B \mathchardef\beta="710C
\mathchardef\gamma="710D \mathchardef\delta="710E
\mathchardef\epsilon="710F \mathchardef\zeta="7110
\mathchardef\eta="7111 \mathchardef\theta="7112 \mathchardef\iota="7113
\mathchardef\kappa="7114 \mathchardef\lambda="7115
\mathchardef\mu="7116 \mathchardef\nu="7117 \mathchardef\xi="7118
\mathchardef\pi="7119 \mathchardef\rho="711A \mathchardef\sigma="711B
\mathchardef\tau="711C \mathchardef\upsilon="711D
\mathchardef\phi="711E \mathchardef\chi="711F \mathchardef\psi="7120
\mathchardef\omega="7121 \mathchardef\varepsilon="7122
\mathchardef\vartheta="7123 \mathchardef\varpi="7124
\mathchardef\varrho="7125 \mathchardef\varsigma="7126
\mathchardef\varphi="7127 \mathchardef\imath="717B
\mathchardef\jmath="717C \mathchardef\ell="7160 \mathchardef\wp="717D
\mathchardef\partial="7140 \mathchardef\flat="715B
\mathchardef\natural="715C \mathchardef\sharp="715D

%************** Miscellaneous changes ***********************

% The dimension \p@ (1pt) is replaced with \rp@ (relative pt, defined below),
% whose size is determined by the base type size of the document.

\def\angle{{\vbox{\ialign{$\m@th\scriptstyle##$\crcr
     \not\mathrel{\mkern14mu}\crcr
     \noalign{\nointerlineskip}
     \mkern2.5mu\leaders\hrule height.34\rp@\hfill\mkern2.5mu\crcr}}}}
\def\vdots{\vbox{\baselineskip4\rp@ \lineskiplimit\z@
     \kern6\rp@\hbox{.}\hbox{.}\hbox{.}}}
\def\ddots{\mathinner{\mkern1mu\raise7\rp@\vbox{\kern7\rp@\hbox{.}}\mkern2mu
     \raise4\rp@\hbox{.}\mkern2mu\raise\rp@\hbox{.}\mkern1mu}}
\def\overrightarrow#1{\vbox{\ialign{##\crcr
     \rightarrowfill\crcr
     \noalign{\kern-\rp@\nointerlineskip}
     $\hfil\displaystyle{#1}\hfil$\crcr}}}
\def\overleftarrow#1{\vbox{\ialign{##\crcr
     \leftarrowfill\crcr
     \noalign{\kern-\rp@\nointerlineskip}
     $\hfil\displaystyle{#1}\hfil$\crcr}}}
\def\overbrace#1{\mathop{\vbox{\ialign{##\crcr
     \noalign{\kern3\rp@}
     \downbracefill\crcr
     \noalign{\kern3\rp@\nointerlineskip}
     $\hfil\displaystyle{#1}\hfil$\crcr}}}\limits}
\def\underbrace#1{\mathop{\vtop{\ialign{##\crcr
     $\hfil\displaystyle{#1}\hfil$\crcr
     \noalign{\kern3\rp@\nointerlineskip}
     \upbracefill\crcr
     \noalign{\kern3\rp@}}}}\limits}
\def\big#1{{\hbox{$\left#1\vbox to8.5\rp@ {}\right.\n@space$}}}
\def\Big#1{{\hbox{$\left#1\vbox to11.5\rp@ {}\right.\n@space$}}}
\def\bigg#1{{\hbox{$\left#1\vbox to14.5\rp@ {}\right.\n@space$}}}
\def\Bigg#1{{\hbox{$\left#1\vbox to17.5\rp@ {}\right.\n@space$}}}
\def\@vereq#1#2{\lower.5\rp@\vbox{\baselineskip\z@skip\lineskip-.5\rp@
     \ialign{$\m@th#1\hfil##\hfil$\crcr#2\crcr=\crcr}}}
\def\rlh@#1{\vcenter{\hbox{\ooalign{\raise2\rp@
     \hbox{$#1\rightharpoonup$}\crcr
     $#1\leftharpoondown$}}}}
\def\bordermatrix#1{\begingroup\m@th
     \setbox\z@\vbox{%
          \def\cr{\crcr\noalign{\kern2\rp@\global\let\cr\endline}}%
          \ialign{$##$\hfil\kern2\rp@\kern\p@renwd
               &\thinspace\hfil$##$\hfil&&\quad\hfil$##$\hfil\crcr
               \omit\strut\hfil\crcr
               \noalign{\kern-\baselineskip}%
               #1\crcr\omit\strut\cr}}%
     \setbox\tw@\vbox{\unvcopy\z@\global\setbox\@ne\lastbox}%
     \setbox\tw@\hbox{\unhbox\@ne\unskip\global\setbox\@ne\lastbox}%
     \setbox\tw@\hbox{$\kern\wd\@ne\kern-\p@renwd\left(\kern-\wd\@ne
          \global\setbox\@ne\vbox{\box\@ne\kern2\rp@}%
          \vcenter{\kern-\ht\@ne\unvbox\z@\kern-\baselineskip}%
          \,\right)$}%
     \null\;\vbox{\kern\ht\@ne\box\tw@}\endgroup}
\def\endinsert{\egroup
     \if@mid\dimen@\ht\z@
          \advance\dimen@\dp\z@
          \advance\dimen@12\rp@
          \advance\dimen@\pagetotal
          \ifdim\dimen@>\pagegoal\@midfalse\p@gefalse\fi
     \fi
     \if@mid\bigskip\box\z@
          \bigbreak
     \else\insert\topins{\penalty100 \splittopskip\z@skip
               \splitmaxdepth\maxdimen\floatingpenalty\z@
               \ifp@ge\dimen@\dp\z@
                    \vbox to\vsize{\unvbox\z@\kern-\dimen@}%
               \else\box\z@\nobreak\bigskip
               \fi}%
     \fi
     \endgroup}

% \normalbaselines is removed from \cases and \matrix.

\def\cases#1{\left\{\,\vcenter{\m@th
     \ialign{$##\hfil$&\quad##\hfil\crcr#1\crcr}}\right.}
\def\matrix#1{\null\,\vcenter{\m@th
     \ialign{\hfil$##$\hfil&&\quad\hfil$##$\hfil\crcr
          \mathstrut\crcr
          \noalign{\kern-\baselineskip}
          #1\crcr
          \mathstrut\crcr
          \noalign{\kern-\baselineskip}}}\,}

% \raggedbottom modified slightly

\newif\ifraggedbottom

\def\raggedbottom{\ifraggedbottom\else
     \advance\topskip by\z@ plus60pt \raggedbottomtrue\fi}%
\def\normalbottom{\ifraggedbottom
     \advance\topskip by\z@ plus-60pt \raggedbottomfalse\fi}

%************************************************************
%*
%*              Miscellaneous definitions
%*
%************************************************************
\message{hacks,}

%************** Hack registers ******************************

\toksdef\toks@i=1 \toksdef\toks@ii=2

%************** Basic macros ********************************

\def\TeX{T\kern-.1667em \lower.5ex \hbox{E}\kern-.125em X\null}
\def\jyTeX{{\leavevmode
     \raise.587ex \hbox{\it\j}\kern-.1em \lower.048ex \hbox{\it y}\kern-.12em
     \TeX}}

\let\then=\iftrue
\def\ifnoarg#1\then{\def\hack@{#1}\ifx\hack@\empty}
\def\ifundefined#1\then{%
     \expandafter\ifx\csname\expandafter\blank\string#1\endcsname\relax}
\def\useif#1\then{\csname#1\endcsname}
\def\usename#1{\csname#1\endcsname}
\def\useafter#1#2{\expandafter#1\csname#2\endcsname}

% Modify so that I can have \loop's within \loop's?
\long\def\loop#1\repeat{\def\@iterate{#1\expandafter\@iterate\fi}\@iterate
     \let\@iterate=\relax}
%\long\def\loop#1\repeat{\def\@loopbody{#1}\@iterate}
%\def\@iterate{\@loopbody\let\next=\@iterate\else\let\next=\relax\fi\next}

\let\TeXend=\end
\def\begin#1{\begingroup\def\@@blockname{#1}\usename{begin#1}}
\def\end#1{\usename{end#1}\def\hack@{#1}%
     \ifx\@@blockname\hack@
          \endgroup
     \else\err@badgroup\hack@\@@blockname
     \fi}
\def\@@blockname{}

\def\defaultoption[#1]#2{%
     \def\hack@{\ifx\hack@ii[\toks@={#2}\else\toks@={#2[#1]}\fi\the\toks@}%
     \futurelet\hack@ii\hack@}

\def\markup#1{\let\@@marksf=\empty
     \ifhmode\edef\@@marksf{\spacefactor=\the\spacefactor\relax}\/\fi
     ${}^{\hbox{\subscriptfonts#1}}$\@@marksf}

%************** Time registers ******************************

\newtoks\shortyear
\newtoks\militaryhour
\newtoks\standardhour
\newtoks\minute
\newtoks\amorpm

\def\settime{\count@=\time\divide\count@ by60
     \militaryhour=\expandafter{\number\count@}%
     {\multiply\count@ by-60 \advance\count@ by\time
          \xdef\hack@{\ifnum\count@<10 0\fi\number\count@}}%
     \minute=\expandafter{\hack@}%
     \ifnum\count@<12
          \amorpm={am}
     \else\amorpm={pm}
          \ifnum\count@>12 \advance\count@ by-12 \fi
     \fi
     \standardhour=\expandafter{\number\count@}%
     \def\hack@19##1##2{\shortyear={##1##2}}%
          \expandafter\hack@\the\year}

\def\monthword#1{%
     \ifcase#1
          $\bullet$\err@badcountervalue{monthword}%
          \or January\or February\or March\or April\or May\or June%
          \or July\or August\or September\or October\or November\or December%
     \else$\bullet$\err@badcountervalue{monthword}%
     \fi}

\def\monthabbr#1{%
     \ifcase#1
          $\bullet$\err@badcountervalue{monthabbr}%
          \or Jan\or Feb\or Mar\or Apr\or May\or Jun%
          \or Jul\or Aug\or Sep\or Oct\or Nov\or Dec%
     \else$\bullet$\err@badcountervalue{monthabbr}%
     \fi}

\def\militarytime{\the\militaryhour:\the\minute}
\def\standardtime{\the\standardhour:\the\minute}

%************** Number styles *******************************

\def\@setnumstyle#1#2{\expandafter\global\expandafter\expandafter
     \expandafter\let\expandafter\expandafter
     \csname @\expandafter\blank\string#1style\endcsname
     \csname#2\endcsname}
\def\numstyle#1{\usename{@\expandafter\blank\string#1style}#1}
\def\ifblank#1\then{\useafter\ifx{@\expandafter\blank\string#1}\blank}

\def\blank#1{}

\def\Roman#1{\expandafter\uppercase\expandafter{\romannumeral#1}}
\def\alphabetic#1{%
     \ifcase#1
          $\bullet$\err@badcountervalue{alphabetic}%
          \or a\or b\or c\or d\or e\or f\or g\or h\or i\or j\or k\or l\or m%
          \or n\or o\or p\or q\or r\or s\or t\or u\or v\or w\or x\or y\or z%
     \else$\bullet$\err@badcountervalue{alphabetic}%
     \fi}
\def\Alphabetic#1{\expandafter\uppercase\expandafter{\alphabetic{#1}}}
\def\symbols#1{%
     \ifcase#1
          $\bullet$\err@badcountervalue{symbols}%
          \or*\or\dag\or\ddag\or\S\or$\|$%
          \or**\or\dag\dag\or\ddag\ddag\or\S\S\or$\|\|$%
     \else$\bullet$\err@badcountervalue{symbols}%
     \fi}

%************** String macros *******************************

\catcode`\^^?=13 \def^^?{\relax}

\def\trimleading#1\to#2{\edef#2{#1}%
     \expandafter\@trimleading\expandafter#2#2^^?^^?}
\def\@trimleading#1#2#3^^?{\ifx#2^^?\def#1{}\else\def#1{#2#3}\fi}

\def\trimtrailing#1\to#2{\edef#2{#1}%
     \expandafter\@trimtrailing\expandafter#2#2^^? ^^?\relax}
\def\@trimtrailing#1#2 ^^?#3{\ifx#3\relax\toks@={}%
     \else\def#1{#2}\toks@={\trimtrailing#1\to#1}\fi
     \the\toks@}

\def\trim#1\to#2{\trimleading#1\to#2\trimtrailing#2\to#2}

\catcode`\^^?=15

%************** List macros *********************************

\long\def\additemL#1\to#2{\toks@={\^^\{#1}}\toks@ii=\expandafter{#2}%
     \xdef#2{\the\toks@\the\toks@ii}}

\long\def\additemR#1\to#2{\toks@={\^^\{#1}}\toks@ii=\expandafter{#2}%
     \xdef#2{\the\toks@ii\the\toks@}}

\def\getitemL#1\to#2{\expandafter\@getitemL#1\hack@#1#2}
\def\@getitemL\^^\#1#2\hack@#3#4{\def#4{#1}\def#3{#2}}

%************************************************************
%*
%*             Font-related macros
%*
%************************************************************
\message{font macros,}

%************** Font set-up *********************************

\newdimen\rp@
\newcount\@@sizeindex \@@sizeindex=0
\newcount\@@factori
\newcount\@@factorii
\newcount\@@factoriii
\newcount\@@factoriv

\countdef\maxfam=18
\newfam\itfam
\newfam\bffam
\newfam\bfsfam
\newfam\bmitfam

\def\@mathfontinit{\count@=4
     \loop\textfont\count@=\nullfont
          \scriptfont\count@=\nullfont
          \scriptscriptfont\count@=\nullfont
          \ifnum\count@<\maxfam\advance\count@ by\@ne
     \repeat}

\def\@fontstyleinit{%
     \def\it{\err@fontnotavailable\it}%
     \def\bf{\err@fontnotavailable\bf}%
     \def\bfs{\err@bfstobf}%
     \def\bmit{\err@fontnotavailable\bmit}%
     \def\sc{\err@fontnotavailable\sc}%
     \def\sl{\err@sltoit}%
     \def\ss{\err@fontnotavailable\ss}%
     \def\tt{\err@fontnotavailable\tt}}

\def\@parameterinit#1{\rm\rp@=.1em \@getscaling{#1}%
     \let\^^\=\@doscaling\scalingskipslist
     \setbox\strutbox=\hbox{\vrule
          height.708\baselineskip depth.292\baselineskip width\z@}}

\def\@getfactor#1#2#3#4{\@@factori=#1 \@@factorii=#2
     \@@factoriii=#3 \@@factoriv=#4}

\def\@getscaling#1{\count@=#1 \advance\count@ by-\@@sizeindex\@@sizeindex=#1
     \ifnum\count@<0
          \let\@mulordiv=\divide
          \let\@divormul=\multiply
          \multiply\count@ by\m@ne
     \else\let\@mulordiv=\multiply
          \let\@divormul=\divide
     \fi
     \edef\@@scratcha{\ifcase\count@                {1}{1}{1}{1}\or
          {1}{7}{23}{3}\or     {2}{5}{3}{1}\or      {9}{89}{13}{1}\or
          {6}{25}{6}{1}\or     {8}{71}{14}{1}\or    {6}{25}{36}{5}\or
          {1}{7}{53}{4}\or     {12}{125}{108}{5}\or {3}{14}{53}{5}\or
          {6}{41}{17}{1}\or    {13}{31}{13}{2}\or   {9}{107}{71}{2}\or
          {11}{139}{124}{3}\or {1}{6}{43}{2}\or     {10}{107}{42}{1}\or
          {1}{5}{43}{2}\or     {5}{69}{65}{1}\or    {11}{97}{91}{2}\fi}%
     \expandafter\@getfactor\@@scratcha}

\def\@doscaling#1{\@mulordiv#1by\@@factori\@divormul#1by\@@factorii
     \@mulordiv#1by\@@factoriii\@divormul#1by\@@factoriv}

%************* Size-changing commands ***********************

\newskip\headskip
\newskip\footskip

\def\typesize=#1pt{\count@=#1 \advance\count@ by-10
     \ifcase\count@
          \@setsizex\or\err@badtypesize\or
          \@setsizexii\or\err@badtypesize\or
          \@setsizexiv
     \else\err@badtypesize
     \fi}

\def\@setsizex{\getixpt
     \def\subsubscriptfonts{\vpt}%
          \def\subsubscriptsize{\vpt\@parameterinit{-8}}%
     \def\subscriptfonts{\viipt}\def\subscriptsize{\viipt\@parameterinit{-4}}%
     \def\footnotefonts{\viiipt}\def\footnotesize{\viiipt\@parameterinit{-2}}%
     \def\smallfonts{\ixpt}\def\smallsize{\ixpt\@parameterinit{-1}}%
     \def\normalfonts{\xpt}\def\normalsize{\xpt\@parameterinit{0}}%
     \def\bigfonts{\xiipt}\def\bigsize{\xiipt\@parameterinit{2}}%
     \def\Bigfonts{\xivpt}\def\Bigsize{\xivpt\@parameterinit{4}}%
     \def\biggfonts{\xviipt}\def\biggsize{\xviipt\@parameterinit{6}}%
     \def\Biggfonts{\xxipt}\def\Biggsize{\xxipt\@parameterinit{8}}%
     \def\tinyfonts{\vpt}\def\tinysize{\vpt\@parameterinit{-8}}%
     \def\HUGEFONTS{\xxvpt}\def\HUGESIZE{\xxvpt\@parameterinit{10}}%
     \normalsize\fixedskipslist}

\def\@setsizexii{\getxipt
     \def\subsubscriptfonts{\vipt}%
          \def\subsubscriptsize{\vipt\@parameterinit{-6}}%
     \def\subscriptfonts{\viiipt}%
          \def\subscriptsize{\viiipt\@parameterinit{-2}}%
     \def\footnotefonts{\xpt}\def\footnotesize{\xpt\@parameterinit{0}}%
     \def\smallfonts{\xipt}\def\smallsize{\xipt\@parameterinit{1}}%
     \def\normalfonts{\xiipt}\def\normalsize{\xiipt\@parameterinit{2}}%
     \def\bigfonts{\xivpt}\def\bigsize{\xivpt\@parameterinit{4}}%
     \def\Bigfonts{\xviipt}\def\Bigsize{\xviipt\@parameterinit{6}}%
     \def\biggfonts{\xxipt}\def\biggsize{\xxipt\@parameterinit{8}}%
     \def\Biggfonts{\xxvpt}\def\Biggsize{\xxvpt\@parameterinit{10}}%
     \def\tinyfonts{\vpt}\def\tinysize{\vpt\@parameterinit{-8}}%
     \def\HUGEFONTS{\xxvpt}\def\HUGESIZE{\xxvpt\@parameterinit{10}}%
     \normalsize\fixedskipslist}

\def\@setsizexiv{\getxiiipt
     \def\subsubscriptfonts{\viipt}%
          \def\subsubscriptsize{\viipt\@parameterinit{-4}}%
     \def\subscriptfonts{\xpt}\def\subscriptsize{\xpt\@parameterinit{0}}%
     \def\footnotefonts{\xiipt}\def\footnotesize{\xiipt\@parameterinit{2}}%
     \def\smallfonts{\xiiipt}\def\smallsize{\xiiipt\@parameterinit{3}}%
     \def\normalfonts{\xivpt}\def\normalsize{\xivpt\@parameterinit{4}}%
     \def\bigfonts{\xviipt}\def\bigsize{\xviipt\@parameterinit{6}}%
     \def\Bigfonts{\xxipt}\def\Bigsize{\xxipt\@parameterinit{8}}%
     \def\biggfonts{\xxvpt}\def\biggsize{\xxvpt\@parameterinit{10}}%
     \def\Biggfonts{\err@sizetoolarge\Biggfonts\HUGEFONTS}%
          \def\Biggsize{\err@sizetoolarge\Biggsize\HUGESIZE}%
     \def\tinyfonts{\vpt}\def\tinysize{\vpt\@parameterinit{-8}}%
     \def\HUGEFONTS{\xxvpt}\def\HUGESIZE{\xxvpt\@parameterinit{10}}%
     \normalsize\fixedskipslist}

\def\subsubscriptfonts{\vpt} \def\subsubscriptsize{\vpt\@parameterinit{-8}}
\def\subscriptfonts{\viipt}  \def\subscriptsize{\viipt\@parameterinit{-4}}
\def\footnotefonts{\viiipt}  \def\footnotesize{\viiipt\@parameterinit{-2}}
\def\smallfonts{\err@sizenotavailable\smallfonts}
                             \def\smallsize{\ixpt\@parameterinit{-1}}
\def\normalfonts{\xpt}       \def\normalsize{\xpt\@parameterinit{0}}
\def\bigfonts{\xiipt}        \def\bigsize{\xiipt\@parameterinit{2}}
\def\Bigfonts{\xivpt}        \def\Bigsize{\xivpt\@parameterinit{4}}
\def\biggfonts{\xviipt}      \def\biggsize{\xviipt\@parameterinit{6}}
\def\Biggfonts{\xxipt}       \def\Biggsize{\xxipt\@parameterinit{8}}
\def\tinyfonts{\vpt}         \def\tinysize{\vpt\@parameterinit{-8}}
\def\HUGEFONTS{\xxvpt}       \def\HUGESIZE{\xxvpt\@parameterinit{10}}

%************************************************************
%*
%*             Document layout
%*
%************************************************************
\message{document layout,}

%************** Page format *********************************

\newtoks\everyoutput \everyoutput={}
\newdimen\depthofpage
\newcount\pagenum \pagenum=0

\newdimen\oddtopmargin  \newdimen\eventopmargin
\newdimen\oddleftmargin \newdimen\evenleftmargin
\newtoks\oddhead        \newtoks\evenhead
\newtoks\oddfoot        \newtoks\evenfoot

\def\topmargin{\afterassignment\@seteventop\oddtopmargin}
\def\leftmargin{\afterassignment\@setevenleft\oddleftmargin}
\def\head{\afterassignment\@setevenhead\oddhead}
\def\foot{\afterassignment\@setevenfoot\oddfoot}

\def\@seteventop{\eventopmargin=\oddtopmargin}
\def\@setevenleft{\evenleftmargin=\oddleftmargin}
\def\@setevenhead{\evenhead=\oddhead}
\def\@setevenfoot{\evenfoot=\oddfoot}

\def\pagenumstyle#1{\@setnumstyle\pagenum{#1}}

\newif\ifdraft
\def\draft{\drafttrue\leftmargin=.5in \overfullrule=5pt }

\def\outputstyle#1{\global\expandafter\let\expandafter
          \@outputstyle\csname#1output\endcsname
     \usename{#1setup}}

\output={\@outputstyle}

\def\normaloutput{\the\everyoutput
     \global\advance\pagenum by\@ne
     \ifodd\pagenum
          \voffset=\oddtopmargin \hoffset=\oddleftmargin
     \else\voffset=\eventopmargin \hoffset=\evenleftmargin
     \fi
     \advance\voffset by-1in  \advance\hoffset by-1in
     \count0=\pagenum
     \expandafter\shipout\pagebox
     \ifnum\outputpenalty>-\@MM\else\dosupereject\fi}

\newdimen\fullhsize
\newbox\leftpage
\newcount\leftpagenum
\newcount\outputpagenum \outputpagenum=0
\let\leftorright=L

\def\twoupoutput{\the\everyoutput
     \global\advance\pagenum by\@ne
     \if L\leftorright
          \global\setbox\leftpage=\leftline{\pagebox}%
          \global\leftpagenum=\pagenum
          \global\let\leftorright=R%
     \else\global\advance\outputpagenum by\@ne
          \ifodd\outputpagenum
               \voffset=\oddtopmargin \hoffset=\oddleftmargin
          \else\voffset=\eventopmargin \hoffset=\evenleftmargin
          \fi
          \advance\voffset by-1in  \advance\hoffset by-1in
          \count0=\leftpagenum \count1=\pagenum
          \shipout\vbox{\hbox to\fullhsize
               {\box\leftpage\hfil\leftline{\pagebox}}}%
          \global\let\leftorright=L%
     \fi
     \ifnum\outputpenalty>-\@MM
     \else\dosupereject
          \if R\leftorright
               \globaldefs=\@ne\head={\hfil}\foot={\hfil}\globaldefs=\z@
               \null\newpage
          \fi
     \fi}

\def\pagebox{\vbox{\makeheadline\pagebody\makefootline}}

\def\makeheadline{%
     \vbox to\z@{\baselinestretch=\@m
          \vskip\topskip\vskip-.708\baselineskip\vskip-\headskip
          \line{\vbox to\ht\strutbox{}%
               \ifodd\pagenum\the\oddhead\else\the\evenhead\fi}%
          \vss}%
     \nointerlineskip}

\def\pagebody{\vbox to\vsize{%
     \boxmaxdepth\maxdepth
     \ifvoid\topins\else\unvbox\topins\fi
     \depthofpage=\dp255
     \unvbox255
     \ifraggedbottom\kern-\depthofpage\vfil\fi
     \ifvoid\footins
     \else\vskip\skip\footins
          \footnoterule
          \unvbox\footins
          \vskip-\footnoteskip
     \fi}}

\def\makefootline{\baselineskip=\footskip
     \line{\ifodd\pagenum\the\oddfoot\else\the\evenfoot\fi}}

%************** Sectioning commands *************************

\newskip\abovechapterskip
\newskip\belowchapterskip
\newskip\abovesectionskip
\newskip\belowsectionskip
\newskip\abovesubsectionskip
\newskip\belowsubsectionskip

\def\chapterstyle#1{\global\expandafter\let\expandafter\@chapterstyle
     \csname#1text\endcsname}
\def\sectionstyle#1{\global\expandafter\let\expandafter\@sectionstyle
     \csname#1text\endcsname}
\def\subsectionstyle#1{\global\expandafter\let\expandafter\@subsectionstyle
     \csname#1text\endcsname}

\def\chapter#1{%
     \ifdim\lastskip=17sp \else\chapterbreak\vskip\abovechapterskip\fi
     \@chapterstyle{\ifblank\chapternumstyle\then
          \else\newchapternum=\next\chapternumformat\ \fi#1}%
     \nobreak\vskip\belowchapterskip\vskip17sp }

\def\section#1{%
     \ifdim\lastskip=17sp \else\sectionbreak\vskip\abovesectionskip\fi
     \@sectionstyle{\ifblank\sectionnumstyle\then
          \else\newsectionnum=\next\sectionnumformat\ \fi#1}%
     \nobreak\vskip\belowsectionskip\vskip17sp }

\def\subsection#1{%
     \ifdim\lastskip=17sp \else\subsectionbreak\vskip\abovesubsectionskip\fi
     \@subsectionstyle{\ifblank\subsectionnumstyle\then
          \else\newsubsectionnum=\next\subsectionnumformat\ \fi#1}%
     \nobreak\vskip\belowsubsectionskip\vskip17sp }

%************** Text formatting commands ********************

\let\TeXunderline=\underline
\let\TeXoverline=\overline
\def\underline#1{\relax\ifmmode\TeXunderline{#1}\else
     $\TeXunderline{\hbox{#1}}$\fi}
\def\overline#1{\relax\ifmmode\TeXoverline{#1}\else
     $\TeXoverline{\hbox{#1}}$\fi}

\def\baselinestretch{\afterassignment\@baselinestretch\count@}
\def\@baselinestretch{\baselineskip=\normalbaselineskip
     \divide\baselineskip by\@m\baselineskip=\count@\baselineskip
     \setbox\strutbox=\hbox{\vrule
          height.708\baselineskip depth.292\baselineskip width\z@}%
     \bigskipamount=\the\baselineskip
          plus.25\baselineskip minus.25\baselineskip
     \medskipamount=.5\baselineskip
          plus.125\baselineskip minus.125\baselineskip
     \smallskipamount=.25\baselineskip
          plus.0625\baselineskip minus.0625\baselineskip}

\def\\{\ifhmode\ifnum\lastpenalty=-\@M\else\hfil\penalty-\@M\fi\fi
     \ignorespaces}
\def\newpage{\vfil\break}

\def\lefttext#1{\par{\@text\leftskip=\z@\rightskip=\centering
     \noindent#1\par}}
\def\righttext#1{\par{\@text\leftskip=\centering\rightskip=\z@
     \noindent#1\par}}
\def\centertext#1{\par{\@text\leftskip=\centering\rightskip=\centering
     \noindent#1\par}}
\def\@text{\parindent=\z@ \parfillskip=\z@ \everypar={}%
     \spaceskip=.3333em \xspaceskip=.5em
     \def\\{\ifhmode\ifnum\lastpenalty=-\@M\else\penalty-\@M\fi\fi
          \ignorespaces}}

\def\beginleft{\par\@text\leftskip=\z@ \rightskip=\centering}
     
\def\beginright{\par\@text\leftskip=\centering\rightskip=\z@ }
     
\def\begincenter{\par\@text\leftskip=\centering\rightskip=\centering}

\def\beginnarrow{\defaultoption[\parindent]\@beginnarrow}
\def\@beginnarrow[#1]{\par\advance\leftskip by#1\advance\rightskip by#1}

\begingroup
\catcode`\[=1 \catcode`\{=11 \gdef\beginignore[\endgroup\bgroup
     \catcode`\e=0 \catcode`\\=12 \catcode`\{=11 \catcode`\f=12 \let\or=\relax
     \let\nd{ignor=\fi \let\}=\egroup
     \iffalse}
\endgroup

\long\def\marginnote#1{\leavevmode
     \edef\@marginsf{\spacefactor=\the\spacefactor\relax}%
     \ifdraft\strut\vadjust{%
          \hbox to\z@{\hskip\hsize\hskip.1in
               \vbox to\z@{\vskip-\dp\strutbox
                    \marginnoteformat
                    \vskip-\ht\strutbox
                    \noindent\strut#1\par
                    \vss}%
               \hss}}%
     \fi
     \@marginsf}

%************** The \bye command ****************************

\newtoks\everybye \everybye={\par\vfil}
\outer\def\bye{\the\everybye
     \footnotecheck
     \prelabelcheck
     \streamcheck
     \supereject
     \TeXend}

%************************************************************
%*
%*             Footnotes
%*
%************************************************************
\message{footnotes,}

\newcount\footnotenum \footnotenum=0
\newskip\footnoteskip
\let\@footnotelist=\empty

\def\footnotenumstyle#1{\@setnumstyle\footnotenum{#1}%
     \useafter\ifx{@footnotenumstyle}\symbols
          \global\let\@footup=\empty
     \else\global\let\@footup=\markup
     \fi}

\def\footnote{\footnotecheck\defaultoption[]\@footnote}
\def\@footnote[#1]{\@footnotemark[#1]\@footnotetext}

\def\footnotemark{\defaultoption[]\@footnotemark}
\def\@footnotemark[#1]{\let\@footsf=\empty
     \ifhmode\edef\@footsf{\spacefactor=\the\spacefactor\relax}\/\fi
     \ifnoarg#1\then
          \global\advance\footnotenum by\@ne
          \@footup{\footnotenumformat}%
          \edef\@@foota{\footnotenum=\the\footnotenum\relax}%
          \expandafter\additemR\expandafter\@footup\expandafter
               {\@@foota\footnotenumformat}\to\@footnotelist
          \global\let\@footnotelist=\@footnotelist
     \else\markup{#1}%
          \additemR\markup{#1}\to\@footnotelist
          \global\let\@footnotelist=\@footnotelist
     \fi
     \@footsf}

\def\footnotetext{%
     \ifx\@footnotelist\empty\err@extrafootnotetext\else\@footnotetext\fi}
\def\@footnotetext{%
     \getitemL\@footnotelist\to\@@foota
     \global\let\@footnotelist=\@footnotelist
     \insert\footins\bgroup
     \footnoteformat
     \splittopskip=\ht\strutbox\splitmaxdepth=\dp\strutbox
     \interlinepenalty=\interfootnotelinepenalty\floatingpenalty=\@MM
     \noindent\llap{\@@foota}\strut
     \bgroup\aftergroup\@footnoteend
     \let\@@scratcha=}
\def\@footnoteend{\strut\par\vskip\footnoteskip\egroup}

\def\footnoterule{\normalfonts
     \kern-.3em \hrule width2in height.04em \kern .26em }

\def\footnotecheck{%
     \ifx\@footnotelist\empty
     \else\err@extrafootnotemark
          \global\let\@footnotelist=\empty
     \fi}

%************************************************************
%*
%*             Labelling macros
%*
%************************************************************
\message{labels,}

\let\@@labeldef=\xdef
\newif\if@labelfile
\newwrite\@labelfile
\let\@prelabellist=\empty

\def\label#1#2{\trim#1\to\@@labarg\edef\@@labtext{#2}%
     \edef\@@labname{lab@\@@labarg}%
     \useafter\ifundefined\@@labname\then\else\@yeslab\fi
     \useafter\@@labeldef\@@labname{#2}%
     \ifstreaming
          \expandafter\toks@\expandafter\expandafter\expandafter
               {\csname\@@labname\endcsname}%
          \immediate\write\streamout{\noexpand\label{\@@labarg}{\the\toks@}}%
     \fi}
\def\@yeslab{%
     \useafter\ifundefined{if\@@labname}\then
          \err@labelredef\@@labarg
     \else\useif{if\@@labname}\then
               \err@labelredef\@@labarg
          \else\global\usename{\@@labname true}%
               \useafter\ifundefined{pre\@@labname}\then
               \else\useafter\ifx{pre\@@labname}\@@labtext
                    \else\err@badlabelmatch\@@labarg
                    \fi
               \fi
               \if@labelfile
               \else\global\@labelfiletrue
                    \immediate\write\sixt@@n{--> Creating file \jobname.lab}%
                    \immediate\openout\@labelfile=\jobname.lab
               \fi
               \immediate\write\@labelfile
                    {\noexpand\prelabel{\@@labarg}{\@@labtext}}%
          \fi
     \fi}

\def\putlab#1{\trim#1\to\@@labarg\edef\@@labname{lab@\@@labarg}%
     \useafter\ifundefined\@@labname\then\@nolab\else\usename\@@labname\fi}
\def\@nolab{%
     \useafter\ifundefined{pre\@@labname}\then
          \undefinedlabelformat
          \err@needlabel\@@labarg
          \useafter\xdef\@@labname{\undefinedlabelformat}%
     \else\usename{pre\@@labname}%
          \useafter\xdef\@@labname{\usename{pre\@@labname}}%
     \fi
     \useafter\newif{if\@@labname}%
     \expandafter\additemR\@@labarg\to\@prelabellist}

\def\prelabel#1{\useafter\gdef{prelab@#1}}

\def\ifundefinedlabel#1\then{%
     \expandafter\ifx\csname lab@#1\endcsname\relax}
\def\useiflab#1\then{\csname iflab@#1\endcsname}

\def\prelabelcheck{{%
     \def\^^\##1{\useiflab{##1}\then\else\err@undefinedlabel{##1}\fi}%
     \@prelabellist}}

%************************************************************
%*
%*             Equation numbering
%*
%************************************************************
\message{equation numbering,}

\newcount\chapternum
\newcount\sectionnum
\newcount\subsectionnum
\newcount\equationnum
\newcount\subequationnum
\newcount\figurenum
\newcount\subfigurenum
\newcount\tablenum
\newcount\subtablenum

\newif\if@subeqncount
\newif\if@subfigcount
\newif\if@subtblcount

\def\newchapternum{\newsectionnum=\z@\@resetnum\chapternum}
\def\newsectionnum{\newsubsectionnum=\z@\@resetnum\sectionnum}
\def\newsubsectionnum{\newequationnum=\z@\newfigurenum=\z@\newtablenum=\z@
     \@resetnum\subsectionnum}
\def\newequationnum{\newsubequationnum=\z@\@resetnum\equationnum}
\def\newsubequationnum{\@resetnum\subequationnum}
\def\newfigurenum{\newsubfigurenum=\z@\@resetnum\figurenum}
\def\newsubfigurenum{\@resetnum\subfigurenum}
\def\newtablenum{\newsubtablenum=\z@\@resetnum\tablenum}
\def\newsubtablenum{\@resetnum\subtablenum}

\def\@resetnum#1{\global\advance#1by1 \edef\next{\the#1\relax}\global#1}

\newchapternum=0

\def\chapternumstyle#1{\@setnumstyle\chapternum{#1}}
\def\sectionnumstyle#1{\@setnumstyle\sectionnum{#1}}
\def\subsectionnumstyle#1{\@setnumstyle\subsectionnum{#1}}
\def\equationnumstyle#1{\@setnumstyle\equationnum{#1}}
\def\subequationnumstyle#1{\@setnumstyle\subequationnum{#1}%
     \ifblank\subequationnumstyle\then\global\@subeqncountfalse\fi
     \ignorespaces}
\def\figurenumstyle#1{\@setnumstyle\figurenum{#1}}
\def\subfigurenumstyle#1{\@setnumstyle\subfigurenum{#1}%
     \ifblank\subfigurenumstyle\then\global\@subfigcountfalse\fi
     \ignorespaces}
\def\tablenumstyle#1{\@setnumstyle\tablenum{#1}}
\def\subtablenumstyle#1{\@setnumstyle\subtablenum{#1}%
     \ifblank\subtablenumstyle\then\global\@subtblcountfalse\fi
     \ignorespaces}

\def\eqnlabel#1{%
     \if@subeqncount
          \newsubequationnum=\next
     \else\newequationnum=\next
          \ifblank\subequationnumstyle\then
          \else\global\@subeqncounttrue
               \newsubequationnum=\@ne
          \fi
     \fi
     \label{#1}{\puteqnformat}(\puteqn{#1})%
     \ifdraft\rlap{\hskip.1in{\tt#1}}\fi}

\let\puteqn=\putlab

\def\equation#1#2{\useafter\gdef{eqn@#1}{#2\eqno\eqnlabel{#1}}}
\def\Equation#1{\useafter\gdef{eqn@#1}}

\def\putequation#1{\useafter\ifundefined{eqn@#1}\then
     \err@undefinedeqn{#1}\else\usename{eqn@#1}\fi}

\def\eqnseriesstyle#1{\gdef\@eqnseriesstyle{#1}}
\def\begineqnseries{\subequationnumstyle{\@eqnseriesstyle}%
     \defaultoption[]\@begineqnseries}
\def\@begineqnseries[#1]{\edef\@@eqnname{#1}}
\def\endeqnseries{\subequationnumstyle{blank}%
     \expandafter\ifnoarg\@@eqnname\then
     \else\label\@@eqnname{\puteqnformat}%
     \fi
     \aftergroup\ignorespaces}

\def\figlabel#1{%
     \if@subfigcount
          \newsubfigurenum=\next
     \else\newfigurenum=\next
          \ifblank\subfigurenumstyle\then
          \else\global\@subfigcounttrue
               \newsubfigurenum=\@ne
          \fi
     \fi
     \label{#1}{\putfigformat}\putfig{#1}%
     {\def\marginnoteformat{\tt}\marginnote{#1}}}

\let\putfig=\putlab

\def\figseriesstyle#1{\gdef\@figseriesstyle{#1}}
\def\beginfigseries{\subfigurenumstyle{\@figseriesstyle}%
     \defaultoption[]\@beginfigseries}
\def\@beginfigseries[#1]{\edef\@@figname{#1}}
\def\endfigseries{\subfigurenumstyle{blank}%
     \expandafter\ifnoarg\@@figname\then
     \else\label\@@figname{\putfigformat}%
     \fi
     \aftergroup\ignorespaces}

\def\tbllabel#1{%
     \if@subtblcount
          \newsubtablenum=\next
     \else\newtablenum=\next
          \ifblank\subtablenumstyle\then
          \else\global\@subtblcounttrue
               \newsubtablenum=\@ne
          \fi
     \fi
     \label{#1}{\puttblformat}\puttbl{#1}%
     {\def\marginnoteformat{\tt}\marginnote{#1}}}

\let\puttbl=\putlab

\def\tblseriesstyle#1{\gdef\@tblseriesstyle{#1}}
\def\begintblseries{\subtablenumstyle{\@tblseriesstyle}%
     \defaultoption[]\@begintblseries}
\def\@begintblseries[#1]{\edef\@@tblname{#1}}
\def\endtblseries{\subtablenumstyle{blank}%
     \expandafter\ifnoarg\@@tblname\then
     \else\label\@@tblname{\puttblformat}%
     \fi
     \aftergroup\ignorespaces}

%************************************************************
%*
%*             Reference numbering
%*
%************************************************************
\message{reference numbering,}

\newcount\referencenum \referencenum=0
\newcount\@@prerefcount \@@prerefcount=0
\newcount\@@thisref
\newcount\@@lastref
\newcount\@@loopref
\newcount\@@refseq
\newdimen\refnumindent
\let\@undefreflist=\empty

\def\referencenumstyle#1{\@setnumstyle\referencenum{#1}}

\def\referencestyle#1{\usename{@ref#1}}

\def\@refsequential{%
     \gdef\@refpredef##1{\global\advance\referencenum by\@ne
          \let\^^\=0\label{##1}{\^^\{\the\referencenum}}%
          \useafter\gdef{ref@\the\referencenum}{{##1}{\undefinedlabelformat}}}%
     \gdef\@reference##1##2{%
          \ifundefinedlabel##1\then
          \else\def\^^\####1{\global\@@thisref=####1\relax}\putlab{##1}%
               \useafter\gdef{ref@\the\@@thisref}{{##1}{##2}}%
          \fi}%
     \gdef\endputreferences{%
          \loop\ifnum\@@loopref<\referencenum
                    \advance\@@loopref by\@ne
                    \expandafter\expandafter\expandafter\@printreference
                         \csname ref@\the\@@loopref\endcsname
          \repeat
          \par}}

\def\@refpreordered{%
     \gdef\@refpredef##1{\global\advance\referencenum by\@ne
          \additemR##1\to\@undefreflist}%
     \gdef\@reference##1##2{%
          \ifundefinedlabel##1\then
          \else\global\advance\@@loopref by\@ne
               {\let\^^\=0\label{##1}{\^^\{\the\@@loopref}}}%
               \@printreference{##1}{##2}%
          \fi}
     \gdef\endputreferences{%
          \def\^^\####1{\useiflab{####1}\then
               \else\reference{####1}{\undefinedlabelformat}\fi}%
          \@undefreflist
          \par}}

\def\beginprereferences{\par
     \def\reference##1##2{\global\advance\referencenum by1\@ne
          \let\^^\=0\label{##1}{\^^\{\the\referencenum}}%
          \useafter\gdef{ref@\the\referencenum}{{##1}{##2}}}}
\def\endprereferences{\global\@@prerefcount=\the\referencenum\par}

\def\beginputreferences{\par
     \refnumindent=\z@\@@loopref=\z@
     \loop\ifnum\@@loopref<\referencenum
               \advance\@@loopref by\@ne
               \setbox\z@=\hbox{\referencenum=\@@loopref
                    \referencenumformat\enskip}%
               \ifdim\wd\z@>\refnumindent\refnumindent=\wd\z@\fi
     \repeat
     \putreferenceformat
     \@@loopref=\z@
     \loop\ifnum\@@loopref<\@@prerefcount
               \advance\@@loopref by\@ne
               \expandafter\expandafter\expandafter\@printreference
                    \csname ref@\the\@@loopref\endcsname
     \repeat
     \let\reference=\@reference}

\def\@printreference#1#2{\ifx#2\undefinedlabelformat\err@undefinedref{#1}\fi
     \noindent\ifdraft\rlap{\hskip\hsize\hskip.1in \tt#1}\fi
     \llap{\referencenum=\@@loopref\referencenumformat\enskip}#2\par}

\def\reference#1#2{{\par\refnumindent=\z@\putreferenceformat\noindent#2\par}}

\def\putref#1{\trim#1\to\@@refarg
     \expandafter\ifnoarg\@@refarg\then
          \toks@={\relax}%
     \else\@@lastref=-\@m\def\@@refsep{}\def\@more{\@nextref}%
          \toks@={\@nextref#1,,}%
     \fi\the\toks@}
\def\@nextref#1,{\trim#1\to\@@refarg
     \expandafter\ifnoarg\@@refarg\then
          \let\@more=\relax
     \else\ifundefinedlabel\@@refarg\then
               \expandafter\@refpredef\expandafter{\@@refarg}%
          \fi
          \def\^^\##1{\global\@@thisref=##1\relax}%
          \global\@@thisref=\m@ne
          \setbox\z@=\hbox{\putlab\@@refarg}%
     \fi
     \advance\@@lastref by\@ne
     \ifnum\@@lastref=\@@thisref\advance\@@refseq by\@ne\else\@@refseq=\@ne\fi
     \ifnum\@@lastref<\z@
     \else\ifnum\@@refseq<\thr@@
               \@@refsep\def\@@refsep{,}%
               \ifnum\@@lastref>\z@
                    \advance\@@lastref by\m@ne
                    {\referencenum=\@@lastref\putrefformat}%
               \else\undefinedlabelformat
               \fi
          \else\def\@@refsep{--}%
          \fi
     \fi
     \@@lastref=\@@thisref
     \@more}

%************************************************************
%*
%*             Job streaming
%*
%************************************************************
\message{streaming,}

\newif\ifstreaming

\def\streamto{\defaultoption[\jobname]\@streamto}
\def\@streamto[#1]{\global\streamingtrue
     \immediate\write\sixt@@n{--> Streaming to #1.str}%
     \newwrite\streamout\immediate\openout\streamout=#1.str }

\def\streamfrom{\defaultoption[\jobname]\@streamfrom}
\def\@streamfrom[#1]{\newread\streamin\openin\streamin=#1.str
     \ifeof\streamin
          \expandafter\err@nostream\expandafter{#1.str}%
     \else\immediate\write\sixt@@n{--> Streaming from #1.str}%
          \let\@@labeldef=\gdef
          \ifstreaming
               \edef\@elc{\endlinechar=\the\endlinechar}%
               \endlinechar=\m@ne
               \loop\read\streamin to\@@scratcha
                    \ifeof\streamin
                         \streamingfalse
                    \else\toks@=\expandafter{\@@scratcha}%
                         \immediate\write\streamout{\the\toks@}%
                    \fi
                    \ifstreaming
               \repeat
               \@elc
               \input #1.str
               \streamingtrue
          \else\input #1.str
          \fi
          \let\@@labeldef=\xdef
     \fi}

\def\streamcheck{\ifstreaming
     \immediate\write\streamout{\pagenum=\the\pagenum}%
     \immediate\write\streamout{\footnotenum=\the\footnotenum}%
     \immediate\write\streamout{\referencenum=\the\referencenum}%
     \immediate\write\streamout{\chapternum=\the\chapternum}%
     \immediate\write\streamout{\sectionnum=\the\sectionnum}%
     \immediate\write\streamout{\subsectionnum=\the\subsectionnum}%
     \immediate\write\streamout{\equationnum=\the\equationnum}%
     \immediate\write\streamout{\subequationnum=\the\subequationnum}%
     \immediate\write\streamout{\figurenum=\the\figurenum}%
     \immediate\write\streamout{\subfigurenum=\the\subfigurenum}%
     \immediate\write\streamout{\tablenum=\the\tablenum}%
     \immediate\write\streamout{\subtablenum=\the\subtablenum}%
     \immediate\closeout\streamout
     \fi}

%************************************************************
%*
%*             Error messages
%*
%************************************************************

\def\err@badtypesize{%
     \errhelp={The limited availability of certain fonts requires^^J%
          that the base type size be 10pt, 12pt, or 14pt.^^J}%
     \errmessage{--> Illegal base type size}}

\def\err@badsizechange{\immediate\write\sixt@@n
     {--> Size change not allowed in math mode, ignored}}

\def\err@sizetoolarge#1{\immediate\write\sixt@@n
     {--> \noexpand#1 too big, substituting HUGE}}

\def\err@sizenotavailable#1{\immediate\write\sixt@@n
     {--> Size not available, \noexpand#1 ignored}}

\def\err@fontnotavailable#1{\immediate\write\sixt@@n
     {--> Font not available, \noexpand#1 ignored}}

\def\err@sltoit{\immediate\write\sixt@@n
     {--> Style \noexpand\sl not available, substituting \noexpand\it}%
     \it}

\def\err@bfstobf{\immediate\write\sixt@@n
     {--> Style \noexpand\bfs not available, substituting \noexpand\bf}%
     \bf}

\def\err@badgroup#1#2{%
     \errhelp={The block you have just tried to close was not the one^^J%
          most recently opened.^^J}%
     \errmessage{--> \noexpand\end{#1} doesn't match \noexpand\begin{#2}}}

\def\err@badcountervalue#1{\immediate\write\sixt@@n
     {--> Counter (#1) out of bounds}}

\def\err@extrafootnotemark{\immediate\write\sixt@@n
     {--> \noexpand\footnotemark command
          has no corresponding \noexpand\footnotetext}}

\def\err@extrafootnotetext{%
     \errhelp{You have given a \noexpand\footnotetext command without first
          specifying^^Ja \noexpand\footnotemark.^^J}%
     \errmessage{--> \noexpand\footnotetext command has no corresponding
          \noexpand\footnotemark}}

\def\err@labelredef#1{\immediate\write\sixt@@n
     {--> Label "#1" redefined}}

\def\err@badlabelmatch#1{\immediate\write\sixt@@n
     {--> Definition of label "#1" doesn't match value in \jobname.lab}}

\def\err@needlabel#1{\immediate\write\sixt@@n
     {--> Label "#1" cited before its definition}}

\def\err@undefinedlabel#1{\immediate\write\sixt@@n
     {--> Label "#1" cited but never defined}}

\def\err@undefinedeqn#1{\immediate\write\sixt@@n
     {--> Equation "#1" not defined}}

\def\err@undefinedref#1{\immediate\write\sixt@@n
     {--> Reference "#1" not defined}}

\def\err@nostream#1{%
     \errhelp={You have tried to input a stream file that doesn't exist.^^J}%
     \errmessage{--> Stream file #1 not found}}

%************************************************************
%*
%*             Initialization
%*
%************************************************************
\message{jyTeX initialization}

\everyjob{\immediate\write16{--> jyTeX version \fmtversion}%
     \edef\@@jobname{\jobname}%
%     \openin0=\inputpath jysupp
%     \ifeof0
%     \else\closein0
%          \immediate\write16{--> Additional macros loaded from jysupp.tex}%
%          \jyinput jysupp
%     \fi
%     \openin0=\inputpath jylocal
%     \ifeof0
%     \else\closein0
%          \immediate\write16{--> Additional macros loaded from jylocal.tex}%
%          \jyinput jylocal
%     \fi
     \edef\jobname{\@@jobname}%
     \settime
     \openin0=\jobname.lab
     \ifeof0
     \else\closein0
          \immediate\write16{--> Getting labels from file \jobname.lab}%
          \input\jobname.lab
     \fi}

%************** Spacing *************************************

\def\fixedskipslist{%
     \^^\{\topskip}%
     \^^\{\splittopskip}%
     \^^\{\maxdepth}%
     \^^\{\skip\topins}%
     \^^\{\skip\footins}%
     \^^\{\headskip}%
     \^^\{\footskip}}

\def\scalingskipslist{%
     \^^\{\p@renwd}%
     \^^\{\delimitershortfall}%
     \^^\{\nulldelimiterspace}%
     \^^\{\scriptspace}%
     \^^\{\jot}%
     \^^\{\normalbaselineskip}%
     \^^\{\normallineskip}%
     \^^\{\normallineskiplimit}%
     \^^\{\baselineskip}%
     \^^\{\lineskip}%
     \^^\{\lineskiplimit}%
     \^^\{\bigskipamount}%
     \^^\{\medskipamount}%
     \^^\{\smallskipamount}%
     \^^\{\parskip}%
     \^^\{\parindent}%
     \^^\{\abovedisplayskip}%
     \^^\{\belowdisplayskip}%
     \^^\{\abovedisplayshortskip}%
     \^^\{\belowdisplayshortskip}%
     \^^\{\abovechapterskip}%
     \^^\{\belowchapterskip}%
     \^^\{\abovesectionskip}%
     \^^\{\belowsectionskip}%
     \^^\{\abovesubsectionskip}%
     \^^\{\belowsubsectionskip}}

%************** Document layout *****************************

\def\twoupsetup{%                                % setup for twoup style
     \topmargin=.75in
     \leftmargin=.5in
     \vsize=6.9in
     \hsize=4.75in
     \fullhsize=10in
     \let\draft=\relax}

\outputstyle{normal}                             % page style

\def\marginnoteformat{\subscriptsize             % paragraphing of margin notes
     \hsize=1in \baselinestretch=1000 \everypar={}%
     \tolerance=5000 \hbadness=5000 \parskip=0pt \parindent=0pt
     \leftskip=0pt \rightskip=0pt \raggedright}

\head={\ifdraft\normalfonts\it\hfil DRAFT\hfil   % format of headline
     \llap{\number\day\ \monthword\month\ \militarytime}\else\hfil\fi}
\foot={\hfil\normalfonts\numstyle\pagenum\hfil}  % format of footline

\normalbaselineskip=12pt                         % usual \baselineskip
\normallineskip=0pt                              % usual \lineskip
\normallineskiplimit=0pt                         % usual \lineskiplimit
\normalbaselines                                 % set \baselineskip

\topskip=.85\baselineskip \splittopskip=\topskip \headskip=2\baselineskip
\footskip=\headskip

\pagenumstyle{arabic}                            % counter style

\parskip=0pt                                     % no skip between paragraphs
\parindent=20pt                                  % usual \parindent

\baselinestretch=1000                            % set \big-, \med-, \smallskip

%************** Sectioning **********************************

\chapterstyle{left}                              % position of heading
\chapternumstyle{blank}                          % counter style
\def\chapterbreak{\newpage}                      % break before heading
\abovechapterskip=0pt                            % space before heading
\belowchapterskip=1.5\baselineskip               % space after heading
     plus.38\baselineskip minus.38\baselineskip
\def\chapternumformat{\numstyle\chapternum.}     % format of heading counter

\sectionstyle{left}                              % position of heading
\sectionnumstyle{blank}                          % counter style
\def\sectionbreak{\vskip0pt plus4\baselineskip\penalty-100
     \vskip0pt plus-4\baselineskip}              % break before heading
\abovesectionskip=1.5\baselineskip               % space before heading
     plus.38\baselineskip minus.38\baselineskip
\belowsectionskip=\the\baselineskip              % space after heading
     plus.25\baselineskip minus.25\baselineskip
\def\sectionnumformat{%                          % format of heading counter
     \ifblank\chapternumstyle\then\else\numstyle\chapternum.\fi
     \numstyle\sectionnum.}

\subsectionstyle{left}                           % position of heading
\subsectionnumstyle{blank}                       % counter style
\def\subsectionbreak{\vskip0pt plus4\baselineskip\penalty-100
     \vskip0pt plus-4\baselineskip}              % break before heading
\abovesubsectionskip=\the\baselineskip           % space before heading
     plus.25\baselineskip minus.25\baselineskip
\belowsubsectionskip=.75\baselineskip            % space after heading
     plus.19\baselineskip minus.19\baselineskip
\def\subsectionnumformat{%                       % format of heading counter
     \ifblank\chapternumstyle\then\else\numstyle\chapternum.\fi
     \ifblank\sectionnumstyle\then\else\numstyle\sectionnum.\fi
     \numstyle\subsectionnum.}

%************** Footnotes ***********************************

\footnotenumstyle{symbols}                       % counter style
\footnoteskip=0pt                                % jyTeX spacing parameter
\def\footnotenumformat{\numstyle\footnotenum}    % \footnotemark format
\def\footnoteformat{\footnotesize                % paragraphing of text
     \everypar={}\parskip=0pt \parfillskip=0pt plus1fil
     \leftskip=1em \rightskip=0pt
     \spaceskip=0pt \xspaceskip=0pt
     \def\\{\ifhmode\ifnum\lastpenalty=-10000
          \else\hfil\penalty-10000 \fi\fi\ignorespaces}}

%************** Labels **************************************

\def\undefinedlabelformat{$\bullet$}             % mark for undefined label

%************** Equation numbering **************************

\equationnumstyle{arabic}                        % counter style
\subequationnumstyle{blank}                      % counter style
\figurenumstyle{arabic}                          % counter style
\subfigurenumstyle{blank}                        % counter style
\tablenumstyle{arabic}                           % counter style
\subtablenumstyle{blank}                         % counter style

\eqnseriesstyle{alphabetic}                      % sub-counter style for series
\figseriesstyle{alphabetic}                      % sub-counter style for series
\tblseriesstyle{alphabetic}                      % sub-counter style for series

\def\puteqnformat{\hbox{%                        % equation number format
     \ifblank\chapternumstyle\then\else\numstyle\chapternum.\fi
     \ifblank\sectionnumstyle\then\else\numstyle\sectionnum.\fi
     \ifblank\subsectionnumstyle\then\else\numstyle\subsectionnum.\fi
     \numstyle\equationnum
     \numstyle\subequationnum}}
\def\putfigformat{\hbox{%                        % figure number format
     \ifblank\chapternumstyle\then\else\numstyle\chapternum.\fi
     \ifblank\sectionnumstyle\then\else\numstyle\sectionnum.\fi
     \ifblank\subsectionnumstyle\then\else\numstyle\subsectionnum.\fi
     \numstyle\figurenum
     \numstyle\subfigurenum}}
\def\puttblformat{\hbox{%                        % table number format
     \ifblank\chapternumstyle\then\else\numstyle\chapternum.\fi
     \ifblank\sectionnumstyle\then\else\numstyle\sectionnum.\fi
     \ifblank\subsectionnumstyle\then\else\numstyle\subsectionnum.\fi
     \numstyle\tablenum
     \numstyle\subtablenum}}

%************** Reference numbering *************************

\referencestyle{sequential}                      % referencing method
\referencenumstyle{arabic}                       % counter style
\def\putrefformat{\numstyle\referencenum}        % format of reference citation
\def\referencenumformat{\numstyle\referencenum.} % format of number in list
\def\putreferenceformat{%                        % paragraphing of list
     \everypar={\hangindent=1em \hangafter=1 }%
     \def\\{\hfil\break\null\hskip-1em \ignorespaces}%
     \leftskip=\refnumindent\parindent=0pt \interlinepenalty=1000 }

%************** Font initialization *************************

\normalsize

%*****************************************************************************

\def\fmtversion{2.6M (June 1992)}

\catcode`\@=12
% ------------------ End of jytex.tex -----------------

%\input jytex.tex   % available from hep-th
\typesize=10pt \magnification=1200 \baselineskip17truept
%\baselineskip25truept
\footnotenumstyle{arabic} \hsize=6truein\vsize=8.5truein
%\input Spinharm.lab
%\draft
%\leftmargin=1.25in
%\oddleftmargin=.5in
%\evenleftmargin=1.5in
\sectionnumstyle{blank}
\chapternumstyle{blank}
\chapternum=1
\sectionnum=1
\pagenum=0
%\referencestyle{preordered}
% title style follows

\def\begintitle{\pagenumstyle{blank}\parindent=0pt
\begin{narrow}[0.4in]}
\def\endtitle{\end{narrow}\newpage\pagenumstyle{arabic}}

% exercise style follows

\def\beginexercise{\vskip 20truept\parindent=0pt\begin{narrow}[10
truept]}
\def\endexercise{\vskip 10truept\end{narrow}}

% **************    my jyTeX abbreviations   *****************

\def\eql#1{\eqno\eqnlabel{#1}}
\def\ref{\reference}
\def\peq{\puteqn}
\def\pref{\putref}

\def\mgn{\marginnote}
\def\bex{\begin{exercise}}
\def\eex{\end{exercise}}

% *********************** My definitions ************************

\font\open=msbm10 %scaled\magstep1 % For VAX. Borde p195.

 %scaled\magstep1 % For VAX. Borde p195.
%\font\open=msym10 %scaled\magstep1 % For Arbortxt on PC
%\font\opens=msym8 %scaled\magstep1 % For Arbortxt on PC
  % For Arbortxt on PC, and VAX. Borde p199
 
\def\StretchRtArr#1{{\count255=0\loop\relbar\joinrel\advance\count255 by1
\ifnum\count255<#1\repeat\rightarrow}}
\def\StretchLtArr#1{\,{\leftarrow\!\!\count255=0\loop\relbar
\joinrel\advance\count255 by1\ifnum\count255<#1\repeat}}

\def\StretchLRtArr#1{\,{\leftarrow\!\!\count255=0\loop\relbar\joinrel\advance
\count255 by1\ifnum\count255<#1\repeat\rightarrow\,\,}}

\def\mbox#1{{\leavevmode\hbox{#1}}}

\def\hspace#1{{\phantom{\mbox#1}}}
\def\oR{\mbox{\open\char82}}

\def\oZ{\mbox{\open\char90}}

\def\oN{\mbox{\open\char78}}

\def\al{\alpha}
\def\bal{{\bmit\alpha}} %in jyTeX
 %in jyTeX
\def\bbe{{\bmit\beta}} %in jyTeX
 %in jyTeX
 %in jyTeX
 %in jyTeX
 %in jyTeX
 %in jyTeX
 %in jyTeX
 %in jyTeX
 %in jyTeX
 %in jyTeX
 %in jyTeX
% in jyTeX
% in jyTeX
% in jyTeX
% in jyTeX
% in jyTeX
\def\be{\beta}

\def\de{\delta}
\def\Ga{\Gamma}

\def\la{\lambda}

\def\om{\omega}
\def\Om{\Omega}

\def\si{\sigma}

\def\ze{\zeta}

\def\De{\Delta}

\def\caP{{\cal P}}

\def\sc{{\rm sc }}

\def\zf{$\zeta$--function}

     % Newline

\def\frac#1/#2{\leavevmode\kern.1em
\raise.5ex\hbox{\the\scriptfont0 #1}\kern-.1em/\kern-.15em
\lower.25ex\hbox{\the\scriptfont0 #2}}
\def\sfrac#1/#2{\leavevmode\kern.1em
\raise.5ex\hbox{\the\scriptscriptfont0 #1}\kern-.1em/\kern-.15em
\lower.25ex\hbox{\the\scriptscriptfont0 #2}}

\def\gtorder{\mathrel{\raise.3ex\hbox{$>$}\mkern-14mu
             \lower0.6ex\hbox{$\sim$}}}
\def\ltorder{\mathrel{\raise.3ex\hbox{$<$}\mkern-14mu
             \lower0.6ex\hbox{$\sim$}}}

\def\semidirprod{\rlap{\ss C}\raise1pt\hbox{$\mkern.75mu\times$}}
\def\for{\lower6pt\hbox{$\Big|$}}
\def\fish{\kern-.25em{\phantom{abcde}\over \phantom{abcde}}\kern-.25em}

 %triple
%dot
 %double
%dot
 %double dot
%for small #1

\def\boxit#1{\vbox{\hrule\hbox{\vrule\kern3pt
        \vbox{\kern3pt#1\kern3pt}\kern3pt\vrule}\hrule}}
\def\dalemb#1#2{{\vbox{\hrule height .#2pt
        \hbox{\vrule width.#2pt height#1pt \kern#1pt \vrule
                width.#2pt} \hrule height.#2pt}}}

\def\ol{\overline}
        %double stroke
\def\frac#1#2{{{#1}\over{#2}}}
 %lower covariant deriv.
 %upper covariant deriv.
 %lower covariant deriv semicolon.
    %lower ordinary  deriv.
    %lower ordinary  deriv comma.

\def\noin{\noindent}

      %Connection
    %Connection'
\def\comb#1#2{{\left(#1\atop#2\right)}}

\def\etc{{\it etc. }}

\def\eg{{\it e.g.}}
\def\ie{{\it i.e. }}
\def\cf{{\it cf }}
\def\pa{\partial}

 %gives average <#1>
 %gives thermal average <<#1>>
   %gives bracket <#1|#2>
   %gives comma bracket <#1,#2>
 %gives round bracket (#1,#2)
 %gives round bracket (#1,|#2)
 %gives big bracket <#1|#2>
  %gives
%matrix element <#1|#2|#3>

%gives reduced matrix element
%<#1||#2||#3>

\def\3j#1#2#3#4#5#6{\left\lgroup\matrix{#1&#2&#3\cr#4&#5&#6\cr}
\right\rgroup}

\def\man{{\cal M}}

\def\m?{\mgn{?}}
% KK's defs

\def\pa{\partial}

\def\beq{\begin{eqnarray}}
\def\eeq{\end{eqnarray}}

%  *******************  Journal refs **********************

\def\aop#1#2#3{{\it Ann. Phys.} {\bf {#1}} ({#2}) #3}
\def\cjp#1#2#3{{\it Can. J. Phys.} {\bf {#1}} ({#2}) #3}
\def\cmp#1#2#3{{\it Comm. Math. Phys.} {\bf {#1}} ({#2}) #3}
\def\cqg#1#2#3{{\it Class. Quant. Grav.} {\bf {#1}} ({#2}) #3}

\def\ijmp#1#2#3{{\it Int. J. Mod. Phys.} {\bf {#1}} ({#2}) #3}

\def\jmp#1#2#3{{\it J. Math. Phys.} {\bf {#1}} ({#2}) #3}
\def\jpa#1#2#3{{\it J. Phys.} {\bf A{#1}} ({#2}) #3}
\def\jpc#1#2#3{{\it J. Phys.} {\bf C{#1}} ({#2}) #3}
\def\lnm#1#2#3{{\it Lect. Notes Math.} {\bf {#1}} ({#2}) #3}

\def\np#1#2#3{{\it Nucl. Phys.} {\bf B{#1}} ({#2}) #3}
\def\npa#1#2#3{{\it Nucl. Phys.} {\bf A{#1}} ({#2}) #3}
\def\pl#1#2#3{{\it Phys. Lett.} {\bf {#1}} ({#2}) #3}

\def\prp#1#2#3{{\it Phys. Rep.} {\bf {#1}} ({#2}) #3}
\def\pr#1#2#3{{\it Phys. Rev.} {\bf {#1}} ({#2}) #3}
\def\prA#1#2#3{{\it Phys. Rev.} {\bf A{#1}} ({#2}) #3}

\def\prD#1#2#3{{\it Phys. Rev.} {\bf D{#1}} ({#2}) #3}
\def\prE#1#2#3{{\it Phys. Rev.} {\bf E{#1}} ({#2}) #3}
\def\prl#1#2#3{{\it Phys. Rev. Lett.} {\bf #1} ({#2}) #3}

\def\rmp#1#2#3{{\it Rev. Mod. Phys.} {\bf {#1}} ({#2}) #3}

\def\zfp#1#2#3{{\it Z. f. Phys.} {\bf {#1}} ({#2}) #3}

\def\cras#1#2#3{{\it Comptes Rend. Acad. Sci. (Paris)} {\bf{#1}} (#2) #3}
\def\prs#1#2#3{{\it Proc. Roy. Soc.} {\bf A{#1}} ({#2}) #3}
\def\pcps#1#2#3{{\it Proc. Camb. Phil. Soc.} {\bf{#1}} ({#2}) #3}
\def\mpcps#1#2#3{{\it Math. Proc. Camb. Phil. Soc.} {\bf{#1}} ({#2}) #3}

\def\amsh#1#2#3{{\it Abh. Math. Sem. Ham.} {\bf {#1}} ({#2}) #3}
\def\am#1#2#3{{\it Acta Mathematica} {\bf {#1}} ({#2}) #3}
\def\aim#1#2#3{{\it Adv. in Math.} {\bf {#1}} ({#2}) #3}
\def\ajm#1#2#3{{\it Am. J. Math.} {\bf {#1}} ({#2}) #3}
\def\amm#1#2#3{{\it Am. Math. Mon.} {\bf {#1}} ({#2}) #3}

\def\aom#1#2#3{{\it Ann. of Math.} {\bf {#1}} ({#2}) #3}
\def\cjm#1#2#3{{\it Can. J. Math.} {\bf {#1}} ({#2}) #3}
\def\bams#1#2#3{{\it Bull.Am.Math.Soc.} {\bf {#1}} ({#2}) #3}

\def\cmh#1#2#3{{\it Comm. Math. Helv.} {\bf {#1}} ({#2}) #3}

\def\dmj#1#2#3{{\it Duke Math. J.} {\bf {#1}} ({#2}) #3}
\def\invm#1#2#3{{\it Invent. Math.} {\bf {#1}} ({#2}) #3}

\def\jdg#1#2#3{{\it J. Diff. Geom.} {\bf {#1}} ({#2}) #3}

\def\joa#1#2#3{{\it J. of Algebra} {\bf {#1}} ({#2}) #3}
\def\jram#1#2#3{{\it J. f. reine u. angew. Math.} {\bf {#1}} ({#2}) #3}
\def\jims#1#2#3{{\it J. Indian. Math. Soc.} {\bf {#1}} ({#2}) #3}
\def\jlms#1#2#3{{\it J. Lond. Math. Soc.} {\bf {#1}} ({#2}) #3}
\def\jmpa#1#2#3{{\it J. Math. Pures. Appl.} {\bf {#1}} ({#2}) #3}
\def\ma#1#2#3{{\it Math. Ann.} {\bf {#1}} ({#2}) #3}

\def\mz#1#2#3{{\it Math. Zeit.} {\bf {#1}} ({#2}) #3}
\def\ojm#1#2#3{{\it Osaka J.Math.} {\bf {#1}} ({#2}) #3}

\def\pems#1#2#3{{\it Proc. Edin. Math. Soc.} {\bf {#1}} ({#2}) #3}

\def\plb#1#2#3{{\it Phys. Letts.} {\bf {B#1}} ({#2}) #3}
\def\pla#1#2#3{{\it Phys. Letts.} {\bf {A#1}} ({#2}) #3}
\def\plms#1#2#3{{\it Proc. Lond. Math. Soc.} {\bf {#1}} ({#2}) #3}
\def\pgma#1#2#3{{\it Proc. Glasgow Math. Ass.} {\bf {#1}} ({#2}) #3}
\def\qjm#1#2#3{{\it Quart. J. Math.} {\bf {#1}} ({#2}) #3}
\def\qjpam#1#2#3{{\it Quart. J. Pure and Appl. Math.} {\bf {#1}} ({#2}) #3}

\def\rmjm#1#2#3{{\it Rocky Mountain J. Math.} {\bf {#1}} ({#2}) #3}

\def\tams#1#2#3{{\it Trans.Am.Math.Soc.} {\bf {#1}} ({#2}) #3}

% *******************   Main text *********************
\begin{title}
\vglue 0.5truein
%\righttext {MUTP/96/23}
%\righttext{hep-th/96}
\vskip15truept
%\leftline{\today}
%\vskip 30truept
\centertext {\Bigfonts \bf Relations between the Ehrhart polynomial,} \vskip7truept
\vskip10truept\centertext{\Bigfonts \bf  the heat kernel and Sylvester waves}
 \vskip 20truept
\centertext{J.S.Dowker\footnote{dowker@man.ac.uk}} \vskip 7truept
\centertext{\it Theory Group,} \centertext{\it School of Physics and
Astronomy,} \centertext{\it The University of Manchester,} \centertext{\it
Manchester, England} \vskip 7truept \centertext{}

\vskip 7truept

\vskip40truept
\begin{narrow}
I show for the specific case of the scalar field spectrum on regular tessellations of the
sphere that the first two terms of the heat--kernel expansion are related to the first two
terms of the Ehrhart (quasi)polynomial. In trying to make this relation precise, I consider
degeneracies as partition denumerants and show the connection of the group theory
expressions with Popoviciu's theorem and with the notion of Sylvester waves. General
denumerants are considered and the first wave, $W_1$, \ie the polynomial part, is
written using the $A$--genus multiplicative sequence. It is pointed out that Sylvester in
effect did the same thing and that he had also obtained Ehrhart reciprocity. I derive an
algebraically neat form for the second wave, $W_2$, which involves the combination of
two multiplicative sequences.
\end{narrow}
\vskip 5truept
\vskip 60truept
%\righttext{Typeset in \jyTeX}
\vfil
\end{title}
\pagenum=0
\newpage

\section{\bf 1. Introduction}
Connections between lattice and spectral problems go back over a hundred years. I
mention only the Weyl conjecture and its physical antecedents (see, for example, Baltes
and Hilf, [\pref{BaandH}]).

Two results, which are curiously similar, are the Ehrhart polynomial and the short--time
expansion of the heat--kernel. Actually, rather than the heat--kernel, it is more
appropriate to consider the asymptotic behaviour of a smoothed version, $\ol N(\la)$, of
the exact eigenvalue counting function, $N(\la)=\sum_{k,\la_k\le\la}\,1$.\footnote{ The
eigenvalues, $\la_k$ are ordered linearly,  with $k$ a counting label.} As a power series
in $\la$, this has coefficients simply related to those of the heat--kernel expansion. For
simplicity of exposition, I restrict initially to two dimensions and give the first two terms
(which are all I am interested in) explicitly,
  $$
  \ol N(\om)\sim {1\over4\pi}\bigg(|\man|\,\om^2\pm
  \,|\pa\man|\,\om+ c\,+\ldots\bigg)\,.
  \eql{nasym}
  $$
I have set $\la$  equal to $\om^2$  and redefined $N(\om)=N(\la)$.

These first two terms are those relevant for the Weyl conjecture, which states that
(\peq{nasym}) holds for the {\it exact} counting function. While this may not be true
(actually it isn't for the manifolds I look at), (\peq{nasym}) does hold for the smoothed
version.

Turning to the lattice side, Ehrhart proved that the number of integer lattice points (\ie
elements of $\oZ^2$) in, and on, a rational polygon, \ie a polygon whose vertices are
rational numbers, that has been uniformly dilated by an integer, $l$, is a {\it
quasipolynomial} in $l$.

More precisely let $\caP$ be a rational polygon\footnote{I sometimes refer to a disc with
polygonal boundary (a 2--polytope) as a polygon, for short.} and $\ol \caP$ its closure
(\ie including its boundary) and let $L(\ol P,l)=\sharp(\ol{l\caP}\cap\oZ^2)$ be the
number of integer points in the closure of the dilated polygon
$l\caP=\{(lx,ly):(x,y)\in\caP\,; l\in\oN\}$. Then
  $$
  L(\ol\caP,l)=|\caP|\,l^2+c_1(l)\,l+c_0(t)
  \eql{ehr}
  $$
where $c_1,c_0$ are generally periodic functions in $l$ (see, for example, Beck and
Robins, [{\pref{BandR}] and Wright [\pref{Wright}].

In the special case that the polygon is an {\it integer} polygon the  coefficients
$c_1,c_0$ are constant and, moreover,
  $$
  L(\ol\caP,l)=|\caP|\,l^2+{1\over2}|\pa\caP|\,l+c_0
  \eql{ehr2}
  $$
which bears a remarkable similarity to (\peq{nasym}).

The two expressions cannot be related in all circumstances because not all lattice
counting problems have spectral associations.  This paper, however, is concerned with a
very particular case in which they can be more precisely connected.

In sections 2 and 3, I set up the situation and present an experimental result. The later
sections are intended to be more exact and give some explanation of the result. This
paper thus proceeds from the particular to the general.
\section{\bf 2. The factored sphere}

The situation I refer to is very well known and so I need not spend time on the underlying
mathematical details. They have been dealt with by, for example, Gromes,
[\pref{Gromes}], B\'erard and Besson, [\pref{BandB}], Br\"uning and Heintze,
[\pref{BandHe}] and Chang and Dowker, [\pref{ChandD}]. The manifold, $\man$ is an
orbifold factor of the two--sphere, S$^2/\Ga$, where $\Ga$ is the reflective symmetry
group of one of the regular solids, a finite subgroup of O(3). $\man$ is therefore a
spherical triangle tiling the sphere, or, in the dihedral case, it is a lune (the situation
discussed by Gromes).

Making the special choice of conformal coupling in $1+2$ dimensions means that the
eigenvalues of the corresponding `improved' Laplacian are perfect squares of integers, or
half--odd integers, (for the unit sphere). By an appropriate selection of the action of the
tiling group, $\Ga$, one can arrange Dirichlet or Neumann conditions on the boundary of
the fundamental domain, $\man$. A calculation, by separation of variables if you wish,
\eg\ [\pref{Gromes}] on the lune, gives the eigenvalues, $\la_{\bf m}$, determined by
  $$\eqalign{
 \om_{\bf m}\equiv \sqrt\la_{\bf m}=a+d_1m_1+d_2m_2\,,\quad m_1,m_2=0,1,2,\ldots
 }
 \eql{eigvals}
  $$
where $d_1$ and $d_2$ are the integer {\it degrees} (not necessarily coprime)
associated with the action of $\Ga$, \ie with the particular regular solid (3--polytope).
For example, for the lune of apex angle $\pi/q$ ($q\in\oZ)$, $d_1=q$, $d_2=1$. The
hemisphere corresponds to $q=1$.

The constant, $a$, takes the values
  $$
  a_N={1\over2}\quad {\rm and}\quad a_D=d_1+d_2-a_N
  $$
for Neumann and Dirichlet conditions, [\pref{ChandD}].

One now sees the relation with Ehrhart polynomials. Trivially, counting eigenvalues
amounts to counting integer lattice points. It is no surprise that similar quantities
appear.

Explicitly, consider the (special) rational triangle (2--polytope)
  $$
  \caP_{\bf d}=\bigg\{(x_1,x_2)\in \oR^2:x_i\ge0,d_1x_1+d_2x_2\le1\bigg\}
  $$
with vertices $(0,0)$, $(1/d_1,0)$ and $(0,1/d_2)$ where $d_1$ and $d_2$  are
relatively prime, which is not an essential restriction. Then construct the integer dilation,
$l\caP_{\bf d}$, and count the number, $L(\ol\caP_{\bf d},l)$  of integer lattice points,
$(m_1,m_2)$, inside (including the boundary as described previously). The result is
given by (\peq{ehr}) where now in fact $c_1$ is a constant (because $d_1$ and $d_2$
are coprime). Hence, so far,
  $$
  L(\ol\caP_{\bf d},l)={l^2\over 2d_1d_2}+c_1l +
  c_0(l)\,.
  $$

The complete expression has been obtained by Beck and Robins, [\pref{BandR}], and I
reproduce the essential parts here,
  $$
  L(\ol\caP_{\bf d},l)={1\over2d_1d_2}\bigg(l^2+l\big(d_1+d_2+1\big)+
  {1\over6}\big(d_1^2+d_2^2+3d_1d_2+1\big)
   +\si(l)\bigg)\,.
   \eql{ehr3}
  $$
where $\si(l)$ is a periodic function involving Dedekind sums. The origin, Ehrhart
[\pref{Ehrhart}], should also be consulted.

To make a connection with eigenvalues, the essential counting restriction is
   $$
   d_1m_1+d_2m_2\le l
   $$
which, in terms of eigenvalues, (\peq{eigvals}), reads
  $$
  \om_{\bf m}\le l+a\,.
  $$
so that the identification
  $$
  N(\om)= L(\ol\caP_{\bf d},\om-a)
  \eql{ident}
  $$
can be made {\it but} only when $\om-a$ is an integer. As a function of a real $l$,
$L(\ol\caP_{\bf d},l)$ will provide an interpolation and then the corresponding
interpolated  $N(\om)$ from (\peq{ident}) would give a smoothing of the exact
$N(\om)$. To see how this compares with the smoothed expression (\peq{nasym}), I
put in the geometric values appropriate for a fundamental domain on the 2--sphere to
get for (\peq{nasym}),
  $$
  \ol N(\om)\sim {1\over2d_1d_2}\,\bigg(\om^2\pm
  (d_1+d_2-1)\,\om+ c'\,+\ldots\bigg)\,.
  \eql{nasym2}
  $$
The factor $2d_1d_2$ is the order of the reflective tiling group, $\Ga$ and the $d_1$ and
$d_2$ need not be coprime.

Choosing the upper (Neumann) sign for ease ($a=1/2$), substitution of $l=\om-1/2$ in
(\peq{ehr3}) does not yield (\peq{nasym2}). However, all is not lost. Plotting the exact
lattice counting expression (\peq{ehr3}) at integer values of $l$ produces a familiar
staircase function. The interpolation provided by (\peq{ehr3}) for real $l$ joins the upper
part of the steps. The curve that joins the lower part is given by $L(\ol\caP_{\bf d},l-1)$
and averaging these produces an interpolation (\ie smoothing) that passes half way up
the vertical rises. It is easily confirmed that this symmetrical combination
  $$
  {1\over2}\big( L(\ol\caP_{\bf d},\om-1/2)+L(\ol\caP_{\bf d},\om-3/2)\big)
  \eql{sym}
  $$
reproduces the first two terms of (\peq{nasym2}), with the upper sign. This is the
conclusion so far.

Also one obtains the equality
  $$
  |\man|=4\pi|\caP_{\bf d}|\,.
  $$
\section{\bf 3. Higher dimensions}

To show that this result possibly goes beyond mere numerical coincidence (which is
unlikely anyway in view of the dependence on $d_1$ and $d_2$), I look at three
dimensions more particularly.

Generally, in $d$ dimensions, the measures of the fundamental domain and its boundary
are given by
  $$
  |\man|={1\over2g}|S^d|\,\quad {\rm and}\quad |\pa\man|={b_1\over g}|S^{d-1}|
  $$
in terms of the order, $2g$, of the tiling group, $\Ga$, and the number, $b_1$, of
reflecting great hyperspheres (or reflecting $d$--flats  in the embedding
$(d+1)$--dimensional Euclidean manifold). In terms of the degrees
  $$
  g=\prod_{i=1}^d d_i\quad {\rm and}\quad b_1=\sum_{i=1}^d (d_i-1)+1
  $$
so that, in particular for $d=3$, $g=d_1d_2d_3$ and $b_1=d_1+d_2+d_3-2$.

Writing out the conventional asymptotic behaviour of the smoothed counting function in
this particular situation gives
   $$
  \ol N(\om)\sim {1\over2g\Ga(d)}\bigg({2\over d}\om^d\pm
  b_1\,\om^{d-1}+{d-1\over6}\big(b_1(b_1-1)+b_2\big) \om^{d-2} +\ldots\bigg)\,.
  \eql{nasym4}
  $$
where I have included the next term which involves $b_2$, the number of elements of
the group $\Ga$ that fix a $(d-1)$--flat but not a $d$--flat. In three dimensions,
  $$
  b_2=d_2d_3+d_3d_1+d_1d_2-d_1-d_2-d_3\,.
  $$

In three dimensions the Neumann constant, $a_N$, equals 1 and the suggested
symmetrical (midway) combination is, instead of (\peq{sym}),
   $$
  {1\over2}\big( L(\ol\caP_{\bf d},\om-1)+L(\ol\caP_{\bf d},\om-2)\big)\,.
  \eql{sym2}
  $$

The relevant polynomial terms have been determined by Beck {\it et al}, [\pref{BDR}],
and again I write them out
  $$\eqalign{
   L(\ol\caP_{\bf d},l)=&{1\over 2d_1d_2d_3}\bigg({l^3\over3}+{l^2\over2}
   \big(d_1+d_2+d_3+1\big)+\cr
   &{l\over6}\big(3(d_1+d_2+d_3+d_2d_3+d_3d_1+d_1d_2)
   +d_1^2+d_2^2+d_3^2\big)+\ldots\bigg)\,,
   }
   \eql{ehr4}
  $$
where the omitted terms are periodic.

Evaluating (\peq{sym})  yields the expression
    $$
    {1\over 2d_1d_2d_3}\bigg({\om^3\over3}+b_1{\om^2\over2}+
    {\big(b_1(b_1-1)+b_2+1\big)\om\over6}\,+\ldots\bigg)
    $$
the first two terms of which agree with the `standard' expansion (\peq{nasym4}).  The
third term differs, but only very slightly, which proximity I cannot explain. One might
expect that one should reproduce those polynomial coefficients which are constant, \ie all
those except the final one, proportional to $\om^0$, although the higher heat--kernel
coefficients depend on the propagation equation.

\section{\bf 4. Degeneracies as denumerants}

I now wish to proceed to some exact relations, with the aim of making the previous
conjecture more sensible. They are also more interesting.

Averaging the counting function involves loss of information. I return to  the exact
function which determines, and is determined by, the (square root) eigenvalue spectrum,
$\{\om_k\}$, $k=1,2,\ldots$, coincidences implying degeneracy.

The spectrum can also be specified by the distinct eigen{\it levels}, $\om(l)$ say,
$l=0,1,2,\ldots$, and the corresponding degeneracies, $g(l)$. On the fundamental
domain, leaving the eigenvalues in the form (\peq{eigvals}) expressed in terms of two
integers (or more for higher spheres) one can often avoid explicitly constructing the
degeneracies. However degeneracies have a useful role and appear in other contexts.

  For a given Neumann eigenlevel, $\om(l)$, the number of integer lattice points
$(m_1,m_2)$ satisfying
  $$
  d_1m_1+d_2m_2=l=\om(l)-a\,,\quad {\rm with}\quad a={1\over2}\,,
  \eql{cond}
  $$
gives the degeneracy, \ie
  $$
  g_N(l;{\bf d})=\sharp\bigg( (m_1,m_2)\in \oZ_+^2:d_1m_1+d_2m_2=l\bigg)\,.
  $$

$g_N(l;{\bf d})$ counts the number of lattice points on the hypotenuse of the dilated
rational triangle $l\caP_{\bf d}$ and is often referred to as a restricted partition or a
denumerant, in an older terminology. Bell, [\pref{Bell}], denotes the obvious extension
to higher dimensions by $D(l\mid{\bf d})$ and Sylvester employs a variety of notations
and terminology. His final one seems to be
  $$
  D(l\mid{\bf d})={l\over d_1,d_2,\ldots,d_d,}\,,
  $$
which I shorten to $l/{\bf d}\,,$.\footnote{ I retain the geometric terminology of
`degrees' for the denominator of the denumerant. Sylvester refers to them as the
`components' of the denumerant, among other things. Another possibility would be
`parameters'.}

The study of denumerants is effectively the same as that of Ehrhart (quasi) polynomials
as discussed later. Ehrhart himself, [\pref{Ehrhart}], actually spends most analytical
time on them.

The associated generating function is well known, going back at least to Euler,
[\pref{Sylvester4}],
  $$
   h_N(\si;{\bf d})\equiv\sum_{l=0}^\infty g_N(l;{\bf d})\,\si^l
   ={1\over(1-\si^{d_1})(1-\si^{d_2})}\,.
   \eql{genfun}
  $$

Of course, as a lattice statement,  this holds for any  integers $d_1$ and $d_2$. For the
S$^2/\Ga$ spectral problem, $d_1$ and $d_2$ take only certain values and the
generating function (\peq{genfun}) is just the Molien series for a finite reflection group
action, \eg\ Stanley, [\pref{Stanley}], Meyer [\pref{Meyer}] and was explored by
Laporte, [\pref{Laporte}], from a physical, mode point of view.  The extension to higher
dimensions is immediate. See Br\"uning and Heintze, [\pref{BandHe}], B\'erard and
Besson, [\pref{BandB}].

Dirichlet conditions correspond to a different value for $a$ in (\peq{cond}) and amount
to a shift in $l$ by a constant $d_0$ where
  $$
 d_0=d_1+d_2-1\,,
  $$
and is the number of reflecting hyperplanes in the embedding space.

Then the Dirichlet generating function is
  $$
   h_D(\si;{\bf d})\equiv\sum_{l=0}^\infty g_D(l;{\bf d})\,\si^l=
   {\si^{d_0}\over(1-\si^{d_1})(1-\si^{d_2})}\,,
   \eql{genfun2}
  $$
or
  $$
  g_D(l+d_0;{\bf d})=g_N(l;{\bf d})\,,\quad l=0,1,2,\ldots\,.
  \eql{reln}
  $$

 As an important example, for the lune of angle $\pi/q$, $q\in\oZ$,
   $$\eqalign{
   h_N(\si;q,1)=&{1\over(1-\si^q)(1-\si)}\cr
   h_D(\si;q,1)=&{\si^q\over(1-\si^q)(1-\si)}\cr
   =&h_N(\si;q,1)-{1\over1-\si}\,.
  }
   $$

A classic roots of unity, or trigonometric, calculation  gives the explicit values for the
degeneracies
   $$
   g_N(l;q,1)=\lfloor l/q\rfloor+1=\lceil l/q\rceil\,,\quad  g_D(l;q,1)=\lfloor l/q\rfloor,
   $$
the solutions for the corresponding denumerants, (\peq{cond}), \eg\ Sylvester
[\pref{Sylvester}].

Combining the Dirichlet and Neumann spectra (which amounts to adding spectral
quantities such as the generating functions) gives the spectrum for a doubled
fundamental domain on which the modes are periodic under the pure rotational part of
$\Ga$. Thus for the periodic lune (of angle $2\pi/q$) we get the standard formula,
  $$\eqalign{
  h(\si;q,1)\equiv & h_N(\si;q,1)+ h_D(\si;q,1)={1+\si^q\over(1-\si^q)(1-\si)}\cr
  =&\sum_{t=0}^\infty\big(2\lfloor t/q\rfloor+1\big)\si^t\,,
  }
  \eql{cyclic}
  $$
which exhibits the usual degeneracies sometimes obtained from group characters, \eg\
[\pref{ChandD}], Harmer, [\pref{Harmer}].

For the finite number of periodic, uniform sphere tilings, the generating functions can be
obtained by combining cyclic expressions corresponding to a geometric decomposition of
the action of the subgroup of SO(3) in terms of axes and orders,
[\pref{Meyer,PandM,ChandD}]. Using the orbit--stabiliser theorem, this leads to the
elegant result, [\pref{Dow}],
  $$
  h(\si;{\bf d})={1\over2}\bigg(\sum _q h(\si;q,1)-h(\si;1,1)\bigg)
  $$
where the sum is over all axes of order $q$. As an example, the tetrahedral group has
$d_1=3$, $d_2=4$ and $q=2,3,3$. Simple algebra yields, using just  knowledge of the
$q$s,
  $$
  h(\si;3,4)={1+\si^6\over(1-\si^3)(1-\si^4)}\,,
  $$
from which the values of $d_1$ and $d_2$ could be read off, if they weren't known.

Using the cyclic degeneracies, (\peq{cyclic}), one derives the explicit expression for the
rotational tetrahedral degeneracies,
  $$\eqalign{
  g(l;3,4)=&\bigg\lfloor{l\over2}\bigg\rfloor+2\bigg\lfloor{l\over3}\bigg\rfloor+1-l\cr
  =&{l\over6}+1-\bigg\{{l\over2}\bigg\}-2\bigg\{{l\over3}\bigg\}\,.
  }
  \eql{Tdegen}
  $$

For comparison, I reproduce the corresponding formulae for the octahedral and
icosahedral tilings, [\pref{Dow}],
  $$\eqalign{
  g(l;4,6)=&\bigg\lfloor{l\over2}\bigg\rfloor+\bigg\lfloor{l\over3}\bigg\rfloor
  +\bigg\lfloor{l\over4}\bigg\rfloor+1-t\cr
  =&{l\over12}+1-\bigg\{{l\over2}\bigg\}-\bigg\{{l\over3}\bigg\}-\bigg\{{l\over4}\bigg\}\,,
  }
  \eql{Odegen}
  $$
and
  $$\eqalign{
  g(l;6,10)=&\bigg\lfloor{l\over2}\bigg\rfloor+\bigg\lfloor{l\over3}\bigg\rfloor
  +\bigg\lfloor{l\over5}\bigg\rfloor+1-l\cr
  =&{l\over30}+1-\bigg\{{l\over2}\bigg\}-\bigg\{{l\over3}\bigg\}-\bigg\{{l\over5}\bigg\}\,.
  }
  \eql{Idegen}
  $$

Returning now to strictly lattice considerations, the solution of  (\peq{genfun}) for the
lattice quantity, $g_N(t;{\bf d})$, with any coprime $d_1$ and $d_2$,  is given by
Popoviciu's theorem,
  $$
  g_N(l;{\bf d})={l\over d_1d_2}+1
  -\bigg\{{d_1^{-1}l\over d_2}\bigg\}-\bigg\{{d_2^{-1}l\over d_1}\bigg\}\,,
  \eql{pop}
  $$
where $d_2^{-1}$ is the mod$(d_1)$ inverse of $d_2$ and $d_1^{-1}$ that of $d_1$,
mod$(d_2)$.

I wish to apply this to the tetrahedral case to check (\peq{Tdegen}). I note that
$3^{-1}\big|_{\rm mod 4}=3$ and $4^{-1}\big|_{\rm mod 3}=1$ so that
(\peq{pop}) reads
  $$\eqalign{
  g_N(l;3,4)&={l\over 12}+1
  -\bigg\{{3l\over 4}\bigg\}-\bigg\{{l\over 3}\bigg\}\cr
  &=\bigg\lfloor{3l\over4}\bigg\rfloor+\bigg\lfloor{l\over3}\bigg\rfloor+1-l\,.
  }
  \eql{pop2}
  $$

This is shown to be consistent with the rotational tetrahedral degeneracies,
(\peq{Tdegen}), by constructing the combination
  $$
   g_N(l;3,4)+g_N(l-6;3,4)={l\over6}+{3\over2}-\bigg\{{3l\over 4}+{1\over2}\bigg\}
   -\bigg\{{3l\over 4}\bigg\}-2\bigg\{{l\over 3}\bigg\}
  $$
according to the relation (\peq{reln}). Equivalence results in view of the equality,
  $$
  \bigg\{{3l\over 4}+{1\over2}\bigg\}+\bigg\{{3l\over 4}\bigg\}
  =\bigg\{{l\over2}\bigg\}+{1\over2}\,,
  $$
essentially just a Hermite identity.

This result means that one has explicit, {\it separate} formulae for the Neumann and
Dirichlet mode numbers,  discussed case by case by Laporte, [\pref{Laporte}]. I now
check the other tilings,

For this I require an extension of Popoviciu's theorem which is, (see Beck and Robins,
[\pref{BandR}] ex.1.28) that, if ${\rm gcd}(d_1,d_2)=e$, then
  $$\eqalign{
  g_N(l;{\bf d})&={le\over d_1d_2}-\bigg\{{\de_1l\over d_2}\bigg\}
  -\bigg\{{\de_2l\over d_1}\bigg\}+1\,,\quad {\rm if}\quad e|l\cr
  &=0\,,\quad {\rm otherwise}\,.
  }
  \eql{pop3}
  $$
where $\de_1 d_1/e=1$ mod $(d_2/e)$ and  $\de_2 d_2/e=1$ mod $(d_1/e)$.

In the octahedral and icosahedral cases, the greatest common divisor is 2 and $l$ must
be even for non--zero $g_N$ in (\peq{pop3}). In the octahedral case, the reduced mod
inverses are $\de_1=2$ and $\de_2=1$ so that
  $$\eqalign{
  g_N(l;4,6)&={l\over 12}+1
  -\bigg\{{l\over 4}\bigg\}-\bigg\{{l\over 3}\bigg\}\,,\quad l\quad{\rm even}\cr
  }
  \eql{pop4}
  $$
which agrees with (\peq{Odegen}) because $\big\{l/2\big\}$ is zero for even $l$.

The Dirichlet modes are obtained by applying (\peq{reln}) to give
  $$\eqalign{
  g_D(l;4,6)&={l-9\over 12}+1
  -\bigg\{{l-9\over 4}\bigg\}-\bigg\{{l-9\over 3}\bigg\}\,,\quad l\quad{\rm odd}\cr
  &={l\over 12}+{1\over4}
  -\bigg\{{l-1\over 4}\bigg\}-\bigg\{{l\over 3}\bigg\}\cr
  &={l\over 12}+{1\over2}
  -\bigg\{{l\over 4}\bigg\}-\bigg\{{l\over 3}\bigg\}\cr
  }
  \eql{pop5}
  $$
again agreeing with (\peq{Odegen}), for $l$ odd.

A similar calculation holds for the icosahedral case where $\de_1=2$, $\de_2=2$.
Substitution into (\peq{pop3}) immediately yields equality with the icosohedral
degeneracies, (\peq{Idegen}).

In these two cases the rotational formula gives both the N and D values because these
are associated with even and odd $l$ respectively, whilst they are mixed up for the
tetrahedron. A similar cross linking complication arises when computing spectral
quantities for twisted vector bundles over homogeneous factors of the three--sphere,
[\pref{Dowgta}].

For the tetrahedron, because the degrees are coprime, there exists a Frobenius number
which equals $d_1d_2-d_1-d_2=5$, meaning that there is at least one $N$--mode for
every integer $l$ greater then $5$. This is not true for the octahedron or icosahedron as
the degrees are both even.

\section{\bf 5. Sylvester waves}

It is seen that all the relevant quantities take the standard form of a polynomial in $l$
plus a periodic term, in agreement with an ancient partition theorem of Sylvester's, first
stated in [\pref{Sylvester}] to the effect that the denumerant  takes the form,
  $$
  {l\over{\bf d,}}=P(l)+U\,,
  \eql{sylv3}
  $$
where $P(l)$  is a polynomial in $l$ of degree $d-1$ (I work in $d$ dimensions) and $U$,
the `undulant' part, is periodic and contains roots of unity.

Because of this periodicity, Sylvester refers to the parts of $l/{\bf d,}$ in (\peq{sylv3})
as `waves', depending on all roots of unity (most waves vanish by natural selection).

For example, in the rotational tetrahedral degeneracy, (\peq{Tdegen}), the periodic part
has period 6 (the individual N and D parts having period 12). The leading term of the
polynomial part is classic, [\pref{PandS}].

Sylvester, [\pref{Sylvester,Sylvester2}], gives the rule  for calculating the polynomial
part, $P$, as the coefficient of $1/t$ in the expansion of the quantity
  $$
   {e^{\,lt}\over\prod_{i=1}^d\big(1-e^{-d_it}\big)}\,,
   \eql{sylv}
  $$
which is clearly a polynomial in $l$.

By definition, this is directly expressed as a generalised Bernoulli function,
  $$\eqalign{
  P_d&={1\over(d-1)!\prod_i d_i}B^{(d)}_{d-1}\big(l+\sum_id_i\mid{\bf d}\big)\cr
    &={(-1)^{d-1}\over(d-1)!\prod_i d_i}B^{(d)}_{d-1}\big(-l\mid{\bf d}\big)
    }
    \eql{sylv2}
  $$
which is a polynomial in $l$ whose coefficients are generalised Bernoulli numbers (up to
a factor).

This expression for the polynomial part of the denumerant has been noted by Rubinstein
and Fel, [\pref{RandF}] who also point out the equivalence with a direct evaluation by
Beck, Gessel and Komatsu, [\pref{BGK}] which seems more or less the same as
Sylvester's, [\pref{Sylvester2}].

The representation (\peq{sylv2}) is, in some ways, only cosmetically compact. It has
value when general properties are being investigated and recursion relations \etc can be
brought in. These can also be used to find explicit forms for particular cases, as discussed
by Norlund, [\pref{Norlund}] and the Bernoulli polynomials can be computed in various
ways. However, I prefer here to use the relation, employed in [\pref{Dow20}], with Todd
polynomials, $T_n$, [\pref{Hirzebruch}],
  $$
  B^{(d)}_{r}(x\mid{\bf d})=(-1)^r\,r!
  \sum_{s=0}^{r}(-1)^s{x^s\over s!}T_{r-s}(\si_1,\ldots,\si_{r-s})\,,\quad r\le d\,,
  \eql{todd}
  $$
where the $\si_s$ are symmetric functions of the degrees $d_i$ ($i=1,\ldots,d$). This
form has the advantage of being valid for all $d$. That is, the {\it functional form} in the
$\si_i$ of the right--hand side does not depend on $d$. Accordingly, the expansion
reads,

  $$\eqalign{
  P_d&={(-1)^{d-1}\over (d-1)!\prod_i d_i}B^{(d)}_{d-1}\big(-l\mid{\bf d}\big)\cr
  &={(-1)^{d-1}\over (d-1)!\prod_i d_i}\sum_{k=0}^{d-1}\comb {d-1}k (-l)^k
  B^{(d)}_{d-k-1}[{\bf d}]\cr
  &={1\over\prod_id_i}\sum_{k=0}^{d-1}\,{l^k\over k!}\,T_{d-1-k}
  \big(\si_1,\ldots,\si_{d-1-k}\big)\cr
  &={1\over\prod_id_i}\bigg({l^{d-1}\over(d-1)!}T_0+{l^{d-2}\over(d-2)!}T_1+
  {l^{d-3}\over(d-3)!}T_2+\ldots\bigg)
  }
  \eql{wave4}
  $$

Known forms for the Todd polynomials  make the coefficients explicit, which allows
checking, but the real advantage of using this structure is that the Todd polynomials are
independent of $d$ and one sees that the expression for dimension $d$ is obtained by
integrating that for $(d-1)$, the only `new' term being the constant. I now present
another (equivalent) version of this expansion.

From an algebraic viewpoint a more advantageous variable is $\ol
l=l+{1\over2}\sum_id_i$. Sylvester, [\pref{Sylvester3}], calls this the {\it augmented}
argument and also extends $l$ to a real number so allowing the use of analysis, as, in
fact, I have just done.\mgn{EXTEND}

A basic property of the Bernoulli functions, [\pref{Norlund,Norlund1}], shows that, {\it
expressed as a function of} $\,\ol l$, the polynomial part has the parity property,
  $$
  P_d(-\ol l)=(-1)^{d-1}P_d(\ol l)\,.
  \eql{parity}
  $$

Furthermore, the generating function equation takes on a more symmetrical aspect in
terms of $\ol l$,
    $$
    \sum_{l=0}^\infty {l\over{\bf d},}\,\si^{\ol l}
    ={1\over\prod_i\big(\si^{-d_i/2}-\si^{d_i/2}\big)}\,,
    $$
which shows that, as a function of $\ol l$, the full denumerant also obeys the parity
property (Ehrhart reciprocity). It follows that the undulant part satisfies it also. This was
proved by Ehrhart geometrically who thence {\it derived} the parity behaviour of the
polynomial part. It seems that Sylvester, [\pref{Sylvester3}], had derived full reciprocity
somewhat earlier and would also have had the reciprocity of Dedekind sums, if he had
defined them.

It is instructive to make the expansion (\peq{wave4}) reflect the parity more easily. It
can be rewritten (using [\pref{Norlund}] p.162, (18)),
   $$\eqalign{
  P_d&={1\over (d-1)!\prod_i d_i}B^{(d)}_{d-1}
  \big(\ol l+\sum d_i/2\mid{\bf d}\big)\cr
  &={1\over2^{d-1} (d-1)!\prod_i d_i}\sum_{k=0}^{d-1}\comb {d-1}k \,(2\ol l\,)^k\,
  D^{(d)}_{d-k-1}[{\bf d}]\cr
  }
  \eql{wave5}
  $$
in terms of the computable constants, $D^{(d)}_\nu[{\bf d}]$ which vanish for odd
$\nu$. (This is the statement of reciprocity for this term).

Again, the constants $D^{(d)}_{2\nu}[{\bf d}]$ are essentially the $A_\nu$
polynomials, [\pref{Hirzebruch}], and Hirzebruch helpfully gives  this time, the relation,
  $$
  D^{(d)}_{2\nu}[{\bf d}]={(2\nu)!\over2^{2\nu}}A_{\nu}\big(p_1,\ldots,p_\nu\big)\,,
  \quad 2\nu\le d\,,
  \eql{Areln}
  $$
and the $\ol l$ expansion is,
   $$\eqalign{
  P_d&={1\over\prod_id_i}\bigg({\ol l^{d-1}\over(d-1)!}A_0+
  {1\over16}{\ol l^{d-3}\over(d-3)!}A_1+
  {1\over2^8}{\ol l^{d-5}\over(d-5)!}A_2+\ldots\bigg)\,,
  }
  \eql{wave6}
  $$
which terminates.\mgn{WHEN?}

Hirzebruch lists a few of the $A$ polynomials as functions of the indeterminates, $p_i$,
which are, here, the elementary symmetric functions of the {\it squares} of the degrees,
[\pref{Hirzebruch}], \S{\bf1.3} which, to remind, are $d$ in number. Thus,
  $$
  A_1=-{2\over3}p_1,\,\,A_2={2\over45}\big(7p_1^2-4p_2\big),\,\,A_3
  =-{4\over3^3.5.7}\big(16p_3-44p_2p_1+31p_1^3\big)\,.
  \eql{Agenus}
  $$

I also reproduce some values of the $B$ and $D$ constants as given in [\pref{Norlund}],
p.167.
  $$\eqalign{
   B^{(n)}_0&=1\,,\quad B^{(n)}_1=-{1\over 2}\si_1\,,\quad
  B^{(n)}_2={1\over 6}s_2+{1\over2}\si_2\cr
  D^{(n)}_0&=1\,,\quad D^{(n)}_2=-{1\over 3}s_2\,,\quad
  D^{(n)}_4={7\over 15}s_4+{2\over3}\sum_{i<j}d_i^2\,d_j^2
  }
  \eql{NBD}
  $$
where the $s_p$ are the sums of the $p$--th powers of the degrees, $d_i$. These
expressions are easily shown to agree with those obtained via the Todd  or the
$A$--polynomials. Furthermore, they hold for all $n$ and the conclusion is that one can
relax the conditions in (\peq{todd}) and (\peq{Areln}). Thus the range of $2\nu$ in
(\peq{Areln}) can be extended beyond $d$ and, for example, from (\peq{Agenus}), one
gets $D^{(1)}_4=7d_1^4/15$ agreeing with (\peq{NBD}).\footnote{ More generally,
the coefficient of the term $p_1^\nu$ in $A_\nu$ equals
$(2^{2\nu}/(2\nu)!)\,D_{2\nu}$ where the $D_{2\nu}$ are listed in Table 4 in
[\pref{Norlund1}] and are the coefficients in the expansion of the characteristic series
$x/\sinh x$.}

It should be apparent that there is no need to introduce generalised Bernoulli
polynomials. It is possible to move to the final formulae in (\peq{wave4}) and
(\peq{wave6}) directly in terms of multiplicative sequences. In fact this is what
Sylvester does in his explicit determinations.

The polynomial part, $P$ is that Sylvester wave corresponding to the root of unity, 1. He
denotes it by $W_1$, and in his series of papers on partitions explains, several times, a
systematic computational scheme for its algebraic determination, the most detailed being
[\pref{Sylvester3}]. There he obtains, by rapid direct calculation, an expansion identical
to ({\peq{wave6}), except that the polynomial coefficients are essentially all the
homogenous products (of a given order) expressed in terms of functions of the sums of
even {\it powers} of the degrees. I review the calculational details. A good introduction
to Sylvester's theory can be found in the textbook by Netto, [\pref{Netto}].

Instead of (\peq{sylv}), one now has the more symmetrical form,

  $$
  W_1= {\rm co}_{-1}\,{e^{\,\ol lt}\over\prod_{i=1}^d\big(e^{d_it/2}-e^{-d_it/2}\big)}
  =\,{\rm co}_{-1}\,{e^{\,\ol lt}\over\prod_{i=1}^d2\sinh{1\over2}td_i}\,.
   \eql{sylv}
  $$

Expansion of the numerator yields the coefficient of the power $\ol l^n/n!$, as
   $$\eqalign{
   &{\rm co}_{-1}{t^n\over\prod_i2\sinh{1\over2}td_i}
   ={1\over\prod d_i}\, {\rm co}_{-1}t^{n-d}\prod_i{x_i\over\sinh x_i}\cr
   &= {1\over\prod d_i}\, {\rm co}_{d-1-n}\prod_i{x_i\over\sinh x_i}
   }
   \eql{coeff2}
   $$
where I have set $x_i=td_i/2$ and have encountered a multiplicative sequence, $K$,
with characteristic function $Q(x)=x/\sinh x$ so that
  $$
  Q(x_1)Q(x_2)\ldots = 1+K_2(s_1,s_2)+K_4(s_1,s_2,s_3,s_4)+\ldots
  \eql{mults}
  $$
where $s_k$ is the sum of the $k$th powers of the $x_i$. It is more conventional to use
the symmetrical products, but this choice is more natural, as will be seen. The $K_m$
are homogeneous of degree $m$ in the $x_i$ and so, in view of the definition of the
$x_i$, the right--hand side of (\peq{mults}) is a power series in $t^2$. According to
(\peq{coeff2}) the required coefficient is simply $K_{d-1-n}$. To actually compute this,
it is necessary to expand the particular characteristic function as a power series. Although
everything is classical, I will do this {\it ab initio} from the old product
  $$
   Q(x)={x\over\sinh x}=\prod_{n=1}^\infty{1\over \big(1+x^2/n^2\pi^2\big)}\,.
   \eql{qprod}
  $$
Therefore
  $$\eqalign{
  \log Q(x)&=-\sum_{n=1}^\infty \log\big(1+x^2/n^2\pi^2\big)\cr
  &=-{\ze(2)\over\pi^2}\,x^2+{\ze(4)\over2\pi^4}\,x^4-\ldots
   +(-1)^k{\ze(2k)\over k\,\pi^{2k}}x^{2k}\pm\ldots\,.
  }
  \eql{logq}
  $$

The expansion of $\log Q$ rather than $Q$ is the tactical point here and the construction
of the multiplicative sequence, (\peq{mults}), follows on exponentiation,
  $$\eqalign{
  Q_1\,Q_2\,\ldots &=\exp\bigg(-{\ze(2)\over(2\pi)^2}\,s_2\,t^2
  +{\ze(4)\over2(2\pi)^4}\,s_4\,t^4-\ldots
   +{(-1)^k\ze(2k)\over k\,(2\pi)^{2k}}\,s_{2k}\,t^{2k}\pm\ldots\bigg)\cr
   &\equiv\exp\bigg(-\tau_1\,t^2
  +{1\over2}\tau_2\,t^4-\ldots
   +(-1)^k{1\over k}\tau_k\,t^{2k}\pm\ldots\bigg)
   }
   \eql{prodqs}
  $$
I have introduced Sylvester's notation,
  $$
  \tau_k={\ze(2k)\over (2\pi)^{2k}}s_{2k}
  ={(-1)^{k+1}\over2(2k)!}(-2k)\ze\big(-(2k-1)\big)\,s_{2k}=
  {(-1)^{k+1}\over2(2k)!}\,B_{2k}\,s_{2k}
  \eql{taus}
  $$
on the modern definition of the Bernoulli numbers. In some ways it is better to leave
these, and similar, constants in terms of \zf\ values.

The final step in determining the multiplicative sequence is to expand the right--hand
side of (\peq{prodqs}) as a power series in $t^2$. This can be  done elegantly by
relating this to the generating function of homogeneous products which is, (\eg\
Littlewood, [\pref{Littlewood2}]),
   $$
   1-h_1x+h_2x^2-\ldots=\prod_i{1\over1+x_i\,x}
   $$
where $h_r$ is the sum of all homogeneous products of degree $r$ of the
indeterminates, $x_i$. Taking logs, expanding that on the right--hand side, performing
the sum on $i$ and then exponentiating, to get back to the left--hand side, gives
  $$
  1-h_1x+h_2x^2-\ldots=\exp\big(-s_1\,x+{1\over2}s_2\,x^2-\ldots\big)
  \eql{genh}
  $$
and the $h_r$ are to be considered as functions of the sums of powers, $s_q$, of the
$x_i$. Comparison with (\peq{prodqs}) shows that
  $$
   Q_1\,Q_2\,\ldots =1-H_1\,t^2+H_2\,t^4-H_3\,t^6+\ldots +(-1)^rH_r\,t^{2r}+\ldots
   \eql{multsb}
  $$
where the $H_r$ are the same functions of the $\tau_i$ as the $h_r$ are of the $s_i$. (A
separate symbol is therefore not really needed.)

Solving the recursions that follow from the logarithmic derivative of (\peq{genh}) gives
Brioschis's determinantal expression for the $H_r$,  (Fa\`{a} de Bruno, [\pref{FdeB}]
\S 89, [\pref{Littlewood2}]),
  $$
   H_r={1\over r!}\,\left|\matrix{\tau_1&-1&0&0&0\cr
                                                 \tau_2&\tau_1&-2&0&0\cr
                                                 \tau_3&\tau_2&\tau_1&-3&0\cr
                                                 \ldots\cr
                                                 \tau_r&\tau_{r-1}&\ldots\cr
   }\right|\,.
   \eql{Brioschi}
  $$

Sylvester's final answer is,
  $$\eqalign{
  W_1&={1\over\prod_id_i}\bigg({\ol l^{d-1}\over(d-1)!}H_0-{\ol l^{d-3}\over(d-3)!}H_1+
  {\ol l^{d-5}\over(d-5)!}H_2+\ldots\bigg)\,,
  }
  \eql{wave7}
  $$
which is the same as (\peq{wave6}) expressed slightly differently ($H_0=1$).
 I observe that Sylvester was led naturally to multiplicative sequences.

 Sylvester, [\pref{Sylvester3}], lists six explicit (numerical) polynomials, which are also given by
 Ehrhart, [\pref{Ehrhart}].

\section{\bf6. The second wave}

There is a wave for every root of unity, I discuss in this paper only that one
corresponding to the root $-1$, denoted by $W_2$, which, like $W_1$, can be
algebraically expressed. Sylvester, [\pref{Sylvester2}], provides the following algebraic
statement: $W_2$ equals the coefficient of $1/t$ in the generating function,
  $$
 {(-1)^l e^{lt}\over \prod_i (1-e^{-\al_it})\,\prod_j(1+e^{-\be_jt})}
 \eql{genfun3}
  $$
where the $\al_i$ are the even degrees, and the $\be_j$ are the odd ones.

In order to express this in a form similar to (\peq{sylv2}), I recall the expansion formula
for the generalised Eulerian functions $E^{(n)}_\nu$, \eg\ [\pref{Norlund,Norlund1}],
  $$
  {2^ne^{xt}\over\prod_{j=1}^n (e^{\be_jt}+1)}=\sum_{\nu=0}^\infty {t^\nu\over\nu!}
  E^{(n)}_\nu(x\mid{\bbe})
  $$
to which I add, for handiness, the one for Bernoulli functions, already used in
(\peq{sylv2}),
  $$
  {t^n\prod_i\al_i\,\, e^{xt}\over\prod_{i=1}^n (e^{\al_it}-1)}
  =\sum_{\nu=0}^\infty {t^\nu\over\nu!}
  B^{(n)}_\nu(x\mid{\bal})\,.
  $$

Constructing the product in (\peq{genfun3}) the condition for the removal of the
exponential is
  $$
   \ol l=x+y-{1\over2}\sum_{i=1}^\al\al_i-{1\over2}\sum_{j=1}^\be\be_j\,,
   \quad (\al+\be=d)\,,
   \eql{cond3}
  $$
which allows two simplest (in terms of $\ol l$) solutions,
  $$
  x=\ol l+{1\over2}\sum\al_i\,\quad y={1\over2}\sum\be_j
  $$
and
  $$
  x={1\over2}\sum\al_i\,\quad y=\ol l+{1\over2}\sum\be_j\,,
  $$
the first of which yields the relevant coefficient required by Sylvester's rule as,
  $$
  W_2={(-1)^{l}\over2^\be(\al-1)!\prod\al_i}\sum_{\nu=0}^{\al-1}\comb{\al-1}\nu
  B^{(\al)}_\nu\big(\ol l+{\textstyle \sum}\al_i/2\mid\bal\big)
  \,E^{(\be)}_{\al-1-\nu}[\bbe]\,,
  \eql{w21}
  $$
while the second form leads to,
  $$
  W_2={(-1)^{l}\over2^{d-1}(\al-1)!\prod\al_i}\sum_{\nu=0}^{\al-1}\comb{\al-1}\nu 2^\nu
  D^{(\al)}_{\al-\nu-1}[\bal]\,E^{(\be)}_{\nu}
  \big(\ol l+{\textstyle \sum}\be_j/2\mid\bbe\big)\,.
  \eql{w22}
  $$
The constants $E^{(n)}_{\nu}[\bbe] \equiv E^{(n)}_{\nu}\big({\textstyle
\sum}\be_i/2\mid\bbe\big)$ are zero for odd $\nu$, [\pref{Norlund}].

Either form shows that $W_2$ is a 2--periodic function of $l$, through $(-1)^l$,
multiplied by a polynomial in $\ol l$. For an odd (even) number, $\al$, of even degrees,
the polynomial has even (odd) order. If there is only one even degree, this polynomial is
a constant, while $W_2$ is zero if all degrees are odd. (No degree has a factor of two.)

It is possible to check reciprocity for $W_2$ directly. Changing the sign of $\ol l$ in the
summand gives a factor of $(-1)^\nu$ which, because of the restriction on the index of
the $E$, equals $(-1)^{\al-1}$ which can be taken outside the sum. Also the sign
$(-1)^l$ changes to
$(-1)^{-l-\sum\al_i-\sum\be_j}=(-1)^{l-\sum\be_j}=(-1)^{l-\be}$ so that the total
change is to give a factor of $(-1)^{l-\be-\al+1}=(-1)^{l-d+1}$ in place of $(-1)^l$,
which is the required reciprocity.\mgn{MORE. SYLVESTER}

Just for comparison, I give the formula that follows from  an `$l$--solution',
  $$
  x= l+\sum\al_i\,\quad y=\sum\be_j
  $$
of (\peq{cond3}), {\it viz},
  $$
  W_2={(-1)^{l}\over2^\be(\al-1)!\prod\al_i}\sum_{\nu=0}^{\al-1}\comb{\al-1}\nu
  B^{(\al)}_\nu\big(l+{\textstyle \sum}\al_i\mid\bal\big)
  \,E^{(\be)}_{\al-1-\nu}\big({\textstyle \sum}\be_j\mid\bbe\big)\,.
  \eql{w23}
  $$

This form is given in the interesting paper by Rubinstein and Fel, [\pref{RandF}]. It is,
perhaps, less natural, reciprocity being harder to spot. Rubinstein and Fel also write the
other waves in a similar way more complicatedly using the generalised Eulerian
polynomials defined by Carlitz, [\pref{Carlitz}].

I now derive a much neater expression for $W_2$.  As might be expected, Sylvester
provides a route to an explicit solution, which I describe more compactly. I therefore
begin again and am required to compute the coefficient,
  $$
  W_2=(-1)^l {\rm co}_{-1}
  {e^{\ol l t}\over\prod_i 2\sinh{1\over2}\al_it\,\,\prod_j 2\cosh{1\over2}\be_jt}\,.
  $$

Similar to the $W_1$ computation, expansion of the numerator gives the coefficient of
the power $(-1)^l\,\ol l^n/n!$ as,
  $$\eqalign{
   &{\rm co}_{-1}{t^n\over\prod_i2\sinh{1\over2}t\al_i \prod_j2\cosh{1\over2}t\be_j}\cr
   &= {1\over2^\be\prod \al_i}\, {\rm co}_{\al-1-n}\prod_i{\xi_i\over\sinh \xi_i}
   \prod_j{1\over\cosh \eta_j}\,,
   }
   \eql{coeff8}
   $$
where $\xi_i=t\al_i/2$ and $\eta_j=t\be_j/2$. The product of two multiplicative series
now arises. One has been encountered previously, see (\peq{coeff2}), (\peq{mults}),
and I refer to it as the untwisted series. I now cover the other, defined by the
characteristic function,
  $$
  \Om(y)={1\over\cosh y}=\prod_{n=1}^\infty
  \bigg(1+{y^2\over\big((n-1/2)\pi\big)^2}\bigg)^{-1}\,,
  \eql{com}
  $$
in like manner\footnote{ The products (\peq{com})  and (\peq{qprod}) can be
considered as examples of the Mittag-Leffler theorem which notion allows one to
generalise the approach.}, and call this the twisted series. Thus,
  $$\eqalign{
  \log \Om(y)=-{\ze_R(2,1/2)\over\pi^2}\,y^2+{\ze_R(4,1/2)\over2\pi^4}\,y^4-\ldots
   +(-1)^k{\ze_R(2k,1/2)\over k\,\pi^{2k}}y^{2k}\pm\ldots\,,
  }
  $$
in terms of the Riemann-Hurwitz \zf, $\ze_R(s,w)$.

The similarity of this to (\peq{logq}) shows that the formal algebra is just the same, only
the numerical coefficients differ, with $B_k$ replaced by $(2^{2k}-1)B_k$. So I define,
  $$
  \varsigma_k={\ze_R(2k,1/2)\over (2\pi)^{2k}}s_{2k}=
  {(-1)^{k+1}(2^{2k}-1)\over2(2k)!}\,B_{2k}\,s_{2k}\,,
  \eql{sigmas}
  $$
and the multiplicative sequence this time is (\cf (\peq{multsb})),
  $$
   \Om_1\,\Om_2\,\ldots =1-H_1(\varsigma)\,t^2+H_2(\varsigma)\,t^4-H_3(\varsigma)
   \,t^6+\ldots +(-1)^rH_r(\varsigma)\,t^{2r}+\ldots\,.
   \eql{multsf}
  $$

According to (\peq{coeff8}) a product of multiplicative sequences is required. This can
be done from the series (\peq{multsb}) and (\peq{multsf}) but is best performed first
at the exponential level using  (\peq{prodqs}) and,
  $$\eqalign{
  \Om_1\,\Om_2\,\ldots
   =\exp\bigg(-\varsigma_1\,t^2
  +{1\over2}\varsigma_2\,t^4-\ldots
   +(-1)^k{1\over k}\varsigma_k\,t^{2k}\pm\ldots\bigg)\,.
   }
   \eql{prodoms}
  $$
Doing this, and then expanding I find, as before,
  $$
   Q_1Q_2\ldots\Om_1\,\Om_2\,\ldots =1-H_1(\varsigma+\tau)\,t^2+H_2(\varsigma+\tau)\,t^4
   -\ldots +(-1)^rH_r(\varsigma+\tau)\,t^{2r}+\ldots\,.
   \eql{multsf}
  $$

The final answer for the second wave is then formally similar to that for the first as a
terminating series,
$$\eqalign{
  W_2&=(-1)^l{1\over2^\be\prod_j\al_i}
  \bigg({\ol l^{\al-1}\over(\al-1)!}H_0-{\ol l^{\al-3}\over(\al-3)!}H_1+
  {\ol l^{\al-5}\over(\al-5)!}H_2+\ldots\bigg)\,,
  }
  \eql{wave10}
  $$
where $H_i\equiv H_i(\tau+\varsigma)$. $\tau_i$ defined by (\peq{taus}) with
$s_{2q}$ the sums of powers of the even degrees and $\varsigma_i$ is defined by
(\peq{sigmas}) with, this time, $s_{2q}$ being the sums of powers of the odd degrees.
I could not find this exact form in Sylvester's writings.

From a formal point of view, I note that the functions associated with the twisted
multiple sequence are the generalised Euler polynomials.

From a {\it spectral} perspective, the expansion variable is the eigenvalue, $\om$,
rather than $l$, hence I require
   $$\eqalign{
   (-1)^{d-1}B^{(d)}_{d-1}\big(a-\om\mid{\bf d}\big)&=
   (-1)^{d-1}\sum_{s=0}^{d-1}\comb {d-1}s(-\om)^s
   B^{(d)}_{d-1-s}(a\mid{\bf d})\,.\cr
   }
   \eql{sylv5}
   $$
So that the polynomial, $P$, is
  $$
  P=\sum_{k=0}^{d-1} P_k\,\om^k
  $$
with
  $$\eqalign{
  P_k&={(-1)^{d-1-k}\over(d-1-k)!\,k!\prod_i d_i} B^{(d)}_{d-1-k}(a\mid{\bf d})\cr
  &={2\over\Ga((k+1)/2)}\,C_{(d-k-1)/2}
  }
  $$
in terms of the heat--kernel coefficients, $C_*$, on S$^d/\Ga$ computed in
[\pref{ChandD}].\mgn{CF FORMULA fOR $\ol N$}

\begin{ignore}
For comparison, I write out the general series for the smoothed counting function
  $$
\ol N(\om)=\sum{C\over\Ga()}\om
  $$

It should be said that, as a strategic point, when calculating and employing spectral
objects such as functional determinants or heat--kernel expansions, explicit use of the
degeneracies can often be bypassed by working directly with the generating functions or
with the relevant zeta function, which, in this case, is the Barnes zeta function (see the
eigenvalues, (\peq{eigvals})).
\end{ignore}
\section{\bf 7. Ehrhart quasipolynomials}

Regarding the Ehrhart polynomial, I note that this is what could be termed an
`accumulated degeneracy', defined by,
  $$
  G(l)=\sum_{l'=0}^l g(l')\,,
  \eql{accgen}
  $$
generically expressed. It is commonplace in generating function circles to write the
obvious recursion for $G$ as
  $$
  H(\si)={1\over1-\si}\,h(\si)
  \eql{recurs2}
  $$
for $H(\si)\equiv\sum_{l=0}^\infty G(l)\,\si^l$  \etc

Exactly the present situation of spectra on regular tesselations of the $d$--sphere has
already been considered in [\pref{dowtess2}] (extended to $p$--forms) and some
discussion of the asymptotics were there made, using generating functions.

Looking back at say (\peq{genfun}), the recursion division in (\peq{recurs2})  is
tantamount to adding $1$ to list of degrees, $d_i$, and increasing $d$ by one.\footnote{
This is seen by noting that adding $1 m_3$ to the left--hand side of (\peq{cond}) is
equivalent to adding the results for (\peq{cond}) with the right--hand side equalling
$l,l-1,l-2,\ldots,0$ in turn. See [\pref{Sylvester3}], [\pref{Ehrhart}], [\pref{BandR}].}
The denumerant for this case is the Ehrhart quasipolynomial $L\big(\ol\caP_{\bf
d},l\big)$, \ie
  $$
  L\big(\ol\caP_{\bf d},l\big)={l\over{\bf d},1,}\,.
  $$

The strictly polynomial part, $P$, of the Ehrhart quasipolynomial is given by the
coefficient of $1/t$ in,
  $$
   {e^{\,lt}\over\prod_{i=1}^d\big(1-e^{-d_it}\big)\big(1-e^{-t}\big)}\,,
  $$
which has already been found in (\peq{sylv2}) as  (with $a=(d-1)/2$),
  $$\eqalign{
  P&={(-1)^d\over d!\prod_i d_i}B^{(d+1)}_d\big(a-\om\mid{\bf d},1\big)\cr
  &={(-1)^d\over d!\prod_i d_i}\sum_{k=0}^d\comb dk (-\om)^k
  B^{(d+1)}_{d-k}\big(a\mid{\bf d},1\big)\,.\cr
  }
  \eql{wave1}
  $$

\begin{ignore}
Standard manipulative rules for the Bernoulli polynomials allow the unit degree to be
removed at the expense of a series,
  $$
   B^{(d+1)}_{d-k}\big(a\mid{\bf d},1\big)=\sum_{j=0}^{d-k}
  \comb {d-k}j B_j B^{(d)}_{d-k-j}\big(a\mid{\bf d}\big)\,.
  \eql{wave2}
  $$
\end{ignore}

The relation with Todd polynomials gives
  $$
  B^{(d+1)}_{r}(x\mid{\bf d},1)=(-1)^r\,r!
  \sum_{s=0}^{r}(-1)^s{x^s\over s!}T_{r-s}(\ol\si_1,\ldots,\ol\si_{r-s})\,,\quad r\le d\,,
  \eql{todd2}
  $$
where now the $\ol\si_s$ are symmetric functions of the degrees $d_i$ ($i=1,\ldots,d$)
and 1.

First, I wish to express the symmetric functions, $\ol\si_s$, of $d_i$, $i=1,\ldots,d$ and
1 in terms those, $\si_s$, of just the $d_i$ since these are the variables used by
[\pref{BandR}], although this is not essential. Trivially,  or from the defining
fundamental identities,
  $$\eqalign{
 \prod_{i=1}^d(1+d_it)&=\sum_0^{d}\si_rt^r\cr
 \prod_{i=1}^{d}(1+d_it)(1+t)&=\sum_0^{d}\ol\si_rt^r\,,
 }
  $$
it  follows that
  $$
  \ol\si_s=\si_s+\si_{s-1}\,.
  \eql{recurs3}
  $$

For convenience I list a few resulting Todd polynomials obtained from the usual ones,
  $$
  \ol T_0=1\,,\quad \ol T_1={1\over2}(\si_1+1)\,,\quad \ol T_2={1\over12}(\si_1^2+\si_2+
  3\si_1+1)
  \eql{toddp}
  $$
and write down some Bernoulli functions from (\peq{todd}),
  $$\eqalign{
   B^{(d+1)}_1((d-1)/2\mid{\bf d},1)&={d-1\over2}\ol T_0-\ol T_1\cr
   B^{(d+1)}_2((d-1)/2\mid{\bf d},1)&=2\ol T_2
   -(d-1)\ol T_1+\bigg({d-1\over2}\bigg)^2\ol T_0
   }
  $$

The first two members of the polynomial (\peq{wave1}) now read
  $$\eqalign{
   P&={(-1)^d\over d!\prod_i d_i}\bigg(\om^d
   +{d(2+\si_1-d)\over2}\om^{d-1}+\ldots\bigg)\cr
  &={(-1)^d\over d!\prod_i d_i}\bigg(l^d
   +{d(\si_1+1)\over2}l^{d-1}+\ldots\bigg)
  }
  \eql{series2}
  $$
where, for the combinatorialists, and as a check, I have transformed the series to one in
$l$. These forms should be compared with (\peq{ehr3}) and (\peq{ehr4}).

Of course it is not necessary to go via $\om$ to $l$, one can write (no different to
(\peq{wave1})),
  $$\eqalign{
  P&={(-1)^d\over d!\prod_i d_i}B^{(d+1)}_d\big(l\mid{\bf d},1\big)\cr
  &={(-1)^d\over d!\prod_i d_i}\sum_{k=0}^d\comb dk (-l)^k
  B^{(d+1)}_{d-k}\big({\bf d},1\big)\,.\cr
  &={1\over\prod_id_i}\sum_{k=0}^d\,{l^k\over k!}\,T_{d-k}\big(\ol\si_1,\ldots,\ol\si_{d-k}\big)
  }
  \eql{wave3}
  $$
which is really only a re--expression. Using the forms of the Todd polynomials,
(\peq{toddp}), published coefficients are readily confirmed in rapid fashion.

However, rather than  use these expressions, it is far easier to  take advantage of
Sylvester's explicit expansion, (\peq{wave7}) and, making the adjustments to move to
the Ehrhart polynomial, I get for its polynomial part,
  $$\eqalign{
  P&={1\over\prod_id_i}\bigg({\ol l^{d}\over d!}H_0-{\ol l^{d-2}\over(d-1)!}H_1+
  {\ol l^{d-4}\over(d-2)!}H_2+\ldots\bigg)\,,
  }
  \eql{wave11}
  $$
where the $H_i$ are the same functions as before, (\peq{Brioschi}),  but of {\it new}
$\tau$s defined by,
  $$
  \tau_k={\ze(2k)\over (2\pi)^{2k}}\big(s_{2k}+1\big)=
  {(-1)^{k+1}B_{2k}\over2(2k)!}\,\,\big(s_{2k}+1\big)\,,
  \eql{ntaus}
  $$
obtained by adding 1 to the list of degrees. The $s_q$ are sums of powers of the $d_i$.

Particular cases are easily constructed and the general one machine coded.
\begin{ignore}

 Rather than carry on with this development using the degrees, $d_i$, I prefer to
express everything in terms of the symmetric functions of the `exponents' $m_i=d_i-1$
and 1, since these have a more geometric significance as the Shepherd--Todd numbers,
$b_r$.

I have discussed in [\pref{Dow20}] the conversion of between these numbers based on
the fundamental identity,
  $$\eqalign{
  \prod_{i=1}^d(1+m_it)&={1\over1+t}\sum_0^{d+1}b_rt^r\cr
 }
  $$
and found,
  $$
  \si_s=\sum_{r=0}^s(-1)^r\,b_r\sum_{i=r}^s(-1)^i\comb{d-i}{d-s}
  \eql{conv}
  $$

I set $m_i=d_i-1$ ($i=1,2,\ldots,d$), which are termed `exponents' in the context of
finite reflection groups, and conventionally extend their number by including the last one,
$m_{d+1}=1$.

The relevant set of exponents corresponding to the degrees $(1,{\bf d}) $ then is
$m_j=(m_0,m_1,\ldots,m_d,m_{d+1}$) with $m_0=0$ and $m_{d+1}=1$. I define
the symmetric functions of these exponents to be $\ol b_r$, where $r=0,1,\ldots,d+2$
and then
  $$\eqalign{
  \prod_{i=1}^d\,1\,(1+m_it)(1+t)&=\sum_0^{d+2}\ol b_rt^r\cr
 }
  $$
from which, $\ol b_{d+2}=0$ and $\ol b_r=b_r$, ($r=0,\ldots,d-1$). Because
$\ol\si_s$ and $\ol b_r$ are related as in (\peq{conv}) (with $d\to d+1$), these
conditions translate into the recursion (\peq{recurs3}).

The advantage of using the $b_r$ is that, suppose for a general set of, say, $(d+1)$
degrees a quantity has been found in terms of the $b_r$. Then, if one of the degrees is
changed to 1, the quantity will remain unchanged except the $b_r$ will be those
appropriate to the remaining degrees and also will equal the same quantity in $d$
dimensions, if $d+1$ is replaced by $d$. This allows me to compress the calculations and
use known results. For example, the expressions [\pref{Dow20}],
  $$\eqalign{
  B^{(d)}_2\big((d-1)/2\mid{\bf d}\big)&={1\over6}\big(b_1(b_1-1)+b_2+(3-d)/2\big)\cr
  B^{(d)}_3\big((d-1)/2\mid{\bf d}\big)&={1\over4}b_1\big(b_1-b_2-(3-d)/2\big)
  }
  $$
were derived using the Todd polynomials and the conversion, (\peq{conv}). This means,
in fact, that the  expressions are independent of the upper index $(d)$ on the Bernoulli
functions, so long as this does not exceed the lower index.
\end{ignore}

I can now return to the recipe given in section 2 and construct the combination of
Ehrhart polynomials,
  $$
  {1\over2}\big[A(\om)+A(\om-1)\big]=
  {(-1)^d\over d!\prod_i d_i}\bigg(\om^d
   +{d(\si_1-d+1)\over2}\om^{d-1}+\ldots\bigg)
  $$
which agrees with the first two terms of the asymptotic expansion, (\peq{nasym4}), of
the smoothed counting function, since $b_1=\si_1-d+1$ and confirms the statement
made in section 2, for the first two terms.

\section{\bf 8. Conclusion}

It has been shown that the first two terms in the orthodox asymptotic expansion of the
smoothed counting function can be obtained from a symmetrical combination of two
Ehrhart polynomials in any dimension.

It is, perhaps, no surprise that the heat--kernel coefficients can be obtained from the
Ehrhart polynomial as this encodes all the eigenvalue information.

The technical evaluation of denumerants is greatly eased by the use of multiplicative
sequences. Indeed these arise naturally. A simple expression for the second wave
involving the numbers of all homogeneous products. In a later paper I intend to look at
waves beyond the second one.

\newpage
 \vglue 20truept

 \noin{\bf References.} \vskip5truept

\begin{putreferences}
    \ref{Littlewood2}{Littlewood,D.E. {\it The Theory of Group Characters}
    (Clarendon Press, Oxford, 1950).}
    \ref{Wright}{Wright,E.M. \amm{68}{1961}{144}.}
\ref{Carlitz}{Carlitz,L. \dmj{27}{1960}{401}.}
     \ref{Netto}{Netto,E. {\it Lehrbuch der Combinatorik} 2nd Edn. (Teubner, Leipzig, 1927).}
    \ref{FdeB}{Fa\`{a} de Bruno, F. {\it Th\'eorie des Formes Binaires} (Brero, Turin,1876).}
    \ref{Ehrhart}{Ehrhart,E. \jram{227}{25}{1967}.}
    \ref{Bell}{Bell,E.T.\ajm{65}{1943}{382}.}
    \ref{BandR}{Beck,M. and Robins,S. {\it Computing the Continuous Discretely,}
    (Springer, New York, 2007).}
    \ref{BandR2}{Beck, M. and Robins,S. {\it Discrete and Comp. Geom.} {\bf 27}(2002) 443.}
    \ref{Harmer}{Harmer,M. {\it J.Australian Math.Soc.} {\bf 84}(2008)217.}
    \ref{RandF}{Rubinstein, B.Y. and Fel,L.G., {\it Ramanujan J.} {\bf11}(2006)331.}
    \ref{BGK}{Beck,M., Gessel, I.M. and Komatsu,T. {\it Electronic Journal of Combinatorics}
    {\bf8}(2001) 1.}
    \ref{Sylvester}{Sylvester,J.J. \qjpam{1}{1858}{81}.}
    \ref{Sylvester2}{Sylvester,J.J. \qjpam{1}{1858}{142}.}
    \ref{Sylvester3}{Sylvester,J.J. \ajm{5}{1882}{79}.}
    \ref{Sylvester4}{Sylvester,J.J. \plms{28}{1896}{33}.}
    \ref{Dowgta}{J.S.Dowker, {\it Group theory aspects of spectral problems on spherical factors},
   ArXiv.Math.DG: 0907.1309.}
    \ref{BDR}{Beck,M., Diaz and Robins,S. {\it J.Numb.Theory} {\bf 96} (2002) 1.}
    \ref{PandS}{P\'{o}lya, G. and Szeg\H{o},G. {\it Aufgaben und Lehrs\"atze aus der Analysis}
    (Springer--Verlag, Berlin, 1925).}
    \ref{EOS}{Elizalde,E., Odintsov, S.D. and Saharian, A.A. \prD{79}{2009}{065023}.}
    \ref{Cavalcanti}{Cavalcanti,R.M. \prD{69}{2004}{065015}.}
    \ref{MWK}{Milton, K.A., Wagner,J. and Kirsten,K. \prD{80}{2009}{125028}.}
    \ref{EBM2}{Ellingsen,S.A., Brevik,I. and Milton,K.A. \prE{81}{2010}{065031}.}
    \ref{EBM}{Ellingsen,S.A., Brevik,I. and Milton,K.A. \prE{80}{2009}{021125}.}
    \ref{BEM}{Brevik,I., Ellingsen,S.A. and Milton,K.A. \prE{79}{2009}{041120}.}
    \ref{FKW}{Fulling,S.A, Kaplan L. and Wilson,J.H. \prA{76}{2007}{012118}.}
    \ref{Lukosz}{Lukosz,W, {\it Physica} {\bf 56} (1971) 109; \zfp{258}{1973}{99}
    ;\zfp{262}{1973}{327}.}
    \ref{Gromes}{Gromes, D. \mz{94}{1966}{110}.}
    \ref{FandK1}{Kirsten,K. and Fulling,S.A. \prD{79}{2009}{065019} .}
    \ref{FandK2}{Fucci,G. and Kirsten,K, JHEP (2011), 1103:016.}
    \ref{dowgjms}{Dowker,J.S. {\it Determinants and conformal anomalies
    of GJMS operators on spheres}, ArXiv: 1007.3865.}
    \ref{Dowcascone}{dowker,J.S. \prD{36}{1987}{3095}.}
    \ref{Dowcos}{dowker,J.S. \prD{36}{1987}{3742}.}
    \ref{Dowtherm}{Dowker,J.S. \prD{18}{1978}{1856}.}
    \ref{Dowgeo}{Dowker,J.S. \cqg{11}{1994}{L55}.}
    \ref{ApandD2}{Dowker,J.S. and Apps,J.S. \cqg{12}{1995}{1363}.}
   \ref{HandW}{Hertzberg,M.P. and Wilczek,F. {\it Some calculable contributions to
   Entanglement Entropy}, ArXiv:1007.0993.}
   \ref{KandB}{Kamela,M. and Burgess,C.P. \cjp{77}{1999}{85}.}
   \ref{Dowhyp}{Dowker,J.S. \jpa{43}{2010}{445402}; ArXiv:1007.3865.}
   \ref{LNST}{Lohmayer,R., Neuberger,H, Schwimmer,A. and Theisen,S.
   \plb{685}{2010}{222}.}
   \ref{Allen2}{Allen,B. PhD Thesis, University of Cambridge, 1984.}
   \ref{MyandS}{Myers,R.C. and Sinha,A. {\it Seeing a c-theorem with
   holography}, ArXiv:1006.1263}
   \ref{MyandS2}{Myers,R.C. and Sinha,A. {\it Holographic c-theorems in
   arbitrary dimensions},\break ArXiv: 1011.5819.}
   \ref{RyandT}{Ryu,S. and Takayanagi,T. JHEP {\bf 0608}(2006)045.}
   \ref{CaandH}{Casini,H. and Huerta,M. {\it Entanglement entropy
   for the n--sphere},\break arXiv:1007.1813.}
   \ref{CaandH3}{Casini,H. and Huerta,M. \jpa {42}{2009}{504007}.}
   \ref{Solodukhin}{Solodukhin,S.N. \plb{665}{2008}{305}.}
   \ref{Solodukhin2}{Solodukhin,S.N. \plb{693}{2010}{605}.}
   \ref{CaandW}{Callan,C.G. and Wilczek,F. \plb{333}{1994}{55}.}
   \ref{FandS1}{Fursaev,D.V. and Solodukhin,S.N. \plb{365}{1996}{51}.}
   \ref{FandS2}{Fursaev,D.V. and Solodukhin,S.N. \prD{52}{1995}{2133}.}
   \ref{Fursaev}{Fursaev,D.V. \plb{334}{1994}{53}.}
   \ref{Donnelly2}{Donnelly,H. \ma{224}{1976}{161}.}
   \ref{ApandD}{Apps,J.S. and Dowker,J.S. \cqg{15}{1998}{1121}.}
   \ref{FandM}{Fursaev,D.V. and Miele,G. \prD{49}{1994}{987}.}
   \ref{Dowker2}{Dowker,J.S.\cqg{11}{1994}{L137}.}
   \ref{Dowker1}{Dowker,J.S.\prD{50}{1994}{6369}.}
   \ref{FNT}{Fujita,M.,Nishioka,T. and Takayanagi,T. JHEP {\bf 0809}
   (2008) 016.}
   \ref{Hund}{Hund,F. \zfp{51}{1928}{1}.}
   \ref{Elert}{Elert,W. \zfp {51}{1928}{8}.}
   \ref{Poole2}{Poole,E.G.C. \qjm{3}{1932}{183}.}
   \ref{Bellon}{Bellon,M.P. {\it On the icosahedron: from two to three
   dimensions}, arXiv:0705.3241.}
   \ref{Bellon2}{Bellon,M.P. \cqg{23}{2006}{7029}.}
   \ref{McLellan}{McLellan,A,G. \jpc{7}{1974}{3326}.}
   \ref{Boiteaux}{Boiteaux, M. \jmp{23}{1982}{1311}.}
   \ref{HHandK}{Hage Hassan,M. and Kibler,M. {\it On Hurwitz
   transformations} in {Le probl\`eme de factorisation de Hurwitz}, Eds.,
   A.Ronveaux and D.Lambert (Fac.Univ.N.D. de la Paix, Namur, 1991),
   pp.1-29.}
   \ref{Weeks2}{Weeks,Jeffrey \cqg{23}{2006}{6971}.}
   \ref{LandW}{Lachi\`eze-Rey,M. and Weeks,Jeffrey, {\it Orbifold construction of
   the modes on the Poincar\'e dodecahedral space}, arXiv:0801.4232.}
   \ref{Cayley4}{Cayley,A. \qjpam{58}{1879}{280}.}
   \ref{JMS}{Jari\'c,M.V., Michel,L. and Sharp,R.T. {\it J.Physique}
   {\bf 45} (1984) 1. }
   \ref{AandB}{Altmann,S.L. and Bradley,C.J.  {\it Phil. Trans. Roy. Soc. Lond.}
   {\bf 255} (1963) 199.}
   \ref{CandP}{Cummins,C.J. and Patera,J. \jmp{29}{1988}{1736}.}
   \ref{Sloane}{Sloane,N.J.A. \amm{84}{1977}{82}.}
   \ref{Gordan2}{Gordan,P. \ma{12}{1877}{147}.}
   \ref{DandSh}{Desmier,P.E. and Sharp,R.T. \jmp{20}{1979}{74}.}
   \ref{Kramer}{Kramer,P., \jpa{38}{2005}{3517}.}
   \ref{Klein2}{Klein, F.\ma{9}{1875}{183}.}
   \ref{Hodgkinson}{Hodgkinson,J. \jlms{10}{1935}{221}.}
   \ref{ZandD}{Zheng,Y. and Doerschuk, P.C. {\it Acta Cryst.} {\bf A52}
   (1996) 221.}
   \ref{EPM}{Elcoro,L., Perez--Mato,J.M. and Madariaga,G.
   {\it Acta Cryst.} {\bf A50} (1994) 182.}
    \ref{PSW2}{Prandl,W., Schiebel,P. and Wulf,K.
   {\it Acta Cryst.} {\bf A52} (1999) 171.}
    \ref{FCD}{Fan,P--D., Chen,J--Q. and Draayer,J.P.
   {\it Acta Cryst.} {\bf A55} (1999) 871.}
   \ref{FCD2}{Fan,P--D., Chen,J--Q. and Draayer,J.P.
   {\it Acta Cryst.} {\bf A55} (1999) 1049.}
   \ref{Honl}{H\"onl,H. \zfp{89}{1934}{244}.}
   \ref{PSW}{Patera,J., Sharp,R.T. and Winternitz,P. \jmp{19}{1978}{2362}.}
   \ref{LandH}{Lohe,M.A. and Hurst,C.A. \jmp{12}{1971}{1882}.}
   \ref{RandSA}{Ronveaux,A. and Saint-Aubin,Y. \jmp{24}{1983}{1037}.}
   \ref{JandDeV}{Jonker,J.E. and De Vries,E. \npa{105}{1967}{621}.}
   \ref{Rowe}{Rowe, E.G.Peter. \jmp{19}{1978}{1962}.}
   \ref{KNR}{Kibler,M., N\'egadi,T. and Ronveaux,A. {\it The Kustaanheimo-Stiefel
   transformation and certain special functions} \lnm{1171}{1985}{497}.}
   \ref{GLP}{Gilkey,P.B., Leahy,J.V. and Park,J-H, \jpa{29}{1996}{5645}.}
   \ref{Kohler}{K\"ohler,K.: Equivariant Reidemeister torsion on
   symmetric spaces. Math.Ann. {\bf 307}, 57-69 (1997)}
   \ref{Kohler2}{K\"ohler,K.: Equivariant analytic torsion on ${\bf P^nC}$.
   Math.Ann.{\bf 297}, 553-565 (1993) }
   \ref{Kohler3}{K\"ohler,K.: Holomorphic analytic torsion on Hermitian
   symmetric spaces. J.Reine Angew.Math. {\bf 460}, 93-116 (1995)}
   \ref{Zagierzf}{Zagier,D. {\it Zetafunktionen und Quadratische
   K\"orper}, (Springer--Verlag, Berlin, 1981).}
   \ref{Stek}{Stekholschkik,R. {\it Notes on Coxeter transformations and the McKay
   correspondence.} (Springer, Berlin, 2008).}
   \ref{Pesce}{Pesce,H. \cmh {71}{1996}{243}.}
   \ref{Pesce2}{Pesce,H. {\it Contemp. Math} {\bf 173} (1994) 231.}
   \ref{Sutton}{Sutton,C.J. {\it Equivariant isospectrality
   and isospectral deformations on spherical orbifolds}, ArXiv:math/0608567.}
   \ref{Sunada}{Sunada,T. \aom{121}{1985}{169}.}
   \ref{GoandM}{Gornet,R, and McGowan,J. {\it J.Comp. and Math.}
   {\bf 9} (2006) 270.}
   \ref{Suter}{Suter,R. {\it Manusc.Math.} {\bf 122} (2007) 1-21.}
   \ref{Lomont}{Lomont,J.S. {\it Applications of finite groups} (Academic
   Press, New York, 1959).}
   \ref{DandCh2}{Dowker,J.S. and Chang,Peter {\it Analytic torsion on
   spherical factors and tessellations}, arXiv:math.DG/0904.0744 .}
   \ref{Mackey}{Mackey,G. {\it Induced representations}
   (Benjamin, New York, 1968).}
   \ref{Koca}{Koca, {\it Turkish J.Physics}.}
   \ref{Brylinski}{Brylinski, J-L., {\it A correspondence dual to McKay's}
    ArXiv alg-geom/9612003.}
   \ref{Rossmann}{Rossman,W. {\it McKay's correspondence
   and characters of finite subgroups of\break SU(2)} {\it Progress in Math.}
      Birkhauser  (to appear) .}
   \ref{JandL}{James, G. and Liebeck, M. {\it Representations and
   characters of groups} (CUP, Cambridge, 2001).}
   \ref{IandR}{Ito,Y. and Reid,M. {\it The Mckay correspondence for finite
   subgroups of SL(3,C)} Higher dimensional varieties, (Trento 1994),
   221-240, (Berlin, de Gruyter 1996).}
   \ref{BandF}{Bauer,W. and Furutani, K. {\it J.Geom. and Phys.} {\bf
   58} (2008) 64.}
   \ref{Luck}{L\"uck,W. \jdg{37}{1993}{263}.}
   \ref{LandR}{Lott,J. and Rothenberg,M. \jdg{34}{1991}{431}.}
   \ref{DoandKi} {Dowker.J.S. and Kirsten, K. {\it Analysis and Appl.}
   {\bf 3} (2005) 45.}
   \ref{dowtess1}{Dowker,J.S. \cqg{23}{2006}{1}.}
   \ref{dowtess2}{Dowker,J.S. {\it J.Geom. and Phys.} {\bf 57} (2007) 1505.}
   \ref{MHS}{De Melo,T., Hartmann,L. and Spreafico,M. {\it Reidemeister
   Torsion and analytic torsion of discs}, ArXiv:0811.3196.}
   \ref{Vertman}{Vertman, B. {\it Analytic Torsion of a  bounded
   generalized cone}, ArXiv:0808.0449.}
   \ref{WandY} {Weng,L. and You,Y., {\it Int.J. of Math.}{\bf 7} (1996)
   109.}
   \ref{ScandT}{Schwartz, A.S. and Tyupkin,Yu.S. \np{242}{1984}{436}.}
   \ref{AAR}{Andrews, G.E., Askey,R. and Roy,R. {\it Special functions}
  (CUP, Cambridge, 1999).}
   \ref{Tsuchiya}{Tsuchiya, N.: R-torsion and analytic torsion for spherical
   Clifford-Klein manifolds.: J. Fac.Sci., Tokyo Univ. Sect.1 A, Math.
   {\bf 23}, 289-295 (1976).}
   \ref{Tsuchiya2}{Tsuchiya, N. J. Fac.Sci., Tokyo Univ. Sect.1 A, Math.
   {\bf 23}, 289-295 (1976).}
  \ref{Lerch}{Lerch,M. \am{11}{1887}{19}.}
  \ref{Lerch2}{Lerch,M. \am{29}{1905}{333}.}
  \ref{TandS}{Threlfall, W. and Seifert, H. \ma{104}{1930}{1}.}
  \ref{RandS}{Ray, D.B., and Singer, I. \aim{7}{1971}{145}.}
  \ref{RandS2}{Ray, D.B., and Singer, I. {\it Proc.Symp.Pure Math.}
  {\bf 23} (1973) 167.}
  \ref{Jensen}{Jensen,J.L.W.V. \aom{17}{1915-1916}{124}.}
  \ref{Rosenberg}{Rosenberg, S. {\it The Laplacian on a Riemannian Manifold}
  (CUP, Cambridge, 1997).}
  \ref{Nando2}{Nash, C. and O'Connor, D-J. {\it Int.J.Mod.Phys.}
  {\bf A10} (1995) 1779.}
  \ref{Fock}{Fock,V. \zfp{98}{1935}{145}.}
  \ref{Levy}{Levy,M. \prs {204}{1950}{145}.}
  \ref{Schwinger2}{Schwinger,J. \jmp{5}{1964}{1606}.}
  \ref{Muller}{M\"uller, \lnm{}{}{}.}
  \ref{VMK}{Varshalovich.}
  \ref{DandWo}{Dowker,J.S. and Wolski, A. \prA{46}{1992}{6417}.}
  \ref{Zeitlin1}{Zeitlin,V. {\it Physica D} {\bf 49} (1991).  }
  \ref{Zeitlin0}{Zeitlin,V. {\it Nonlinear World} Ed by
   V.Baryakhtar {\it et al},  Vol.I p.717,  (World Scientific, Singapore, 1989).}
  \ref{Zeitlin2}{Zeitlin,V. \prl{93}{2004}{264501}. }
  \ref{Zeitlin3}{Zeitlin,V. \pla{339}{2005}{316}. }
  \ref{Groenewold}{Groenewold, H.J. {\it Physica} {\bf 12} (1946) 405.}
  \ref{Cohen}{Cohen, L. \jmp{7}{1966}{781}.}
  \ref{AandW}{Argawal G.S. and Wolf, E. \prD{2}{1970}{2161,2187,2206}.}
  \ref{Jantzen}{Jantzen,R.T. \jmp{19}{1978}{1163}.}
  \ref{Moses2}{Moses,H.E. \aop{42}{1967}{343}.}
  \ref{Carmeli}{Carmeli,M. \jmp{9}{1968}{1987}.}
  \ref{SHS}{Siemans,M., Hancock,J. and Siminovitch,D. {\it Solid State
  Nuclear Magnetic Resonance} {\bf 31}(2007)35.}
 \ref{Dowk}{Dowker,J.S. \prD{28}{1983}{3013}.}
 \ref{Heine}{Heine, E. {\it Handbuch der Kugelfunctionen}
  (G.Reimer, Berlin. 1878, 1881).}
  \ref{Pockels}{Pockels, F. {\it \"Uber die Differentialgleichung $\De
  u+k^2u=0$} (Teubner, Leipzig. 1891).}
  \ref{Hamermesh}{Hamermesh, M., {\it Group Theory} (Addison--Wesley,
  Reading. 1962).}
  \ref{Racah}{Racah, G. {\it Group Theory and Spectroscopy}
  (Princeton Lecture Notes, 1951). }
  \ref{Gourdin}{Gourdin, M. {\it Basics of Lie Groups} (Editions
  Fronti\'eres, Gif sur Yvette. 1982.)}
  \ref{Clifford}{Clifford, W.K. \plms{2}{1866}{116}.}
  \ref{Story2}{Story, W.E. \plms{23}{1892}{265}.}
  \ref{Story}{Story, W.E. \ma{41}{1893}{469}.}
  \ref{Poole}{Poole, E.G.C. \plms{33}{1932}{435}.}
  \ref{Dickson}{Dickson, L.E. {\it Algebraic Invariants} (Wiley, N.Y.
  1915).}
  \ref{Dickson2}{Dickson, L.E. {\it Modern Algebraic Theories}
  (Sanborn and Co., Boston. 1926).}
  \ref{Hilbert2}{Hilbert, D. {\it Theory of algebraic invariants} (C.U.P.,
  Cambridge. 1993).}
  \ref{Olver}{Olver, P.J. {\it Classical Invariant Theory} (C.U.P., Cambridge.
  1999.)}
  \ref{AST}{A\v{s}erova, R.M., Smirnov, J.F. and Tolsto\v{i}, V.N. {\it
  Teoret. Mat. Fyz.} {\bf 8} (1971) 255.}
  \ref{AandS}{A\v{s}erova, R.M., Smirnov, J.F. \np{4}{1968}{399}.}
  \ref{Shapiro}{Shapiro, J. \jmp{6}{1965}{1680}.}
  \ref{Shapiro2}{Shapiro, J.Y. \jmp{14}{1973}{1262}.}
  \ref{NandS}{Noz, M.E. and Shapiro, J.Y. \np{51}{1973}{309}.}
  \ref{Cayley2}{Cayley, A. {\it Phil. Trans. Roy. Soc. Lond.}
  {\bf 144} (1854) 244.}
  \ref{Cayley3}{Cayley, A. {\it Phil. Trans. Roy. Soc. Lond.}
  {\bf 146} (1856) 101.}
  \ref{Wigner}{Wigner, E.P. {\it Gruppentheorie} (Vieweg, Braunschweig. 1931).}
  \ref{Sharp}{Sharp, R.T. \ajop{28}{1960}{116}.}
  \ref{Laporte}{Laporte, O. {\it Z. f. Naturf.} {\bf 3a} (1948) 447.}
  \ref{Lowdin}{L\"owdin, P-O. \rmp{36}{1964}{966}.}
  \ref{Ansari}{Ansari, S.M.R. {\it Fort. d. Phys.} {\bf 15} (1967) 707.}
  \ref{SSJR}{Samal, P.K., Saha, R., Jain, P. and Ralston, J.P. {\it
  Testing Isotropy of Cosmic Microwave Background Radiation},
  astro-ph/0708.2816.}
  \ref{Lachieze}{Lachi\'eze-Rey, M. {\it Harmonic projection and
  multipole Vectors}. astro- \break ph/0409081.}
  \ref{CHS}{Copi, C.J., Huterer, D. and Starkman, G.D.
  \prD{70}{2003}{043515}.}
  \ref{Jaric}{Jari\'c, J.P. {\it Int. J. Eng. Sci.} {\bf 41} (2003) 2123.}
  \ref{RandD}{Roche, J.A. and Dowker, J.S. \jpa{1}{1968}{527}.}
  \ref{KandW}{Katz, G. and Weeks, J.R. \prD{70}{2004}{063527}.}
  \ref{Waerden}{van der Waerden, B.L. {\it Die Gruppen-theoretische
  Methode in der Quantenmechanik} (Springer, Berlin. 1932).}
  \ref{EMOT}{Erdelyi, A., Magnus, W., Oberhettinger, F. and Tricomi, F.G. {
  \it Higher Transcendental Functions} Vol.1 (McGraw-Hill, N.Y. 1953).}
  \ref{Dowzilch}{Dowker, J.S. {\it Proc. Phys. Soc.} {\bf 91} (1967) 28.}
  \ref{DandD}{Dowker, J.S. and Dowker, Y.P. {\it Proc. Phys. Soc.}
  {\bf 87} (1966) 65.}
  \ref{DandD2}{Dowker, J.S. and Dowker, Y.P. \prs{}{}{}.}
  \ref{Dowk3}{Dowker,J.S. \cqg{7}{1990}{1241}.}
  \ref{Dowk5}{Dowker,J.S. \cqg{7}{1990}{2353}.}
  \ref{CoandH}{Courant, R. and Hilbert, D. {\it Methoden der
  Mathematischen Physik} vol.1 \break (Springer, Berlin. 1931).}
  \ref{Applequist}{Applequist, J. \jpa{22}{1989}{4303}.}
  \ref{Torruella}{Torruella, \jmp{16}{1975}{1637}.}
  \ref{Weinberg}{Weinberg, S.W. \pr{133}{1964}{B1318}.}
  \ref{Meyerw}{Meyer, W.F. {\it Apolarit\"at und rationale Curven}
  (Fues, T\"ubingen. 1883.) }
  \ref{Ostrowski}{Ostrowski, A. {\it Jahrsb. Deutsch. Math. Verein.} {\bf
  33} (1923) 245.}
  \ref{Kramers}{Kramers, H.A. {\it Grundlagen der Quantenmechanik}, (Akad.
  Verlag., Leipzig, 1938).}
  \ref{ZandZ}{Zou, W.-N. and Zheng, Q.-S. \prs{459}{2003}{527}.}
  \ref{Weeks1}{Weeks, J.R. {\it Maxwell's multipole vectors
  and the CMB}.  astro-ph/0412231.}
  \ref{Corson}{Corson, E.M. {\it Tensors, Spinors and Relativistic Wave
  Equations} (Blackie, London. 1950).}
  \ref{Rosanes}{Rosanes, J. \jram{76}{1873}{312}.}
  \ref{Salmon}{Salmon, G. {\it Lessons Introductory to the Modern Higher
  Algebra} 3rd. edn. \break (Hodges,  Dublin. 1876.)}
  \ref{Milnew}{Milne, W.P. {\it Homogeneous Coordinates} (Arnold. London. 1910).}
  \ref{Niven}{Niven, W.D. {\it Phil. Trans. Roy. Soc.} {\bf 170} (1879) 393.}
  \ref{Scott}{Scott, C.A. {\it An Introductory Account of
  Certain Modern Ideas and Methods in Plane Analytical Geometry,}
  (MacMillan, N.Y. 1896).}
  \ref{Bargmann}{Bargmann, V. \rmp{34}{1962}{300}.}
  \ref{Maxwell}{Maxwell, J.C. {\it A Treatise on Electricity and
  Magnetism} 2nd. edn. (Clarendon Press, Oxford. 1882).}
  \ref{BandL}{Biedenharn, L.C. and Louck, J.D.
  {\it Angular Momentum in Quantum Physics} (Addison-Wesley, Reading. 1981).}
  \ref{Weylqm}{Weyl, H. {\it The Theory of Groups and Quantum Mechanics}
  (Methuen, London. 1931).}
  \ref{Robson}{Robson, A. {\it An Introduction to Analytical Geometry} Vol I
  (C.U.P., Cambridge. 1940.)}
  \ref{Sommerville}{Sommerville, D.M.Y. {\it Analytical Conics} 3rd. edn.
   (Bell, London. 1933).}
  \ref{Coolidge}{Coolidge, J.L. {\it A Treatise on Algebraic Plane Curves}
  (Clarendon Press, Oxford. 1931).}
  \ref{SandK}{Semple, G. and Kneebone. G.T. {\it Algebraic Projective
  Geometry} (Clarendon Press, Oxford. 1952).}
  \ref{AandC}{Abdesselam A., and Chipalkatti, J. {\it The Higher
  Transvectants are redundant}, arXiv:0801.1533 [math.AG] 2008.}
  \ref{Elliott}{Elliott, E.B. {\it The Algebra of Quantics} 2nd edn.
  (Clarendon Press, Oxford. 1913).}
  \ref{Elliott2}{Elliott, E.B. \qjpam{48}{1917}{372}.}
  \ref{Howe}{Howe, R. \tams{313}{1989}{539}.}
  \ref{Clebsch}{Clebsch, A. \jram{60}{1862}{343}.}
  \ref{Prasad}{Prasad, G. \ma{72}{1912}{136}.}
  \ref{Dougall}{Dougall, J. \pems{32}{1913}{30}.}
  \ref{Penrose}{Penrose, R. \aop{10}{1960}{171}.}
  \ref{Penrose2}{Penrose, R. \prs{273}{1965}{171}.}
  \ref{Burnside}{Burnside, W.S. \qjm{10}{1870}{211}. }
  \ref{Lindemann}{Lindemann, F. \ma{23} {1884}{111}.}
  \ref{Backus}{Backus, G. {\it Rev. Geophys. Space Phys.} {\bf 8} (1970) 633.}
  \ref{Baerheim}{Baerheim, R. {\it Q.J. Mech. appl. Math.} {\bf 51} (1998) 73.}
  \ref{Lense}{Lense, J. {\it Kugelfunktionen} (Akad.Verlag, Leipzig. 1950).}
  \ref{Littlewood}{Littlewood, D.E. \plms{50}{1948}{349}.}
  \ref{Fierz}{Fierz, M. {\it Helv. Phys. Acta} {\bf 12} (1938) 3.}
  \ref{Williams}{Williams, D.N. {\it Lectures in Theoretical Physics} Vol. VII,
  (Univ.Colorado Press, Boulder. 1965).}
  \ref{Dennis}{Dennis, M. \jpa{37}{2004}{9487}.}
  \ref{Pirani}{Pirani, F. {\it Brandeis Lecture Notes on
  General Relativity,} edited by S. Deser and K. Ford. (Brandeis, Mass. 1964).}
  \ref{Sturm}{Sturm, R. \jram{86}{1878}{116}.}
  \ref{Schlesinger}{Schlesinger, O. \ma{22}{1883}{521}.}
  \ref{Askwith}{Askwith, E.H. {\it Analytical Geometry of the Conic
  Sections} (A.\&C. Black, London. 1908).}
  \ref{Todd}{Todd, J.A. {\it Projective and Analytical Geometry}.
  (Pitman, London. 1946).}
  \ref{Glenn}{Glenn. O.E. {\it Theory of Invariants} (Ginn \& Co, N.Y. 1915).}
  \ref{DowkandG}{Dowker, J.S. and Goldstone, M. \prs{303}{1968}{381}.}
  \ref{Turnbull}{Turnbull, H.A. {\it The Theory of Determinants,
  Matrices and Invariants} 3rd. edn. (Dover, N.Y. 1960).}
  \ref{MacMillan}{MacMillan, W.D. {\it The Theory of the Potential}
  (McGraw-Hill, N.Y. 1930).}
   \ref{Hobson}{Hobson, E.W. {\it The Theory of Spherical
   and Ellipsoidal Harmonics} (C.U.P., Cambridge. 1931).}
  \ref{Hobson1}{Hobson, E.W. \plms {24}{1892}{55}.}
  \ref{GandY}{Grace, J.H. and Young, A. {\it The Algebra of Invariants}
  (C.U.P., Cambridge, 1903).}
  \ref{FandR}{Fano, U. and Racah, G. {\it Irreducible Tensorial Sets}
  (Academic Press, N.Y. 1959).}
  \ref{TandT}{Thomson, W. and Tait, P.G. {\it Treatise on Natural Philosophy}
   (Clarendon Press, Oxford. 1867).}
  \ref{Brinkman}{Brinkman, H.C. {\it Applications of spinor invariants in
atomic physics}, North Holland, Amsterdam 1956.}
  \ref{Kramers1}{Kramers, H.A. {\it Proc. Roy. Soc. Amst.} {\bf 33} (1930) 953.}
  \ref{DandP2}{Dowker,J.S. and Pettengill,D.F. \jpa{7}{1974}{1527}}
  \ref{Dowk1}{Dowker,J.S. \jpa{}{}{45}.}
  \ref{Dowk2}{Dowker,J.S. \aop{71}{1972}{577}}
  \ref{DandA}{Dowker,J.S. and Apps, J.S. \cqg{15}{1998}{1121}.}
  \ref{Weil}{Weil,A., {\it Elliptic functions according to Eisenstein
  and Kronecker}, Springer, Berlin, 1976.}
  \ref{Ling}{Ling,C-H. {\it SIAM J.Math.Anal.} {\bf5} (1974) 551.}
  \ref{Ling2}{Ling,C-H. {\it J.Math.Anal.Appl.}(1988).}
 \ref{BMO}{Brevik,I., Milton,K.A. and Odintsov, S.D. \aop{302}{2002}{120}.}
 \ref{KandL}{Kutasov,D. and Larsen,F. {\it JHEP} 0101 (2001) 1.}
 \ref{KPS}{Klemm,D., Petkou,A.C. and Siopsis {\it Entropy
 bounds, monoticity properties and scaling in CFT's}. hep-th/0101076.}
 \ref{DandC}{Dowker,J.S. and Critchley,R. \prD{15}{1976}{1484}.}
 \ref{AandD}{Al'taie, M.B. and Dowker, J.S. \prD{18}{1978}{3557}.}
 \ref{Dow1}{Dowker,J.S. \prD{37}{1988}{558}.}
 \ref{Dow30}{Dowker,J.S. \prD{28}{1983}{3013}.}
 \ref{DandK}{Dowker,J.S. and Kennedy,G. \jpa{}{1978}{895}.}
 \ref{Dow2}{Dowker,J.S. \cqg{1}{1984}{359}.}
 \ref{DandKi}{Dowker,J.S. and Kirsten, K. {\it Comm. in Anal. and Geom.
 }{\bf7} (1999) 641.}
 \ref{DandKe}{Dowker,J.S. and Kennedy,G.\jpa{11}{1978}{895}.}
 \ref{Gibbons}{Gibbons,G.W. \pl{60A}{1977}{385}.}
 \ref{Cardy}{Cardy,J.L. \np{366}{1991}{403}.}
 \ref{ChandD}{Chang,P. and Dowker,J.S. \np{395}{1993}{407}.}
 \ref{DandC2}{Dowker,J.S. and Critchley,R. \prD{13}{1976}{224}.}
 \ref{Camporesi}{Camporesi,R. \prp{196}{1990}{1}.}
 \ref{BandM}{Brown,L.S. and Maclay,G.J. \pr{184}{1969}{1272}.}
 \ref{CandD}{Candelas,P. and Dowker,J.S. \prD{19}{1979}{2902}.}
 \ref{Unwin1}{Unwin,S.D. Thesis. University of Manchester. 1979.}
 \ref{Unwin2}{Unwin,S.D. \jpa{13}{1980}{313}.}
 \ref{DandB}{Dowker,J.S. and Banach,R. \jpa{11}{1978}{2255}.}
 \ref{Obhukov}{Obhukov,Yu.N. \pl{109B}{1982}{195}.}
 \ref{Kennedy}{Kennedy,G. \prD{23}{1981}{2884}.}
 \ref{CandT}{Copeland,E. and Toms,D.J. \np {255}{1985}{201}.}
  \ref{CandT2}{Copeland,E. and Toms,D.J. \cqg {3}{1986}{431}.}
 \ref{ELV}{Elizalde,E., Lygren, M. and Vassilevich,
 D.V. \jmp {37}{1996}{3105}.}
 \ref{Malurkar}{Malurkar,S.L. {\it J.Ind.Math.Soc} {\bf16} (1925/26) 130.}
 \ref{Glaisher}{Glaisher,J.W.L. {\it Messenger of Math.} {\bf18}
(1889) 1.} \ref{Anderson}{Anderson,A. \prD{37}{1988}{536}.}
 \ref{CandA}{Cappelli,A. and D'Appollonio,\pl{487B}{2000}{87}.}
 \ref{Wot}{Wotzasek,C. \jpa{23}{1990}{1627}.}
 \ref{RandT}{Ravndal,F. and Tollesen,D. \prD{40}{1989}{4191}.}
 \ref{SandT}{Santos,F.C. and Tort,A.C. \pl{482B}{2000}{323}.}
 \ref{FandO}{Fukushima,K. and Ohta,K. {\it Physica} {\bf A299} (2001) 455.}
 \ref{GandP}{Gibbons,G.W. and Perry,M. \prs{358}{1978}{467}.}
 \ref{Dow4}{Dowker,J.S..}
  \ref{Rad}{Rademacher,H. {\it Topics in analytic number theory,}
Springer-Verlag,  Berlin,1973.}
  \ref{Halphen}{Halphen,G.-H. {\it Trait\'e des Fonctions Elliptiques},
  Vol 1, Gauthier-Villars, Paris, 1886.}
  \ref{CandW}{Cahn,R.S. and Wolf,J.A. {\it Comm.Mat.Helv.} {\bf 51}
  (1976) 1.}
  \ref{Berndt}{Berndt,B.C. \rmjm{7}{1977}{147}.}
  \ref{Hurwitz}{Hurwitz,A. \ma{18}{1881}{528}.}
  \ref{Hurwitz2}{Hurwitz,A. {\it Mathematische Werke} Vol.I. Basel,
  Birkhauser, 1932.}
  \ref{Berndt2}{Berndt,B.C. \jram{303/304}{1978}{332}.}
  \ref{RandA}{Rao,M.B. and Ayyar,M.V. \jims{15}{1923/24}{150}.}
  \ref{Hardy}{Hardy,G.H. \jlms{3}{1928}{238}.}
  \ref{TandM}{Tannery,J. and Molk,J. {\it Fonctions Elliptiques},
   Gauthier-Villars, Paris, 1893--1902.}
  \ref{schwarz}{Schwarz,H.-A. {\it Formeln und
  Lehrs\"atzen zum Gebrauche..},Springer 1893.(The first edition was 1885.)
  The French translation by Henri Pad\'e is {\it Formules et Propositions
  pour L'Emploi...},Gauthier-Villars, Paris, 1894}
  \ref{Hancock}{Hancock,H. {\it Theory of elliptic functions}, Vol I.
   Wiley, New York 1910.}
  \ref{watson}{Watson,G.N. \jlms{3}{1928}{216}.}
  \ref{MandO}{Magnus,W. and Oberhettinger,F. {\it Formeln und S\"atze},
  Springer-Verlag, Berlin 1948.}
  \ref{Klein}{Klein,F. {\it Lectures on the Icosohedron}
  (Methuen, London. 1913).}
  \ref{AandL}{Appell,P. and Lacour,E. {\it Fonctions Elliptiques},
  Gauthier-Villars,
  Paris. 1897.}
  \ref{HandC}{Hurwitz,A. and Courant,C. {\it Allgemeine Funktionentheorie},
  Springer,
  Berlin. 1922.}
  \ref{WandW}{Whittaker,E.T. and Watson,G.N. {\it Modern analysis},
  Cambridge. 1927.}
  \ref{SandC}{Selberg,A. and Chowla,S. \jram{227}{1967}{86}. }
  \ref{zucker}{Zucker,I.J. {\it Math.Proc.Camb.Phil.Soc} {\bf 82 }(1977)
  111.}
  \ref{glasser}{Glasser,M.L. {\it Maths.of Comp.} {\bf 25} (1971) 533.}
  \ref{GandW}{Glasser, M.L. and Wood,V.E. {\it Maths of Comp.} {\bf 25}
  (1971)
  535.}
  \ref{greenhill}{Greenhill,A,G. {\it The Applications of Elliptic
  Functions}, MacMillan. London, 1892.}
  \ref{Weierstrass}{Weierstrass,K. {\it J.f.Mathematik (Crelle)}
{\bf 52} (1856) 346.}
  \ref{Weierstrass2}{Weierstrass,K. {\it Mathematische Werke} Vol.I,p.1,
  Mayer u. M\"uller, Berlin, 1894.}
  \ref{Fricke}{Fricke,R. {\it Die Elliptische Funktionen und Ihre Anwendungen},
    Teubner, Leipzig. 1915, 1922.}
  \ref{Konig}{K\"onigsberger,L. {\it Vorlesungen \"uber die Theorie der
 Elliptischen Funktionen},  \break Teubner, Leipzig, 1874.}
  \ref{Milne}{Milne,S.C. {\it The Ramanujan Journal} {\bf 6} (2002) 7-149.}
  \ref{Schlomilch}{Schl\"omilch,O. {\it Ber. Verh. K. Sachs. Gesell. Wiss.
  Leipzig}  {\bf 29} (1877) 101-105; {\it Compendium der h\"oheren
  Analysis}, Bd.II, 3rd Edn, Vieweg, Brunswick, 1878.}
  \ref{BandB}{Briot,C. and Bouquet,C. {\it Th\`eorie des Fonctions
  Elliptiques}, Gauthier-Villars, Paris, 1875.}
  \ref{Dumont}{Dumont,D. \aim {41}{1981}{1}.}
  \ref{Andre}{Andr\'e,D. {\it Ann.\'Ecole Normale Superior} {\bf 6} (1877)
  265;
  {\it J.Math.Pures et Appl.} {\bf 5} (1878) 31.}
  \ref{Raman}{Ramanujan,S. {\it Trans.Camb.Phil.Soc.} {\bf 22} (1916) 159;
 {\it Collected Papers}, Cambridge, 1927}
  \ref{Weber}{Weber,H.M. {\it Lehrbuch der Algebra} Bd.III, Vieweg,
  Brunswick 190  3.}
  \ref{Weber2}{Weber,H.M. {\it Elliptische Funktionen und algebraische
  Zahlen},
  Vieweg, Brunswick 1891.}
  \ref{ZandR}{Zucker,I.J. and Robertson,M.M.
  {\it Math.Proc.Camb.Phil.Soc} {\bf 95 }(1984) 5.}
  \ref{JandZ1}{Joyce,G.S. and Zucker,I.J.
  {\it Math.Proc.Camb.Phil.Soc} {\bf 109 }(1991) 257.}
  \ref{JandZ2}{Zucker,I.J. and Joyce.G.S.
  {\it Math.Proc.Camb.Phil.Soc} {\bf 131 }(2001) 309.}
  \ref{zucker2}{Zucker,I.J. {\it SIAM J.Math.Anal.} {\bf 10} (1979) 192,}
  \ref{BandZ}{Borwein,J.M. and Zucker,I.J. {\it IMA J.Math.Anal.} {\bf 12}
  (1992) 519.}
  \ref{Cox}{Cox,D.A. {\it Primes of the form $x^2+n\,y^2$}, Wiley,
  New York, 1989.}
  \ref{BandCh}{Berndt,B.C. and Chan,H.H. {\it Mathematika} {\bf42} (1995)
  278.}
  \ref{EandT}{Elizalde,R. and Tort.hep-th/}
  \ref{KandS}{Kiyek,K. and Schmidt,H. {\it Arch.Math.} {\bf 18} (1967) 438.}
  \ref{Oshima}{Oshima,K. \prD{46}{1992}{4765}.}
  \ref{greenhill2}{Greenhill,A.G. \plms{19} {1888} {301}.}
  \ref{Russell}{Russell,R. \plms{19} {1888} {91}.}
  \ref{BandB}{Borwein,J.M. and Borwein,P.B. {\it Pi and the AGM}, Wiley,
  New York, 1998.}
  \ref{Resnikoff}{Resnikoff,H.L. \tams{124}{1966}{334}.}
  \ref{vandp}{Van der Pol, B. {\it Indag.Math.} {\bf18} (1951) 261,272.}
  \ref{Rankin}{Rankin,R.A. {\it Modular forms} C.U.P. Cambridge}
  \ref{Rankin2}{Rankin,R.A. {\it Proc. Roy.Soc. Edin.} {\bf76 A} (1976) 107.}
  \ref{Skoruppa}{Skoruppa,N-P. {\it J.of Number Th.} {\bf43} (1993) 68 .}
  \ref{Down}{Dowker.J.S. {\it Nucl.Phys.}B (Proc.Suppl) ({\bf 104})(2002)153;
  also Dowker,J.S. hep-th/ 0007129.}
  \ref{Eichler}{Eichler,M. \mz {67}{1957}{267}.}
  \ref{Zagier}{Zagier,D. \invm{104}{1991}{449}.}
  \ref{Lang}{Lang,S. {\it Modular Forms}, Springer, Berlin, 1976.}
  \ref{Kosh}{Koshliakov,N.S. {\it Mess.of Math.} {\bf 58} (1928) 1.}
  \ref{BandH}{Bodendiek, R. and Halbritter,U. \amsh{38}{1972}{147}.}
  \ref{Smart}{Smart,L.R., \pgma{14}{1973}{1}.}
  \ref{Grosswald}{Grosswald,E. {\it Acta. Arith.} {\bf 21} (1972) 25.}
  \ref{Kata}{Katayama,K. {\it Acta Arith.} {\bf 22} (1973) 149.}
  \ref{Ogg}{Ogg,A. {\it Modular forms and Dirichlet series} (Benjamin,
  New York,
   1969).}
  \ref{Bol}{Bol,G. \amsh{16}{1949}{1}.}
  \ref{Epstein}{Epstein,P. \ma{56}{1903}{615}.}
  \ref{Petersson}{Petersson.}
  \ref{Serre}{Serre,J-P. {\it A Course in Arithmetic}, Springer,
  New York, 1973.}
  \ref{Schoenberg}{Schoenberg,B., {\it Elliptic Modular Functions},
  Springer, Berlin, 1974.}
  \ref{Apostol}{Apostol,T.M. \dmj {17}{1950}{147}.}
  \ref{Ogg2}{Ogg,A. {\it Lecture Notes in Math.} {\bf 320} (1973) 1.}
  \ref{Knopp}{Knopp,M.I. \dmj {45}{1978}{47}.}
  \ref{Knopp2}{Knopp,M.I. \invm {}{1994}{361}.}
  \ref{LandZ}{Lewis,J. and Zagier,D. \aom{153}{2001}{191}.}
  \ref{DandK1}{Dowker,J.S. and Kirsten,K. {\it Elliptic functions and
  temperature inversion symmetry on spheres} hep-th/.}
  \ref{HandK}{Husseini and Knopp.}
  \ref{Kober}{Kober,H. \mz{39}{1934-5}{609}.}
  \ref{HandL}{Hardy,G.H. and Littlewood, \am{41}{1917}{119}.}
  \ref{Watson}{Watson,G.N. \qjm{2}{1931}{300}.}
  \ref{SandC2}{Chowla,S. and Selberg,A. {\it Proc.Nat.Acad.} {\bf 35}
  (1949) 371.}
  \ref{Landau}{Landau, E. {\it Lehre von der Verteilung der Primzahlen},
  (Teubner, Leipzig, 1909).}
  \ref{Berndt4}{Berndt,B.C. \tams {146}{1969}{323}.}
  \ref{Berndt3}{Berndt,B.C. \tams {}{}{}.}
  \ref{Bochner}{Bochner,S. \aom{53}{1951}{332}.}
  \ref{Weil2}{Weil,A.\ma{168}{1967}{}.}
  \ref{CandN}{Chandrasekharan,K. and Narasimhan,R. \aom{74}{1961}{1}.}
  \ref{Rankin3}{Rankin,R.A. {} {} ().}
  \ref{Berndt6}{Berndt,B.C. {\it Trans.Edin.Math.Soc}.}
  \ref{Elizalde}{Elizalde,E. {\it Ten Physical Applications of Spectral
  Zeta Function Theory}, \break (Springer, Berlin, 1995).}
  \ref{Allen}{Allen,B., Folacci,A. and Gibbons,G.W. \pl{189}{1987}{304}.}
  \ref{Krazer}{Krazer}
  \ref{Elizalde3}{Elizalde,E. {\it J.Comp.and Appl. Math.} {\bf 118}
  (2000) 125.}
  \ref{Elizalde2}{Elizalde,E., Odintsov.S.D, Romeo, A. and Bytsenko,
  A.A and
  Zerbini,S.
  {\it Zeta function regularisation}, (World Scientific, Singapore,
  1994).}
  \ref{Eisenstein}{Eisenstein}
  \ref{Hecke}{Hecke,E. \ma{112}{1936}{664}.}
  \ref{Hecke2}{Hecke,E. \ma{112}{1918}{398}.}
  \ref{Terras}{Terras,A. {\it Harmonic analysis on Symmetric Spaces} (Springer,
  New York, 1985).}
  \ref{BandG}{Bateman,P.T. and Grosswald,E. {\it Acta Arith.} {\bf 9}
  (1964) 365.}
  \ref{Deuring}{Deuring,M. \aom{38}{1937}{585}.}
  \ref{Mordell}{Mordell,J. \prs{}{}{}.}
  \ref{GandZ}{Glasser,M.L. and Zucker, {}.}
  \ref{Landau2}{Landau,E. \jram{}{1903}{64}.}
  \ref{Kirsten1}{Kirsten,K. \jmp{35}{1994}{459}.}
  \ref{Sommer}{Sommer,J. {\it Vorlesungen \"uber Zahlentheorie}
  (1907,Teubner,Leipzig).
  French edition 1913 .}
  \ref{Reid}{Reid,L.W. {\it Theory of Algebraic Numbers},
  (1910,MacMillan,New York).}
  \ref{Milnor}{Milnor, J. {\it Is the Universe simply--connected?},
  IAS, Princeton, 1978.}
  \ref{Milnor2}{Milnor, J. \ajm{79}{1957}{623}.}
  \ref{Opechowski}{Opechowski,W. {\it Physica} {\bf 7} (1940) 552.}
  \ref{Bethe}{Bethe, H.A. \zfp{3}{1929}{133}.}
  \ref{LandL}{Landau, L.D. and Lishitz, E.M. {\it Quantum
  Mechanics} (Pergamon Press, London, 1958).}
  \ref{GPR}{Gibbons, G.W., Pope, C. and R\"omer, H., \np{157}{1979}{377}.}
  \ref{Jadhav}{Jadhav,S.P. PhD Thesis, University of Manchester 1990.}
  \ref{DandJ}{Dowker,J.S. and Jadhav, S. \prD{39}{1989}{1196}.}
  \ref{CandM}{Coxeter, H.S.M. and Moser, W.O.J. {\it Generators and
  relations of finite groups} (Springer. Berlin. 1957).}
  \ref{Coxeter2}{Coxeter, H.S.M. {\it Regular Complex Polytopes},
   (Cambridge University Press, \break Cambridge, 1975).}
  \ref{Coxeter}{Coxeter, H.S.M. {\it Regular Polytopes}.}
  \ref{Stiefel}{Stiefel, E., J.Research NBS {\bf 48} (1952) 424.}
  \ref{BandS}{Brink, D.M. and Satchler, G.R. {\it Angular momentum theory}.
  (Clarendon Press, Oxford. 1962.).}
  %\ref{Racah1}
  \ref{Rose}{Rose}
  \ref{Schwinger}{Schwinger, J. {\it On Angular Momentum}
  in {\it Quantum Theory of Angular Momentum} edited by
  Biedenharn,L.C. and van Dam, H. (Academic Press, N.Y. 1965).}
  \ref{Bromwich}{Bromwich, T.J.I'A. {\it Infinite Series},
  (Macmillan, 1947).}
  \ref{Ray}{Ray,D.B. \aim{4}{1970}{109}.}
  \ref{Ikeda}{Ikeda,A. {\it Kodai Math.J.} {\bf 18} (1995) 57.}
  \ref{Kennedy}{Kennedy,G. \prD{23}{1981}{2884}.}
  \ref{Ellis}{Ellis,G.F.R. {\it General Relativity} {\bf2} (1971) 7.}
  \ref{Dow8}{Dowker,J.S. \cqg{20}{2003}{L105}.}
  \ref{IandY}{Ikeda, A and Yamamoto, Y. \ojm {16}{1979}{447}.}
  \ref{BandI}{Bander,M. and Itzykson,C. \rmp{18}{1966}{2}.}
  \ref{Schulman}{Schulman, L.S. \pr{176}{1968}{1558}.}
  \ref{Bar1}{B\"ar,C. {\it Arch.d.Math.}{\bf 59} (1992) 65.}
  \ref{Bar2}{B\"ar,C. {\it Geom. and Func. Anal.} {\bf 6} (1996) 899.}
  \ref{Vilenkin}{Vilenkin, N.J. {\it Special functions},
  (Am.Math.Soc., Providence, 1968).}
  \ref{Talman}{Talman, J.D. {\it Special functions} (Benjamin,N.Y.,1968).}
  \ref{Miller}{Miller, W. {\it Symmetry groups and their applications}
  (Wiley, N.Y., 1972).}
  \ref{Dow3}{Dowker,J.S. \cmp{162}{1994}{633}.}
  \ref{Cheeger}{Cheeger, J. \jdg {18}{1983}{575}.}
  \ref{Cheeger2}{Cheeger, J. \aom {109}{1979}{259}.}
  \ref{Dow6}{Dowker,J.S. \jmp{30}{1989}{770}.}
  \ref{Dow20}{Dowker,J.S. \jmp{35}{1994}{6076}.}
  \ref{Dowjmp}{Dowker,J.S. \jmp{35}{1994}{4989}.}
  \ref{Dow21}{Dowker,J.S. {\it Heat kernels and polytopes} in {\it
   Heat Kernel Techniques and Quantum Gravity}, ed. by S.A.Fulling,
   Discourses in Mathematics and its Applications, No.4, Dept.
   Maths., Texas A\&M University, College Station, Texas, 1995.}
  \ref{Dow9}{Dowker,J.S. \jmp{42}{2001}{1501}.}
  \ref{Dow7}{Dowker,J.S. \jpa{25}{1992}{2641}.}
  \ref{Warner}{Warner.N.P. \prs{383}{1982}{379}.}
  \ref{Wolf}{Wolf, J.A. {\it Spaces of constant curvature},
  (McGraw--Hill,N.Y., 1967).}
  \ref{Meyer}{Meyer,B. \cjm{6}{1954}{135}.}
  \ref{BandB}{B\'erard,P. and Besson,G. {\it Ann. Inst. Four.} {\bf 30}
  (1980) 237.}
  \ref{PandM}{P\'{o}lya,G. and Meyer,B. \cras{228}{1948}{28}.}
  \ref{Springer}{Springer, T.A. Lecture Notes in Math. vol 585 (Springer,
  Berlin,1977).}
  \ref{SeandT}{Threlfall, H. and Seifert, W. \ma{104}{1930}{1}.}
  \ref{Hopf}{Hopf,H. \ma{95}{1925}{313}. }
  \ref{Dow}{Dowker,J.S. \jpa{5}{1972}{936}.}
  \ref{LLL}{Lehoucq,R., Lachi\'eze-Rey,M. and Luminet, J.--P. {\it
  Astron.Astrophys.} {\bf 313} (1996) 339.}
  \ref{LaandL}{Lachi\'eze-Rey,M. and Luminet, J.--P.
  \prp{254}{1995}{135}.}
  \ref{Schwarzschild}{Schwarzschild, K., {\it Vierteljahrschrift der
  Ast.Ges.} {\bf 35} (1900) 337.}
  \ref{Starkman}{Starkman,G.D. \cqg{15}{1998}{2529}.}
  \ref{LWUGL}{Lehoucq,R., Weeks,J.R., Uzan,J.P., Gausman, E. and
  Luminet, J.--P. \cqg{19}{2002}{4683}.}
  \ref{Dow10}{Dowker,J.S. \prD{28}{1983}{3013}.}
  \ref{BandD}{Banach, R. and Dowker, J.S. \jpa{12}{1979}{2527}.}
  \ref{Jadhav2}{Jadhav,S. \prD{43}{1991}{2656}.}
  \ref{Gilkey}{Gilkey,P.B. {\it Invariance theory,the heat equation and
  the Atiyah--Singer Index theorem} (CRC Press, Boca Raton, 1994).}
  \ref{BandY}{Berndt,B.C. and Yeap,B.P. {\it Adv. Appl. Math.}
  {\bf29} (2002) 358.}
  \ref{HandR}{Hanson,A.J. and R\"omer,H. \pl{80B}{1978}{58}.}
  \ref{Hill}{Hill,M.J.M. {\it Trans.Camb.Phil.Soc.} {\bf 13} (1883) 36.}
  \ref{Cayley}{Cayley,A. {\it Quart.Math.J.} {\bf 7} (1866) 304.}
  \ref{Seade}{Seade,J.A. {\it Anal.Inst.Mat.Univ.Nac.Aut\'on
  M\'exico} {\bf 21} (1981) 129.}
  \ref{CM}{Cisneros--Molina,J.L. {\it Geom.Dedicata} {\bf84} (2001)
  \ref{Goette1}{Goette,S. \jram {526} {2000} 181.}
  207.}
  \ref{NandO}{Nash,C. and O'Connor,D--J, \jmp {36}{1995}{1462}.}
  \ref{Dows}{Dowker,J.S. \aop{71}{1972}{577}; Dowker,J.S. and Pettengill,D.F.
  \jpa{7}{1974}{1527}; J.S.Dowker in {\it Quantum Gravity}, edited by
  S. C. Christensen (Hilger,Bristol,1984)}
  \ref{Jadhav2}{Jadhav,S.P. \prD{43}{1991}{2656}.}
  \ref{Dow11}{Dowker,J.S. \cqg{21}{2004}4247.}
  \ref{Dow12}{Dowker,J.S. \cqg{21}{2004}4977.}
  \ref{Dow13}{Dowker,J.S. \jpa{38}{2005}1049.}
  \ref{Zagier}{Zagier,D. \ma{202}{1973}{149}}
  \ref{RandG}{Rademacher, H. and Grosswald,E. {\it Dedekind Sums},
  (Carus, MAA, 1972).}
  \ref{Berndt7}{Berndt,B, \aim{23}{1977}{285}.}
  \ref{HKMM}{Harvey,J.A., Kutasov,D., Martinec,E.J. and Moore,G.
  {\it Localised Tachyons and RG Flows}, hep-th/0111154.}
  \ref{Beck}{Beck,M., {\it Dedekind Cotangent Sums}, {\it Acta Arithmetica}
  {\bf 109} (2003) 109-139 ; math.NT/0112077.}
  \ref{McInnes}{McInnes,B. {\it APS instability and the topology of the brane
  world}, hep-th/0401035.}
  \ref{BHS}{Brevik,I, Herikstad,R. and Skriudalen,S. {\it Entropy Bound for the
  TM Electromagnetic Field in the Half Einstein Universe}; hep-th/0508123.}
  \ref{BandO}{Brevik,I. and Owe,C.  \prD{55}{4689}{1997}.}
  \ref{Kenn}{Kennedy,G. Thesis. University of Manchester 1978.}
  \ref{KandU}{Kennedy,G. and Unwin S. \jpa{12}{L253}{1980}.}
  \ref{BandO1}{Bayin,S.S.and Ozcan,M.
  \prD{48}{2806}{1993}; \prD{49}{5313}{1994}.}
  \ref{Chang}{Chang, P., {\it Quantum Field Theory on Regular Polytopes}.
   Thesis. University of Manchester, 1993.}
  \ref{Barnesa}{Barnes,E.W. {\it Trans. Camb. Phil. Soc.} {\bf 19} (1903) 374.}
  \ref{Barnesb}{Barnes,E.W. {\it Trans. Camb. Phil. Soc.}
  {\bf 19} (1903) 426.}
  \ref{Stanley1}{Stanley,R.P. \joa {49Hilf}{1977}{134}.}
  \ref{Stanley}{Stanley,R.P. \bams {1}{1979}{475}.}
  \ref{Hurley}{Hurley,A.C. \pcps {47}{1951}{51}.}
  \ref{IandK}{Iwasaki,I. and Katase,K. {\it Proc.Japan Acad. Ser} {\bf A55}
  (1979) 141.}
  \ref{IandT}{Ikeda,A. and Taniguchi,Y. {\it Osaka J. Math.} {\bf 15} (1978)
  515.}
  \ref{GandM}{Gallot,S. and Meyer,D. \jmpa{54}{1975}{259}.}
  \ref{Flatto}{Flatto,L. {\it Enseign. Math.} {\bf 24} (1978) 237.}
  \ref{OandT}{Orlik,P and Terao,H. {\it Arrangements of Hyperplanes},
  Grundlehren der Math. Wiss. {\bf 300}, (Springer--Verlag, 1992).}
  \ref{Shepler}{Shepler,A.V. \joa{220}{1999}{314}.}
  \ref{SandT}{Solomon,L. and Terao,H. \cmh {73}{1998}{237}.}
  \ref{Vass}{Vassilevich, D.V. \plb {348}{1995}39.}
  \ref{Vass2}{Vassilevich, D.V. \jmp {36}{1995}3174.}
  \ref{CandH}{Camporesi,R. and Higuchi,A. {\it J.Geom. and Physics}
  {\bf 15} (1994) 57.}
  \ref{Solomon2}{Solomon,L. \tams{113}{1964}{274}.}
  \ref{Solomon}{Solomon,L. {\it Nagoya Math. J.} {\bf 22} (1963) 57.}
  \ref{Obukhov}{Obukhov,Yu.N. \pl{109B}{1982}{195}.}
  \ref{BGH}{Bernasconi,F., Graf,G.M. and Hasler,D. {\it The heat kernel
  expansion for the electromagnetic field in a cavity}; math-ph/0302035.}
  \ref{Baltes}{Baltes,H.P. \prA {6}{1972}{2252}.}
  \ref{BaandH}{Baltes.H.P and Hilf,E.R. {\it Spectra of Finite Systems}
  (Bibliographisches Institut, Mannheim, 1976).}
  \ref{Ray}{Ray,D.B. \aim{4}{1970}{109}.}
  \ref{Hirzebruch}{Hirzebruch,F. {\it Topological methods in algebraic
  geometry} (Springer-- Verlag,\break  Berlin, 1978). }
  \ref{BBG}{Bla\v{z}i\'c,N., Bokan,N. and Gilkey, P.B. {\it Ind.J.Pure and
  Appl.Math.} {\bf 23} (1992) 103.}
  \ref{WandWi}{Weck,N. and Witsch,K.J. {\it Math.Meth.Appl.Sci.} {\bf 17}
  (1994) 1017.}
  \ref{Norlund}{N\"orlund,N.E. \am{43}{1922}{121}.}
   \ref{Norlund1}{N\"orlund,N.E. {\it Differenzenrechnung} (Springer--Verlag, 1924, Berlin.)}
  \ref{Duff}{Duff,G.F.D. \aom{56}{1952}{115}.}
  \ref{DandS}{Duff,G.F.D. and Spencer,D.C. \aom{45}{1951}{128}.}
  \ref{BGM}{Berger, M., Gauduchon, P. and Mazet, E. {\it Lect.Notes.Math.}
  {\bf 194} (1971) 1. }
  \ref{Patodi}{Patodi,V.K. \jdg{5}{1971}{233}.}
  \ref{GandS}{G\"unther,P. and Schimming,R. \jdg{12}{1977}{599}.}
  \ref{MandS}{McKean,H.P. and Singer,I.M. \jdg{1}{1967}{43}.}
  \ref{Conner}{Conner,P.E. {\it Mem.Am.Math.Soc.} {\bf 20} (1956).}
  \ref{Gilkey2}{Gilkey,P.B. \aim {15}{1975}{334}.}
  \ref{MandP}{Moss,I.G. and Poletti,S.J. \plb{333}{1994}{326}.}
  \ref{BKD}{Bordag,M., Kirsten,K. and Dowker,J.S. \cmp{182}{1996}{371}.}
  \ref{RandO}{Rubin,M.A. and Ordonez,C. \jmp{25}{1984}{2888}.}
  \ref{BaandD}{Balian,R. and Duplantier,B. \aop {112}{1978}{165}.}
  \ref{Kennedy2}{Kennedy,G. \aop{138}{1982}{353}.}
  \ref{DandKi2}{Dowker,J.S. and Kirsten, K. {\it Analysis and Appl.}
 {\bf 3} (2005) 45.}
  \ref{Dow40}{Dowker,J.S. \cqg{23}{2006}{1}.}
  \ref{BandHe}{Br\"uning,J. and Heintze,E. {\it Duke Math.J.} {\bf 51} (1984)
   959.}
  \ref{Dowl}{Dowker,J.S. {\it Functional determinants on M\"obius corners};
    Proceedings, `Quantum field theory under
    the influence of external conditions', 111-121,Leipzig 1995.}
  \ref{Dowqg}{Dowker,J.S. in {\it Quantum Gravity}, edited by
  S. C. Christensen (Hilger, Bristol, 1984).}
  \ref{Dowit}{Dowker,J.S. \jpa{11}{1978}{347}.}
  \ref{Kane}{Kane,R. {\it Reflection Groups and Invariant Theory} (Springer,
  New York, 2001).}
  \ref{Sturmfels}{Sturmfels,B. {\it Algorithms in Invariant Theory}
  (Springer, Vienna, 1993).}
  \ref{Bourbaki}{Bourbaki,N. {\it Groupes et Alg\`ebres de Lie}  Chap.III, IV
  (Hermann, Paris, 1968).}
  \ref{SandTy}{Schwarz,A.S. and Tyupkin, Yu.S. \np{242}{1984}{436}.}
  \ref{Reuter}{Reuter,M. \prD{37}{1988}{1456}.}
  \ref{EGH}{Eguchi,T. Gilkey,P.B. and Hanson,A.J. \prp{66}{1980}{213}.}
  \ref{DandCh}{Dowker,J.S. and Chang,Peter, \prD{46}{1992}{3458}.}
  \ref{APS}{Atiyah M., Patodi and Singer,I.\mpcps{77}{1975}{43}.}
  \ref{Donnelly}{Donnelly.H. {\it Indiana U. Math.J.} {\bf 27} (1978) 889.}
  \ref{Katase}{Katase,K. {\it Proc.Jap.Acad.} {\bf 57} (1981) 233.}
  \ref{Gilkey3}{Gilkey,P.B.\invm{76}{1984}{309}.}
  \ref{Degeratu}{Degeratu.A. {\it Eta--Invariants and Molien Series for
  Unimodular Groups}, Thesis MIT, 2001.}
  \ref{Seeley}{Seeley,R. \ijmp {A\bf18}{2003}{2197}.}
  \ref{Seeley2}{Seeley,R. .}
  \ref{melrose}{Melrose}
  \ref{DandW}{Douglas,R.G. and Wojciekowski,K.P. \cmp{142}{1991}{139}.}
  \ref{Dai}{Dai,X. \tams{354}{2001}{107}.}
\end{putreferences}

\bye